\documentclass[12pt]{article}
\usepackage{amsmath,amssymb}
\usepackage[enableskew]{youngtab}
\usepackage{xypic}
\usepackage{eepic}
\usepackage{longtable}

\newtheorem{mytheorem}{Theorem}
\newtheorem{myproposition}[mytheorem]{Proposition}
\newtheorem{mylemma}[mytheorem]{Lemma}
\newtheorem{mycorollary}[mytheorem]{Corollary}
\newenvironment{mydefinition}{\vspace{2mm}\noindent\textbf{Definition}}{\vspace{2mm}}
\newenvironment{myremark}{\vspace{2mm}\noindent\small\textbf{Remark}}{\normalsize\vspace{2mm}}
\newenvironment{proof}{\textit{Proof.}}{\mbox{}\hfill$\square$\vspace{2mm}}

\newcommand{\Fbox}{\mbox}
\newcommand{\vol}{\operatorname{vol}}
\newcommand{\rank}{\operatorname{rk}}
\newcommand{\Aut}{\operatorname{Aut}}
\newcommand{\Dih}{\operatorname{Dih}}
\newcommand{\Hasse}{\operatorname{Hasse}}
\newcommand{\ab}{\mathfrak a}
\newcommand{\borel}{\mathfrak b}
\newcommand{\lie}{\mathfrak g}
\newcommand{\cartan}{\mathfrak h}
\newcommand{\anc}{{\not\perp\theta}}
\newcommand{\sfA}{\mathsf A}
\newcommand{\sfB}{\mathsf B}
\newcommand{\sfC}{\mathsf C}
\newcommand{\sfD}{\mathsf D}
\newcommand{\sfE}{\mathsf E}
\newcommand{\sfF}{\mathsf F}
\newcommand{\sfG}{\mathsf G}

\newcommand{\sfX}{\mathsf X}
\newcommand{\RR}{\mathbb R}
\newcommand{\ZZ}{\mathbb Z}

\newcounter{linearcounter}
\newcommand{\linear}[7]{\setlength{\unitlength}{4mm}
\begin{picture}(0,7)(.1,0)
\path(0,0)(0,\value{linearcounter})
\addtocounter{linearcounter}{1}
\multiput(0,0)(0,1){\value{linearcounter}}{\circle*{.2}}
\put(0,0.5){\makebox(0,0)[r]{\footnotesize#1}}
\put(0,1.5){\makebox(0,0)[r]{\footnotesize#2}}
\put(0,2.5){\makebox(0,0)[r]{\footnotesize#3}}
\put(0,3.5){\makebox(0,0)[r]{\footnotesize#4}}
\put(0,4.5){\makebox(0,0)[r]{\footnotesize#5}}
\put(0,5.5){\makebox(0,0)[r]{\footnotesize#6}}
\put(0,6.5){\makebox(0,0)[r]{\footnotesize#7}}
\end{picture}}

\newcommand{\Aloz}{\setlength{\unitlength}{4mm}
\begin{picture}(0,4)(.1,0)
\path(0,0)(0,1)(1,2)(0,3)\path(0,1)(-1,2)(0,3)(0,4)
\multiput(0,0)(0,1){2}{\circle*{.2}}\multiput(0,3)(0,1){2}{\circle*{.2}}
\multiput(-1,2)(2,0){2}{\circle*{.2}}
\put(0,0.5){\makebox(0,0)[r]{\footnotesize$s_0$}}
\put(-.5,1.5){\makebox(0,0)[tr]{\footnotesize$s_1$}}
\put(-.5,2.5){\makebox(0,0)[br]{\footnotesize$s_5$}}
\put(0,3.5){\makebox(0,0)[r]{\footnotesize$s_0$}}
\end{picture}}

\newcommand{\Cloz}{\setlength{\unitlength}{4mm}
\begin{picture}(0,6)(.1,0)
\path(0,0)(0,2)(1,3)(0,4)\path(0,2)(-1,3)(0,4)(0,6)
\multiput(0,0)(0,1){3}{\circle*{.2}}\multiput(0,4)(0,1){3}{\circle*{.2}}
\multiput(-1,3)(2,0){2}{\circle*{.2}}
\put(0,0.5){\makebox(0,0)[r]{\footnotesize$s_0$}}
\put(0,1.5){\makebox(0,0)[r]{\footnotesize$s_1$}}
\put(-.5,2.5){\makebox(0,0)[tr]{\footnotesize$s_0$}}
\put(-.5,3.5){\makebox(0,0)[br]{\footnotesize$s_2$}}
\put(0,4.5){\makebox(0,0)[r]{\footnotesize$s_1$}}
\put(0,5.5){\makebox(0,0)[r]{\footnotesize$s_0$}}
\end{picture}}

\newcommand{\Bloz}{\setlength{\unitlength}{4mm}
\begin{picture}(0,6)(.1,0)
\path(0,0)(0,2)(1,3)(0,4)\path(0,2)(-1,3)(0,4)(0,6)
\multiput(0,0)(0,1){3}{\circle*{.2}}\multiput(0,4)(0,1){3}{\circle*{.2}}
\multiput(-1,3)(2,0){2}{\circle*{.2}}
\put(0,0.5){\makebox(0,0)[r]{\footnotesize$s_0$}}
\put(0,1.5){\makebox(0,0)[r]{\footnotesize$s_2$}}
\put(-.5,2.5){\makebox(0,0)[tr]{\footnotesize$s_1$}}
\put(-.5,3.5){\makebox(0,0)[br]{\footnotesize$s_3$}}
\put(0,4.5){\makebox(0,0)[r]{\footnotesize$s_2$}}
\put(0,5.5){\makebox(0,0)[r]{\footnotesize$s_0$}}
\end{picture}}

\newcommand{\Eloz}{\setlength{\unitlength}{4mm}
\begin{picture}(0,10)(.1,0)
\path(0,0)(0,4)(1,5)(0,6)\path(0,4)(-1,5)(0,6)(0,10)
\multiput(0,0)(0,1){5}{\circle*{.2}}\multiput(0,6)(0,1){5}{\circle*{.2}}
\multiput(-1,5)(2,0){2}{\circle*{.2}}
\put(0,0.5){\makebox(0,0)[r]{\footnotesize$s_0$}}
\put(0,1.5){\makebox(0,0)[r]{\footnotesize$s_1$}}
\put(0,2.5){\makebox(0,0)[r]{\footnotesize$s_2$}}
\put(0,3.5){\makebox(0,0)[r]{\footnotesize$s_3$}}
\put(-.5,4.5){\makebox(0,0)[tr]{\footnotesize$s_4$}}
\put(-.5,5.5){\makebox(0,0)[br]{\footnotesize$s_5$}}
\put(0,6.5){\makebox(0,0)[r]{\footnotesize$s_3$}}
\put(0,7.5){\makebox(0,0)[r]{\footnotesize$s_2$}}
\put(0,8.5){\makebox(0,0)[r]{\footnotesize$s_1$}}
\put(0,9.5){\makebox(0,0)[r]{\footnotesize$s_0$}}
\end{picture}}

\newcommand{\Clong}{\setlength{\unitlength}{4mm}
\begin{picture}(6,10)(-2.9,0)
\path(0,0)(0,2)(2,4)(-1,7)\path(1,3)(-2,6)(0,8)(0,10)
\path(0,2)(-1,3)(1,5)\path(-1,5)(1,7)(0,8)
\multiput(0,0)(0,1){2}{\circle*{.2}}\multiput(0,9)(0,1){2}{\circle*{.2}}
\multiput(0,2)(1,1){3}{\circle*{.2}}\multiput(-1,3)(1,1){3}{\circle*{.2}}
\multiput(-1,5)(1,1){3}{\circle*{.2}}\multiput(-2,6)(1,1){3}{\circle*{.2}}
\put(0,.5){\makebox(0,0)[r]{\footnotesize$s_0$}}
\put(0,1.5){\makebox(0,0)[r]{\footnotesize$s_1$}}
\put(-.5,2.5){\makebox(0,0)[tr]{\footnotesize$s_0$}}
\put(-.5,3.5){\makebox(0,0)[br]{\footnotesize$s_2$}}
\put(-.5,4.5){\makebox(0,0)[tr]{\footnotesize$s_1$}}
\put(-1.5,5.5){\makebox(0,0)[tr]{\footnotesize$s_0$}}
\put(-1.5,6.5){\makebox(0,0)[br]{\footnotesize$s_3$}}
\put(-.5,7.5){\makebox(0,0)[br]{\footnotesize$s_2$}}
\put(0,8.5){\makebox(0,0)[r]{\footnotesize$s_1$}}
\put(0,9.5){\makebox(0,0)[r]{\footnotesize$s_0$}}
\end{picture}}

\newcommand{\AoneAff}[2]{\begin{picture}(0,2.6)(0,-.8)
\put(0,0){\circle*{0.32}}
\put(0,-1){\circle*{0.32}} \put(0,0){\line(0,-1){1}}
\put(0,-1){\circle{0.6}}
\put(0.133,-.5){\makebox(0,0)[l]{\tiny$\infty$}}
\put(0,.333){\makebox(0,0)[b]{\footnotesize#1}}
\put(0,-1.333){\makebox(0,0)[t]{\footnotesize#2}} 
\end{picture}}

\newcommand{\AtwoAff}[3]{\begin{picture}(1,2.6)(0,-.8)
\multiput(0,0)(1,0){2}{\circle*{0.32}}
\put(0,0){\line(1,0){1}}
\put(0.5,-1){\circle*{0.32}} \put(0,0){\line(1,-2){0.5}} \put(0.5,-1){\line(1,2){0.5}}
\put(0.5,-1){\circle{0.6}}
\put(0,.333){\makebox(0,0)[b]{\footnotesize#1}}
\put(1,.333){\makebox(0,0)[b]{\footnotesize#2}}
\put(0.5,-1.333){\makebox(0,0)[t]{\footnotesize#3}} 
\end{picture}}

\newcommand{\AlAffev}[5]{\begin{picture}(5,2.6)(0,-.8)
\multiput(0,0)(2,0){2}{\circle*{0.32}} \multiput(3,0)(2,0){2}{\circle*{0.32}}
\put(0,0){\line(1,0){0.5}} \put(1.5,0){\line(1,0){2}} \put(4.5,0){\line(1,0){0.5}}
\put(2.5,-1){\circle*{0.32}} \put(0,0){\line(5,-2){2.5}} \put(2.5,-1){\line(5,2){2.5}}
\put(2.5,-1){\circle{0.6}}
\put(1,0){\makebox(0,0){\tiny$\cdots$}} \put(4,0){\makebox(0,0){\tiny$\cdots$}}
\put(0,.333){\makebox(0,0)[b]{\footnotesize#1}}
\put(2,.333){\makebox(0,0)[b]{\footnotesize#2}}
\put(3,.333){\makebox(0,0)[b]{\footnotesize#3}}
\put(5,.333){\makebox(0,0)[b]{\footnotesize#4}}
\put(2.5,-1.333){\makebox(0,0)[t]{\footnotesize#5}} 
\end{picture}}

\newcommand{\AlAffod}[4]{\begin{picture}(4,2.6)(0,-.8)
\multiput(0,0)(2,0){3}{\circle*{0.32}}
\put(0,0){\line(1,0){0.5}} \put(1.5,0){\line(1,0){1}} \put(3.5,0){\line(1,0){0.5}}
\put(2,-1){\circle*{0.32}} \put(0,0){\line(2,-1){2}} \put(2,-1){\line(2,1){2}}
\put(2,-1){\circle{0.6}}
\put(1,0){\makebox(0,0){\tiny$\cdots$}} \put(3,0){\makebox(0,0){\tiny$\cdots$}}
\put(0,.333){\makebox(0,0)[b]{\footnotesize#1}}
\put(2,.333){\makebox(0,0)[b]{\footnotesize#2}}
\put(4,.333){\makebox(0,0)[b]{\footnotesize#3}}
\put(2,-1.333){\makebox(0,0)[t]{\footnotesize#4}} 
\end{picture}}

\newcommand{\BthreeAff}[4]{\begin{picture}(1.894,2.494)(-.894,-.347)
\multiput(0,0)(1,0){2}{\circle*{0.32}}
\put(0,.04){\line(1,0){1}} \put(0,-.04){\line(1,0){1}}
\multiput(-.894,.447)(0,-.894){2}{\circle*{0.32}}
\put(-.894,-.447){\circle{0.6}}
\put(0,0){\line(-2,1){1}}\put(0,0){\line(-2,-1){1}}
\put(-.894,.78){\makebox(0,0)[b]{\footnotesize#4}}
\put(-.894,-.78){\makebox(0,0)[t]{\footnotesize#2}}
\put(.5,0){\makebox(0,0){$>$}} 
\put(0,.333){\makebox(0,0)[b]{\footnotesize#1}}
\put(1,.333){\makebox(0,0)[b]{\footnotesize#3}}
\end{picture}}

\newcommand{\BlAff}[8][>]{\begin{picture}(5.894,2.494)(-.894,-.347)
\multiput(0,0)(1,0){2}{\circle*{0.32}} \multiput(3,0)(1,0){3}{\circle*{0.32}}
\put(0,0){\line(1,0){1.5}} \put(2.5,0){\line(1,0){1.5}}
\put(4,.04){\line(1,0){1}} \put(4,-.04){\line(1,0){1}}
\multiput(-.894,.447)(0,-.894){2}{\circle*{0.32}}
\put(-.894,-.447){\circle{0.6}}
\put(0,0){\line(-2,1){1}}\put(0,0){\line(-2,-1){1}}
\put(-.894,.78){\makebox(0,0)[b]{\footnotesize#8}}
\put(-.894,-.78){\makebox(0,0)[t]{\footnotesize#2}}
\put(4.5,0){\makebox(0,0){$#1$}} 
\put(2,0){\makebox(0,0){\tiny$\cdots$}}
\put(0,.333){\makebox(0,0)[b]{\footnotesize#3}}
\put(1,.333){\makebox(0,0)[b]{\footnotesize#4}}
\put(3,.333){\makebox(0,0)[b]{\footnotesize#5}}
\put(4,.333){\makebox(0,0)[b]{\footnotesize#6}}
\put(5,.333){\makebox(0,0)[b]{\footnotesize#7}} 
\end{picture}}

\newcommand{\ClAff}[7][>]{\begin{picture}(6,1.2)(-1,-.3)
\multiput(-1,0)(1,0){3}{\circle*{0.32}} \multiput(3,0)(1,0){3}{\circle*{0.32}}
\put(-1,0){\circle{0.6}}
\put(0,0){\line(1,0){1.5}} \put(2.5,0){\line(1,0){1.5}}
\put(4,.04){\line(1,0){1}} \put(4,-.04){\line(1,0){1}}
\put(-1,.04){\line(1,0){1}} \put(-1,-.04){\line(1,0){1}}
\put(4.5,0){\makebox(0,0){$<$}} \put(-.5,0){\makebox(0,0){$#1$}} 
\put(2,0){\makebox(0,0){\tiny$\cdots$}}
\put(0,.333){\makebox(0,0)[b]{\footnotesize#2}}
\put(1,.333){\makebox(0,0)[b]{\footnotesize#3}}
\put(3,.333){\makebox(0,0)[b]{\footnotesize#4}}
\put(4,.333){\makebox(0,0)[b]{\footnotesize#5}}
\put(5,.333){\makebox(0,0)[b]{\footnotesize#6}} 
\put(-1,.333){\makebox(0,0)[b]{\footnotesize#7}} 
\end{picture}}

\newcommand{\DlAff}[8]{\begin{picture}(5.788,2.494)(-.894,-.347)
\multiput(0,0)(1,0){2}{\circle*{0.32}} \multiput(3,0)(1,0){2}{\circle*{0.32}}
\put(0,0){\line(1,0){1.5}} \put(2.5,0){\line(1,0){1.5}}
\multiput(4.894,.447)(0,-.894){2}{\circle*{0.32}}
\put(4,0){\line(2,1){1}} \put(4,0){\line(2,-1){1}}
\multiput(-.894,.447)(0,-.894){2}{\circle*{0.32}}
\put(-.894,-.447){\circle{0.6}}
\put(0,0){\line(-2,1){1}}\put(0,0){\line(-2,-1){1}}
\put(-.894,.78){\makebox(0,0)[b]{\footnotesize#1}}
\put(-.894,-.78){\makebox(0,0)[t]{\footnotesize#8}}
\put(2,0){\makebox(0,0){\tiny$\cdots$}}
\put(0,.333){\makebox(0,0)[b]{\footnotesize#2}}
\put(1,.333){\makebox(0,0)[b]{\footnotesize#3}}
\put(3,.333){\makebox(0,0)[b]{\footnotesize#4}}
\put(4,.333){\makebox(0,0)[b]{\footnotesize#5}}
\put(4.894,.78){\makebox(0,0)[b]{\footnotesize#6}}
\put(4.894,-.78){\makebox(0,0)[t]{\footnotesize#7}}
\end{picture}}

\newcommand{\EaAff}[7]{\begin{picture}(4,3.1)(0,-.3)
\multiput(0,1)(1,0){5}{\circle*{0.32}}\put(2,0){\circle*{0.32}}\put(2,-1){\circle*{0.32}}
\put(0,1){\line(1,0){4}}\put(2,-1){\line(0,1){2}}
\put(2,-1){\circle{0.6}}
\put(0,1.333){\makebox(0,0)[b]{\footnotesize#1}}
\put(1.667,0){\makebox(0,0)[r]{\footnotesize#2}}
\put(1,1.333){\makebox(0,0)[b]{\footnotesize#3}}
\put(2,1.333){\makebox(0,0)[b]{\footnotesize#4}}
\put(3,1.333){\makebox(0,0)[b]{\footnotesize#5}}
\put(4,1.333){\makebox(0,0)[b]{\footnotesize#6}}
\put(1.667,-1){\makebox(0,0)[r]{\footnotesize#7}}
\end{picture}}

\newcommand{\EbAff}[8]{\begin{picture}(6,2.6)(-1,.8)
\multiput(-1,1)(1,0){7}{\circle*{0.32}}\put(2,0){\circle*{0.32}}
\put(-1,1){\circle{0.6}}
\put(-1,1){\line(1,0){6}}\put(2,0){\line(0,1){1}}
\put(0,1.333){\makebox(0,0)[b]{\footnotesize#1}}
\put(2,-.333){\makebox(0,0)[t]{\footnotesize#2}}
\put(1,1.333){\makebox(0,0)[b]{\footnotesize#3}}
\put(2,1.333){\makebox(0,0)[b]{\footnotesize#4}}
\put(3,1.333){\makebox(0,0)[b]{\footnotesize#5}}
\put(4,1.333){\makebox(0,0)[b]{\footnotesize#6}}
\put(5,1.333){\makebox(0,0)[b]{\footnotesize#7}}
\put(-1,1.333){\makebox(0,0)[b]{\footnotesize#8}}
\end{picture}}

\newcommand{\EcAff}[9]{\begin{picture}(8,2.6)(-.5,.8)
\multiput(0,1)(1,0){8}{\circle*{0.32}}\put(2,0){\circle*{0.32}}
\put(7,1){\circle{0.6}}
\put(0,1){\line(1,0){7}}\put(2,0){\line(0,1){1}}
\put(0,1.333){\makebox(0,0)[b]{\footnotesize#1}}
\put(2,-.333){\makebox(0,0)[t]{\footnotesize#2}}
\put(1,1.333){\makebox(0,0)[b]{\footnotesize#3}}
\put(2,1.333){\makebox(0,0)[b]{\footnotesize#4}}
\put(3,1.333){\makebox(0,0)[b]{\footnotesize#5}}
\put(4,1.333){\makebox(0,0)[b]{\footnotesize#6}}
\put(5,1.333){\makebox(0,0)[b]{\footnotesize#7}}
\put(6,1.333){\makebox(0,0)[b]{\footnotesize#8}}
\put(7,1.333){\makebox(0,0)[b]{\footnotesize#9}}
\end{picture}}

\newcommand{\FAff}[6][>]{\begin{picture}(4,1.2)(-1,-.3)
\multiput(-1,0)(1,0){5}{\circle*{0.32}}
\put(-1,0){\circle{0.6}}
\put(-1,0){\line(1,0){2}} \put(2,0){\line(1,0){1}}
\put(1,.04){\line(1,0){1}} \put(1,-.04){\line(1,0){1}}
\put(1.5,0){\makebox(0,0){$#1$}} 
\put(0,.333){\makebox(0,0)[b]{\footnotesize#2}}
\put(1,.333){\makebox(0,0)[b]{\footnotesize#3}}
\put(2,.333){\makebox(0,0)[b]{\footnotesize#4}}
\put(3,.333){\makebox(0,0)[b]{\footnotesize#5}}
\put(-1,.333){\makebox(0,0)[b]{\footnotesize#6}}
\end{picture}}

\newcommand{\GAff}[4][<]{\begin{picture}(2,1.2)(0,-.3)
\multiput(0,0)(1,0){3}{\circle*{0.32}}
\put(2,0){\circle{0.6}}
\put(0,0){\line(1,0){2}}
\put(0,.08){\line(1,0){1}} \put(0,-.08){\line(1,0){1}}
\put(.5,0){\makebox(0,0){$#1$}} 
\put(0,.333){\makebox(0,0)[b]{\footnotesize#2}}
\put(1,.333){\makebox(0,0)[b]{\footnotesize#3}}
\put(2,.333){\makebox(0,0)[b]{\footnotesize#4}}
\end{picture}}

\newcommand{\Al}[4]{\begin{picture}(3,0)(0,-1)
\multiput(0,0)(1,0){4}{\circle*{0.32}}
\put(0,0){\line(1,0){3}}
\put(0,.333){\makebox(0,0)[b]{\footnotesize#1}}
\put(1,.333){\makebox(0,0)[b]{\footnotesize#2}}
\put(2,.333){\makebox(0,0)[b]{\footnotesize#3}}
\put(3,.333){\makebox(0,0)[b]{\footnotesize#4}}
\end{picture}}

\newcommand{\Bl}[4]{\begin{picture}(3,0)(0,-1)
\multiput(0,0)(1,0){4}{\circle*{0.32}}
\put(0,0){\line(1,0){2}}
\put(2,0.07){\line(1,0){1}} \put(2,-.07){\line(1,0){1}}
\put(2.5,0){\makebox(0,0){\footnotesize$>$}} 
\put(0,.333){\makebox(0,0)[b]{\footnotesize#1}}
\put(1,.333){\makebox(0,0)[b]{\footnotesize#2}}
\put(2,.333){\makebox(0,0)[b]{\footnotesize#3}}
\put(3,.333){\makebox(0,0)[b]{\footnotesize#4}}
\end{picture}}

\newcommand{\Cl}[4]{\begin{picture}(3,0)(0,-1)
\multiput(0,0)(1,0){4}{\circle*{0.32}}
\put(0,0){\line(1,0){2}}
\put(2,0.07){\line(1,0){1}} \put(2,-.07){\line(1,0){1}}
\put(2.5,0){\makebox(0,0){\footnotesize$<$}} 
\put(0,.333){\makebox(0,0)[b]{\footnotesize#1}}
\put(1,.333){\makebox(0,0)[b]{\footnotesize#2}}
\put(2,.333){\makebox(0,0)[b]{\footnotesize#3}}
\put(3,.333){\makebox(0,0)[b]{\footnotesize#4}}
\end{picture}}

\newcommand{\Dl}[4]{\begin{picture}(1.894,.894)(0,-1)
\multiput(0,0)(1,0){2}{\circle*{0.32}}
\put(0,0){\line(1,0){1}}
\multiput(1.894,.447)(0,-.894){2}{\circle*{0.32}}
\put(1,0){\line(2,1){1}} \put(1,0){\line(2,-1){1}} 
\put(0,.333){\makebox(0,0)[b]{\footnotesize#1}}
\put(1,.333){\makebox(0,0)[b]{\footnotesize#2}}
\put(1.894,.78){\makebox(0,0)[b]{\footnotesize#3}}
\put(1.894,-.78){\makebox(0,0)[t]{\footnotesize#4}}
\end{picture}}

\newcommand{\Fl}[4]{\begin{picture}(3,0)(0,-1)
\multiput(0,0)(1,0){4}{\circle*{0.32}}
\put(0,0){\line(1,0){1}}\put(2,0){\line(1,0){1}}
\put(1,0.07){\line(1,0){1}} \put(1,-.07){\line(1,0){1}}
\put(1.5,0){\makebox(0,0){\footnotesize$>$}} 
\put(0,.333){\makebox(0,0)[b]{\footnotesize#1}}
\put(1,.333){\makebox(0,0)[b]{\footnotesize#2}}
\put(2,.333){\makebox(0,0)[b]{\footnotesize#3}}
\put(3,.333){\makebox(0,0)[b]{\footnotesize#4}}
\end{picture}}

\newcommand{\Agenl}[4]{\raisebox{-.8\unitlength}{\Fbox{
\begin{picture}(5,1)(-.5,-.3)
\multiput(0,0)(1,0){2}{\circle*{0.15}} \multiput(3,0)(1,0){2}{\circle*{0.15}}
\put(0,0){\line(1,0){1.5}} \put(2.5,0){\line(1,0){1.5}}
\put(2,0){\makebox(0,0){$\cdots$}}
\put(0,.2){\makebox(0,0)[b]{#1}}
\put(1,.2){\makebox(0,0)[b]{#2}}
\put(3,.2){\makebox(0,0)[b]{#3}}
\put(4,.2){\makebox(0,0)[b]{#4}} 
\end{picture}}}}
\newcommand{\Afuenf}[5]{
\begin{picture}(5,1)(-.5,-.3)
\multiput(0,0)(1,0){5}{\circle*{0.15}}
\put(0,0){\line(1,0){4}}
\put(0,.2){\makebox(0,0)[b]{#1}}
\put(1,.2){\makebox(0,0)[b]{#2}}
\put(2,.2){\makebox(0,0)[b]{#3}}
\put(3,.2){\makebox(0,0)[b]{#4}}
\put(4,.2){\makebox(0,0)[b]{#5}} 
\end{picture}}

\newcommand{\Bgenl}[5]{\raisebox{-.8\unitlength}{\Fbox{
\begin{picture}(6,1)(-.5,-.3)
\multiput(0,0)(1,0){2}{\circle*{0.15}} \multiput(3,0)(1,0){3}{\circle*{0.15}}
\put(0,0){\line(1,0){1.5}} \put(2.5,0){\line(1,0){1.5}}
\put(4,0.03){\line(1,0){1}} \put(4,-.03){\line(1,0){1}}
\put(4.5,0){\makebox(0,0){\LARGE$>$}} 
\put(2,0){\makebox(0,0){$\cdots$}}
\put(0,.2){\makebox(0,0)[b]{#1}}
\put(1,.2){\makebox(0,0)[b]{#2}}
\put(3,.2){\makebox(0,0)[b]{#3}}
\put(4,.2){\makebox(0,0)[b]{#4}}
\put(5,.2){\makebox(0,0)[b]{#5}} 
\end{picture}}}}
\newcommand{\Bfuenf}[5]{
\begin{picture}(5,1)(-.5,-.3)
\multiput(0,0)(1,0){5}{\circle*{0.15}}
\put(0,0){\line(1,0){3}}
\put(3,0.03){\line(1,0){1}} \put(3,-.03){\line(1,0){1}}
\put(3.5,0){\makebox(0,0){\LARGE$>$}} 
\put(0,.2){\makebox(0,0)[b]{#1}}
\put(1,.2){\makebox(0,0)[b]{#2}}
\put(2,.2){\makebox(0,0)[b]{#3}}
\put(3,.2){\makebox(0,0)[b]{#4}}
\put(4,.2){\makebox(0,0)[b]{#5}} 
\end{picture}}

\newcommand{\Cgenl}[5]{\raisebox{-.8\unitlength}{\Fbox{
\begin{picture}(6,1)(-.5,-.3)
\multiput(0,0)(1,0){2}{\circle*{0.15}} \multiput(3,0)(1,0){3}{\circle*{0.15}}
\put(0,0){\line(1,0){1.5}} \put(2.5,0){\line(1,0){1.5}}
\put(4,0.03){\line(1,0){1}} \put(4,-.03){\line(1,0){1}}
\put(4.5,0){\makebox(0,0){\LARGE$<$}} 
\put(2,0){\makebox(0,0){$\cdots$}}
\put(0,.2){\makebox(0,0)[b]{#1}}
\put(1,.2){\makebox(0,0)[b]{#2}}
\put(3,.2){\makebox(0,0)[b]{#3}}
\put(4,.2){\makebox(0,0)[b]{#4}}
\put(5,.2){\makebox(0,0)[b]{#5}} 
\end{picture}}}}
\newcommand{\Cfuenf}[5]{
\begin{picture}(5,1)(-.5,-.3)
\multiput(0,0)(1,0){5}{\circle*{0.15}}
\put(0,0){\line(1,0){3}}
\put(3,0.03){\line(1,0){1}} \put(3,-.03){\line(1,0){1}}
\put(3.5,0){\makebox(0,0){\LARGE$<$}} 
\put(0,.2){\makebox(0,0)[b]{#1}}
\put(1,.2){\makebox(0,0)[b]{#2}}
\put(2,.2){\makebox(0,0)[b]{#3}}
\put(3,.2){\makebox(0,0)[b]{#4}}
\put(4,.2){\makebox(0,0)[b]{#5}} 
\end{picture}}

\newcommand{\Dgenl}[6]{\raisebox{-1.8\unitlength}{\Fbox{
\begin{picture}(5.894,2.1)(-.5,-1.05)
\multiput(0,0)(1,0){2}{\circle*{0.15}} \multiput(3,0)(1,0){2}{\circle*{0.15}}
\put(0,0){\line(1,0){1.5}} \put(2.5,0){\line(1,0){1.5}}
\multiput(4.894,.447)(0,-.894){2}{\circle*{0.15}}
\put(4,0){\line(2,1){0.894}} \put(4,0){\line(2,-1){0.894}} 
\put(2,0){\makebox(0,0){$\cdots$}}
\put(0,.2){\makebox(0,0)[b]{#1}}
\put(1,.2){\makebox(0,0)[b]{#2}}
\put(3,.2){\makebox(0,0)[b]{#3}}
\put(4,.2){\makebox(0,0)[b]{#4}}
\put(4.894,.65){\makebox(0,0)[b]{#5}}
\put(4.894,-.65){\makebox(0,0)[t]{#6}}
\end{picture}}}}
\newcommand{\Dfuenf}[5]{
\begin{picture}(3.894,1.7)(-.5,-.85)
\multiput(0,0)(1,0){3}{\circle*{0.15}}
\put(0,0){\line(1,0){2}}
\multiput(2.894,.447)(0,-.894){2}{\circle*{0.15}}
\put(2,0){\line(2,1){0.894}}\put(2,0){\line(2,-1){0.894}} 
\put(0,.2){\makebox(0,0)[b]{#1}}
\put(1,.2){\makebox(0,0)[b]{#2}}
\put(2,.2){\makebox(0,0)[b]{#3}}
\put(2.894,.65){\makebox(0,0)[b]{#4}}
\put(2.894,-.65){\makebox(0,0)[t]{#5}}
\end{picture}}

\newcommand{\Esechs}[6]{\raisebox{-2\unitlength}{\Fbox{
\begin{picture}(5,2.4)(-.5,-.7)
\multiput(0,1)(1,0){5}{\circle*{0.15}}\put(2,0){\circle*{0.15}}
\put(0,1){\line(1,0){4}}\put(2,0){\line(0,1){1}}
\put(0,1.2){\makebox(0,0)[b]{#1}}
\put(2,-.2){\makebox(0,0)[t]{#2}}
\put(1,1.2){\makebox(0,0)[b]{#3}}
\put(2,1.2){\makebox(0,0)[b]{#4}}
\put(3,1.2){\makebox(0,0)[b]{#5}}
\put(4,1.2){\makebox(0,0)[b]{#6}}
\end{picture}}}}

\newcommand{\Esieben}[7]{\raisebox{-2\unitlength}{\Fbox{
\begin{picture}(6,2.4)(-.5,-.7)
\multiput(0,1)(1,0){6}{\circle*{0.15}}\put(2,0){\circle*{0.15}}
\put(0,1){\line(1,0){5}}\put(2,0){\line(0,1){1}}
\put(0,1.2){\makebox(0,0)[b]{#1}}
\put(2,-.2){\makebox(0,0)[t]{#2}}
\put(1,1.2){\makebox(0,0)[b]{#3}}
\put(2,1.2){\makebox(0,0)[b]{#4}}
\put(3,1.2){\makebox(0,0)[b]{#5}}
\put(4,1.2){\makebox(0,0)[b]{#6}}
\put(5,1.2){\makebox(0,0)[b]{#7}}
\end{picture}}}}

\newcommand{\Eacht}[8]{\raisebox{-2\unitlength}{\Fbox{
\begin{picture}(7,2.4)(-.5,-.7)
\multiput(0,1)(1,0){7}{\circle*{0.15}}\put(2,0){\circle*{0.15}}
\put(0,1){\line(1,0){6}}\put(2,0){\line(0,1){1}}
\put(0,1.2){\makebox(0,0)[b]{#1}}
\put(2,-.2){\makebox(0,0)[t]{#2}}
\put(1,1.2){\makebox(0,0)[b]{#3}}
\put(2,1.2){\makebox(0,0)[b]{#4}}
\put(3,1.2){\makebox(0,0)[b]{#5}}
\put(4,1.2){\makebox(0,0)[b]{#6}}
\put(5,1.2){\makebox(0,0)[b]{#7}}
\put(6,1.2){\makebox(0,0)[b]{#8}}
\end{picture}}}}

\newcommand{\Fvier}[4]{\raisebox{-.6\unitlength}{\Fbox{
\begin{picture}(4,1)(-.5,-.3)
\multiput(0,0)(1,0){4}{\circle*{0.15}}
\put(0,0){\line(1,0){1}} \put(2,0){\line(1,0){1}}
\put(1,0.03){\line(1,0){1}} \put(1,-.03){\line(1,0){1}}
\put(1.5,0){\makebox(0,0){\LARGE$>$}} 
\put(0,.2){\makebox(0,0)[b]{#1}}
\put(1,.2){\makebox(0,0)[b]{#2}}
\put(2,.2){\makebox(0,0)[b]{#3}}
\put(3,.2){\makebox(0,0)[b]{#4}}
\end{picture}}}}

\newcommand{\Gzwei}[2]{\raisebox{-.6\unitlength}{\Fbox{
\begin{picture}(2,1)(-.5,-.3)
\multiput(0,0)(1,0){2}{\circle*{0.15}}
\put(0,0){\line(1,0){1}}
\put(0,.06){\line(1,0){1}} \put(0,-.06){\line(1,0){1}}
\put(.5,0){\makebox(0,0){\LARGE$<$}} 
\put(0,.2){\makebox(0,0)[b]{#1}}
\put(1,.2){\makebox(0,0)[b]{#2}}
\end{picture}}}}

\newcommand{\Axxxx}{\mbox{$\Al1111$}}
\newcommand{\Axxxo}{\mbox{$\Al1110$}}
\newcommand{\Axxoo}{\mbox{$\Al1100$}}
\newcommand{\Axooo}{\mbox{$\Al1000$}}
\newcommand{\Aoxxx}{\mbox{$\Al0111$}}
\newcommand{\Aoxxo}{\mbox{$\Al0110$}}
\newcommand{\Aoxoo}{\mbox{$\Al0100$}}
\newcommand{\Aooxx}{\mbox{$\Al0011$}}
\newcommand{\Aooxo}{\mbox{$\Al0010$}}
\newcommand{\Aooox}{\mbox{$\Al0001$}}
\newcommand{\Bxwww}{\mbox{$\Bl1222$}}
\newcommand{\Bxxww}{\mbox{$\Bl1122$}}
\newcommand{\Bxxxw}{\mbox{$\Bl1112$}}
\newcommand{\Bxxxx}{\mbox{$\Bl1111$}}
\newcommand{\Bxxxo}{\mbox{$\Bl1110$}}
\newcommand{\Bxxoo}{\mbox{$\Bl1100$}}
\newcommand{\Bxooo}{\mbox{$\Bl1000$}}
\newcommand{\Boxww}{\mbox{$\Bl0122$}}
\newcommand{\Boxxw}{\mbox{$\Bl0112$}}
\newcommand{\Booxw}{\mbox{$\Bl0012$}}
\newcommand{\Cwwwx}{\mbox{$\Cl2221$}}
\newcommand{\Cxwwx}{\mbox{$\Cl1221$}}
\newcommand{\Cxxwx}{\mbox{$\Cl1121$}}
\newcommand{\Cxxxx}{\mbox{$\Cl1111$}}
\newcommand{\Cowwx}{\mbox{$\Cl0221$}}
\newcommand{\Coxwx}{\mbox{$\Cl0121$}}
\newcommand{\Coxxx}{\mbox{$\Cl0111$}}
\newcommand{\Coowx}{\mbox{$\Cl0021$}}
\newcommand{\Cooxx}{\mbox{$\Cl0011$}}
\newcommand{\Cooox}{\mbox{$\Cl0001$}}
\newcommand{\Dxwxx}{\mbox{$\Dl1211$}}
\newcommand{\Dxxxx}{\mbox{$\Dl1111$}}
\newcommand{\Dxxox}{\mbox{$\Dl1101$}}
\newcommand{\Dxxxo}{\mbox{$\Dl1110$}}
\newcommand{\Dxxoo}{\mbox{$\Dl1100$}}
\newcommand{\Dxooo}{\mbox{$\Dl1000$}}
\newcommand{\Doxxx}{\mbox{$\Dl0111$}}
\newcommand{\Doxox}{\mbox{$\Dl0101$}}
\newcommand{\Doxxo}{\mbox{$\Dl0110$}}
\newcommand{\Dooox}{\mbox{$\Dl0001$}}
\newcommand{\Dooxo}{\mbox{$\Dl0010$}}
\newcommand{\Fwhfw}{\mbox{$\Fl2342$}}
\newcommand{\Fxhfw}{\mbox{$\Fl1342$}}
\newcommand{\Fxwfw}{\mbox{$\Fl1242$}}
\newcommand{\Fxwww}{\mbox{$\Fl1222$}}
\newcommand{\Fxxww}{\mbox{$\Fl1122$}}
\newcommand{\Foxww}{\mbox{$\Fl0122$}}
\newcommand{\Fxwhw}{\mbox{$\Fl1232$}}
\newcommand{\Fxwhx}{\mbox{$\Fl1231$}}
\newcommand{\Fxwwx}{\mbox{$\Fl1221$}}
\newcommand{\Fxwwo}{\mbox{$\Fl1220$}}

\newcommand{\refla}{\mbox{\footnotesize$s_0$}}
\newcommand{\reflb}{\mbox{\footnotesize$s_1$}}
\newcommand{\reflc}{\mbox{\footnotesize$s_2$}}
\newcommand{\refld}{\mbox{\footnotesize$s_3$}}
\newcommand{\refle}{\mbox{\footnotesize$s_4$}}

\newcommand{\reflg}{\mbox{\footnotesize$s_6$}}
\newcommand{\reflh}{\mbox{\footnotesize$s_7$}}
\newcommand{\refli}{\mbox{\footnotesize$s_8$}}
\newcommand{\reflj}{\mbox{\footnotesize$s_9$}}
\newcommand{\reflk}{\mbox{\footnotesize$s_{10}$}}
\newcommand{\refll}{\mbox{\footnotesize$s_{11}$}}

\begin{document}

\title{\textbf{Abelian ideals in a Borel subalgebra}\\
\textbf{of a complex simple Lie algebra}}
\author{\large Ruedi Suter\\[2mm]
\small Departement Mathematik, ETH Z\"urich,\\
\small ETH Zentrum, 8092 Z\"urich, Switzerland\\
\small\texttt{suter@math.ethz.ch}}
\def\today{}
\maketitle

\begin{abstract}
Let $\lie$ be a complex simple Lie algebra and $\borel$ a fixed Borel subalgebra
of $\lie$. We shall describe the abelian ideals of $\borel$ in a uniform way,
that is, independent of the classification of complex simple Lie algebras.
\end{abstract}

\section{Introduction}
Let $\lie$ be a complex simple Lie algebra and $\borel$ a fixed Borel subalgebra
of $\lie$.

This paper has three purposes. First, the maximal dimension among the abelian
ideals in $\borel$ is determined purely in terms of certain invariants
(dual Coxeter number and numbers of positive roots) of the root system
and some associated root subsystems.
To assuage any possible curiosity we already list the maximal dimensions together
with their computations for the five exceptional types. The whole picture
will be revealed in the table on page~\pageref{maxdim}.
\begin{align*}
g_{\sfE_6}-1+N_{\sfA_5}-N_{\sfA_4}&=12-1+15-10=16\\
g_{\sfE_7}-1+N_{\sfD_6}-N_{\sfD_5}&=18-1+20-10=27\\
g_{\sfE_8}-1+N_{\sfA_7}-N_{\sfA_6}&=30-1+28-21=36\\
g_{\sfF_4}-1+N_{\sfA_1}-N_{\varnothing\,\,}&=\phantom{0}9-1+\phantom{0}1
-\phantom{0}0=\phantom{0}9\\
g_{\sfG_2}-1+N_{\varnothing\,\,}-N_{\varnothing\,\,}&=\phantom{0}4-1+\phantom{0}0
-\phantom{0}0=\phantom{0}3
\end{align*}

Second, we answer a question of Panyushev and R\"ohrle, asking for a uniform
explanation for the one-to-one correspondence between the maximal abelian
ideals in $\borel$ and the long simple roots. More generally, in our approach
all positive long roots will emerge in a natural way.\footnote{In the very recent
preprint \cite{Pa} Panyushev answers this question. He finds the natural mapping
from the set of nonzero abelian ideals onto the set of positive long roots, too.
He then proves that each fibre of this mapping is a poset having a unique maximal
element and a unique minimal element and asks for further investigating the poset
structure and, in particular, to find a general description of the maximal element
of each fibre. Our main theorem solves this problem.}

Third, we keep the promise\footnote{The announcement was made in \cite{Su} as well as in
a talk one year ago in the algebra and topology seminar at ETH Zurich. I decided
to post the present version of the paper rather than keep it back any longer.}
of giving a generalization and explanation
of the symmetry property of a certain subposet of Young's lattice
(the lattice of integer partitions) that
was observed in \cite{Su} and which we now recall. For that
consider the subposet $\mathbb Y_N$ of Young's lattice
induced by the Young diagrams whose
(largest) hook lengths are at most $N-1$. One sees easily that the poset
$\mathbb Y_N$ has $2^{N-1}$ elements. This follows for instance by
associating to each such diagram an integer between $0$ and $2^{N-1}-1$ by
the following procedure: in each column of the diagram
write the figure~$1$ at the bottom and fill
the rest by $0$; then read the binary number along the rim.
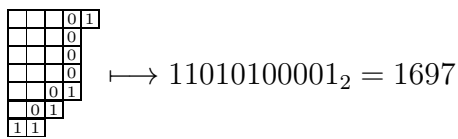
\begin{center}
$\raisebox{-7mm}{\tiny
\young(\ \ \ 01,\ \ \ 0,\ \ \ 0,\ \ \ 0,\ \ 01,\ 01,11)}
\longmapsto11010100001_2=1697$
\end{center}
The main result of \cite{Su} states that
the Hasse graph of $\mathbb Y_N$
(considered as an undirected graph) has the dihedral group $\Dih_N$ of order $2N$ as its
automorphism group provided $N\geqslant3$. The following figure exemplifies this
fact for~$N=5$.
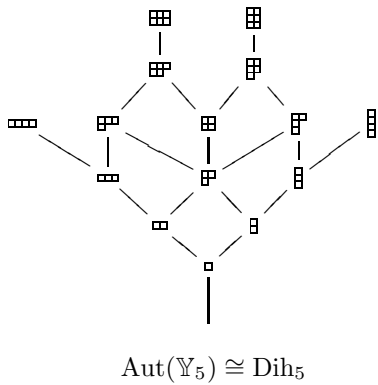
\begin{center}
\mbox{}{\tiny\Yboxdim1mm
\xymatrix@R=2mm@C=2mm{&&&{\yng(3,3)}&&{\yng(2,2,2)}\\
&&&{\yng(3,2)}\ar@{-}[u]&&{\yng(2,2,1)}\ar@{-}[u]\\
{\yng(4)}&&{\yng(3,1)}\ar@{-}[ur]&&{\yng(2,2)}\ar@{-}[ul]\ar@{-}[ur]
&&{\yng(2,1,1)}\ar@{-}[ul]&&{\yng(1,1,1,1)}\\
&&{\yng(3)}\ar@{-}[ull]\ar@{-}[u]&&{\yng(2,1)}\ar@{-}[ull]
\ar@{-}[u]\ar@{-}[urr]&&{\yng(1,1,1)}\ar@{-}[u]\ar@{-}[urr]\\
&&&{\yng(2)}\ar@{-}[ul]\ar@{-}[ur]
&&{\yng(1,1)}\ar@{-}[ul]\ar@{-}[ur]\\
&&&&{\yng(1)}\ar@{-}[ul]\ar@{-}[ur]\\\\
&&&&{\mbox{}}\ar@{-}[uu]}}\\[1ex]
{\footnotesize \ \ \ \ \ $\Aut(\mathbb Y_5)\cong\Dih_5$}
\end{center}
The case dealt with in \cite{Su} is now seen as the
$\mathsf A_{N-1}$ case, i.\,e., associated with the Lie algebra
$\lie=\mathfrak{sl}_N(\mathbb C)$. It is so to say the most spectacular
case. The reason is that its affine Coxeter-Dynkin graph is a cycle
of length $N$. Its dihedral symmetry induces a dihedral symmetry on a
certain simplicial complex $\mathcal C$. The Hasse graph of $\mathbb Y_N$ can be
geometrically realized as the $1$-skeleton of the cell
complex dual to $\mathcal C$.

Let me close this section with a brief historical narrative.
In his 1905 paper \cite{Sch} in Crelle's journal I.~Schur proved that the maximum number of linearly
independent commuting $N$\/$\times$\/$N$ matrices is $\bigl\lfloor\frac{N^2}{4}
\bigr\rfloor+1$. In 1944 Jacobson \cite{Ja} gave a simplified derivation of Schur's result.
In the next year A.~Malcev \cite{Mal} determined the commutative subalgebras of maximum dimension
of the semisimple complex Lie groups, or equivalently, their Lie algebras. The next
entry in this short historical outline is Kostant's paper \cite{K1} published
in 1965. There he gave a connexion of Malcev's result with the maximal eigenvalue of the Laplacian acting
on the exterior powers $\bigwedge^k\lie$ of the adjoint representation.
Kostant \cite{K2} again, in 1998, reconsidered the theme of abelian
ideals in a Borel subalgebra of $\lie$ and reported inter alia about Peterson's
proof that the number of abelian ideals in a fixed Borel subalgebra of
$\lie$ is $2^{\operatorname{rank}\lie}$. This, in Kostant's words,
utterly surprising and ingenious proof involves the affine Weyl group.

A natural generalization of Peterson's approach from abelian
to $\operatorname{ad}$-nilpo\-tent ideals was developed recently
by several authors in various collaborations \cite{AKOP,CP1,CP2,CP3,KOP},
see also~\cite{Sh}, and for Kostant's results \cite{CMP}.

\section{Notations and tools}
Basic facts concerning root systems can be found in the standard references
\cite{B,Hu} and also in \cite{Br,Hi}.
Let us fix a complex simple Lie algebra $\lie$ of rank $l$ together with a
Borel subalgebra $\borel$ and a Cartan subalgebra $\cartan
\subseteq\borel$. Associated with these data there are quite a number
of further objects whose notations are provided in the following table.
Most of them are standard (but sometimes there are different conventions).
We list the most important notations used here for the reader's convenience.

\begin{longtable}{ll}
$A$&(closed) fundamental alcove,\\
$\operatorname{Alt}_N$&alternating group of degree $N$,\\
$C$&(closed) dominant chamber,\\
$\Dih_N$&dihedral group of order $2N$,\\
$F_i$&facets of type $i$ ($i=0,\dots,l$) of the fundamental alcove $A$,\\
$g$&dual Coxeter number,\\
$\lie_\varphi$&root subspace,\\
$h$&Coxeter number,\\
$\cartan_\RR^\ast$, $\cartan_\RR$&real vector space spanned by the roots, and its predual,\\
$H_i$&hyperplanes supporting the facets $F_i$ ($i=0,\dots,l$) of $A$,\\
$l$&$=\operatorname{rk}\lie$,\\
$\ell$&length function on $\widehat W$ or on $W$,\\
$L$&$L(\varphi)=\frac{2\,(\theta-\varphi|\rho)}{(\theta|\theta)}$\,,\\
$m_i$&exponents ($i=1,\dots,l$),\\
$n_i$&marks ($i=1,\dots,l$), $\theta=\sum\limits_{i=1}^ln_i\,\alpha_i$,
in addition, $n_0=1$,\\
$N_\sfX$&number of positive roots for a root system of type $\sfX$,\\
$s_i$&simple reflections ($i=1,\dots,l$), in addition,
$s_0:\cartan_\RR^*\to\cartan_\RR^*$,\\
&$s_0(\lambda)=\lambda-\langle\lambda,\theta^\vee\rangle\,\theta
+g\theta$,\\
$s_\theta$&reflection along the highest root, $s_\theta(\lambda)
=s_0(\lambda)-s_0(0)$,\\
$\operatorname{Sym}_N$&symmetric group of degree $N$,\\
$w_\circ^\varphi$, $\widehat w_\circ^\varphi$&longest elements in
$W_{\perp\varphi}$, $\widehat W_{\perp\varphi}$, respectively,\\
$W$, $\widehat W$&finite Weyl group, affine Weyl group,\\
$W_{\perp\varphi}$&$=\operatorname{gp}\bigl(s_i\bigm|\alpha_i\perp\varphi
\mbox{ ($i={1,\dots,l}$)}\bigr)$,\\
$\widehat W_{\perp\varphi}$&$\widehat W_{\perp\varphi}=W_{\perp\varphi}$
if $\theta\not\perp\varphi$,
$\widehat W_{\perp\varphi}=\operatorname{gp}\bigl(s_0,W_{\perp\varphi}\bigr)$
if $\theta\perp\varphi$,\\
\\
$\alpha_i$&simple roots ($i=1,\dots,l$), in addition, $\alpha_0=-\theta$,\\
$\theta$&highest root,\\
$\varpi_i$&fundamental weights ($i=1,\dots,l$),\\
$\stackrel{\vee}{\varpi}_i$&$=\frac{1}{\Vert\alpha_i\Vert^2}\,\varpi_i$,\\
$\Pi$&$=\{\alpha_1,\dots,\alpha_l\}$ set of simple roots,\\
$\rho$&half the sum of positive roots,\\
$\varphi^\vee$&coroot corresponding to $\varphi\in\Phi$,\\
$\Phi$, $\Phi_\pm$&root system, set of positive/negative roots,\\
$\Phi(\ab)$&positive roots such that $\ab=\bigoplus\limits_{\varphi\in
\Phi(\ab)}\lie_\varphi$,\\
$\Phi_w$&$=\Phi_+\cap w\Phi_-$,\\
$\widehat\Phi$, $\widehat\Phi_\pm$&affine root system, set of positive/negative
affine roots,\\
$\widehat\Phi_{\widehat w}$&$=\widehat\Phi_+\cap\widehat w\widehat\Phi_-$,\\
$(\phantom{\alpha}|\phantom{\alpha})$&canonical bilinear form on $\cartan_\RR^*$,\\
$\Vert\phantom{\alpha}\Vert$&norm from $(\phantom{\alpha}|\phantom{\alpha})$,\\
$\langle\phantom{\alpha},\phantom{\alpha}\rangle$&natural pairing,\\
$[\phantom{n}]$&$[n]=\frac{1-t^n}{1-t}$\,.
\end{longtable}

We denote by $\Phi_+\subseteq\cartan^\ast$
the set of positive roots. Here the convention is that the root
subspaces of $\borel$ belong to positive roots, i.\,e.,
$\borel=\cartan\oplus\bigoplus\limits_{\varphi\in\Phi_+}
\lie_\varphi$, where $\lie_\varphi$ is the ($1$-dimensional) root
subspace on which $\cartan$ acts by the weight $\varphi$, that is,
$\lie_\varphi=\bigl\{X\in\lie\bigm|[H,X]=\varphi(H)\,X\ \forall\,
H\in\cartan\bigr\}$.
As further pieces of notation we write $\Phi_-=-\Phi_+$ for the set of
negative roots, $\Phi=\Phi_+\amalg\Phi_-$ for the root system of $\lie$
relative to $\cartan$, and $\Pi\subseteq\Phi_+$ for the root basis.
Recall that $\Pi$ consists of the roots in $\Phi_+$ that lie on the
edges of the polyhedral (in fact, simplicial) cone spanned by the vectors
in $\Phi_+$. Each positive root is a linear combination of the vectors
in $\Pi$ with nonnegative integral coefficients. The Weyl group of $\Phi$
will be denoted by $W$. More about Weyl groups and some geometry associated
with them will be recalled at the appropriate place below.

Now let $\ab\unlhd\borel$ be an ideal. It is
$\operatorname{ad}\cartan$-stable and hence compatible with
the root space decomposition. If we further require that $\ab$
lies in the nilpotent radical $\mathfrak n=[\borel,\borel]$, we get
that $\ab$ is of the form $\ab=\bigoplus\limits_{\varphi\in\Psi}
\lie_\varphi$ for some subset $\Psi\subseteq\Phi_+$ of positive
roots. The ideal property of $\ab$ translates into the
condition for $\Psi$ that $\Psi\stackrel{.}{+}\Phi_+:=
(\Psi+\Phi_+)\cap\Phi_+\subseteq\Psi$. If, in addition, $\ab$ is
supposed to be abelian (so that $\ab\subseteq
[\borel,\borel]$ holds automatically),
we must have $\Psi\stackrel{.}{+}\Psi:=
(\Psi+\Psi)\cap\Phi_+=\varnothing$. It is plain that there is the following
bijection.
\begin{align*}
\left\{\mbox{\begin{tabular}{@{}l@{}}
subsets $\Psi\subseteq\Phi_+$ such that\\
$\Psi\stackrel{.}{+}\Phi_+\subseteq\Psi$ and $\Psi\stackrel{.}{+}\Psi
=\varnothing$
\end{tabular}}\right\}
&\stackrel{\cong}{\longleftrightarrow}
\left\{\mbox{\begin{tabular}{@{}l@{}}
abelian ideals $\ab\unlhd\borel$
\end{tabular}}\right\}\\
\Psi\ &\longmapsto\ \ab_\Psi:=\bigoplus_{\varphi\in\Psi}
\lie_\varphi
\end{align*}

\paragraph{\textit{The inner product.}}
Before we can go on and state Konstant's theorem, which will be an
essential tool for our approach, we recall the canonical inner product
on the real vector space $\cartan_\RR^\ast$ spanned by the (finite) irreducible (reduced)
root system $\Phi$. This inner product will be denoted by
$(\phantom{\alpha}|\phantom{\alpha})$ and the associated Euclidean
norm by $\Vert\phantom{\alpha}\Vert$. It is characterized by being
$W$-invariant and satisfying the normalization
\mbox{$\Vert\rho+\theta\Vert^2-\Vert\rho\Vert^2=1$} where $\rho$ is half
the sum of positive roots and $\theta$ is the highest root.

\begin{myremark}
The canonical inner product is the restriction to $\cartan_\RR^*$ of
the symmetric bilinear form dual to the Killing form of $\lie$.
There are several alternative descriptions of the same normalization.
Here is a short list.
\begin{enumerate}
\item[] $\Vert\rho+\theta\Vert^2-\Vert\rho\Vert^2=1$, i.\,e., the eigenvalue
of the Casimir operator associated to the Killing form is $1$ for the
adjoint representation;
\item[] $\Vert\theta\Vert^{-2}=g$, the dual Coxeter number;
\item[] $\Vert\rho\Vert^2=\frac{1}{24}\dim\lie$, the
``strange formula'' of Freudenthal and de Vries;
\item[] $\sum\limits_{\varphi\in\Phi}\Vert\varphi\Vert^2=\rank\lie$,
a formula due to G.~Brown;
\item[] $\Vert\theta\Vert^2+\sum\limits_{i=1}^ln_i\,\Vert\alpha_i\Vert^2=1$,
where $n_1,\dots,n_l$ are the marks and $\alpha_1,\dots,\alpha_l$ the
simple roots. (The formula looks funnier if one substitutes $\sum\limits_{i=1}^l
n_i\,\alpha_i$ for $\theta$.) One can show the formula by writing
$\Vert\theta\Vert^{-2}=g$ and using the connexion between
the Coxeter number and the dual Coxeter number. Another derivation will be
given in the remark on page~\pageref{connexiongh}.
\end{enumerate}
\end{myremark}

\begin{mydefinition}
For $\Psi\subseteq\Phi_+$ we define its root sum $\langle\Psi\rangle:=
\sum\limits_{\varphi\in\Psi}\varphi$.
\end{mydefinition}

\begin{mylemma}[Kostant]
Let $\Psi_i\subseteq\Phi_+$ with $\Psi_i\stackrel{.}{+}\Phi_+\subseteq\Psi_i$
($i={1,2}$) (ideals) such that $\langle\Psi_1\rangle=\langle\Psi_2\rangle$.
Then $\Psi_1=\Psi_2$.
\end{mylemma}
\begin{proof}
Let $\Psi:=\Psi_1\cap\Psi_2$. Assume to the contrary that $\Psi_1\neq\Psi_2$.
Then since $\langle\Psi_1\rangle=\langle\Psi_2\rangle$ both $\Psi_1-\Psi$ and
$\Psi_2-\Psi$ are nonempty. Let $\varphi_i\in\Psi_i-\Psi$ ($i={1,2}$).
We must have $(\varphi_1|\varphi_2)\leqslant0$. Otherwise $\varphi_1-\varphi_2$
would be a root which can be assumed positive by possibly interchanging the
indices $1$ and $2$. By the ideal property
$\Psi_i\stackrel{.}{+}\Phi_+\subseteq\Psi_i$ we then have
$\varphi_1=\varphi_2+(\varphi_1-\varphi_2)\in\Psi_2$, a contradiction.
Thus $(\varphi_1|\varphi_2)\leqslant0$. Hence since
$\langle\Psi_1-\Psi\rangle=\langle\Psi_2-\Psi\rangle$ we obtain
$$0\leqslant\bigl\Vert\langle\Psi_i-\Psi\rangle\bigr\Vert^2=
\bigl(\langle\Psi_1-\Psi\rangle\bigm|\langle\Psi_2-\Psi\rangle\bigr)
\leqslant0$$
and so $\Psi=\Psi_1=\Psi_2$.
\end{proof}

Like the previous lemma the following theorem is due to Kostant and was
published in 1965.

\begin{mytheorem}[Kostant]\label{criterion}
Let $\Psi\subseteq\Phi_+$ be a set of positive
roots. Let $\ab_\Psi:=\bigoplus\limits_{\varphi\in\Psi}
\lie_\varphi\subseteq\borel$ be the corresponding subspace. Then one always
has the inequality
$$\Bigl\Vert\rho+
\textstyle\sum\limits_{\varphi\in\Psi}\varphi\Bigr\Vert^2-
\Vert\rho\Vert^2\leqslant|\Psi|$$
with equality if and only if $\ab_\Psi$ is an abelian ideal of $\borel$
(and every abelian ideal of $\borel$ is of this form).
\end{mytheorem}

In particular, one recovers the normalization
$\Vert\rho+\theta\Vert^2-\Vert\rho\Vert^2=1$ because $\lie_\theta\unlhd\borel$
is an abelian ideal.

\paragraph{\textit{Reflections and Weyl groups.}}
We will see that each sum $\rho+\sum\limits_{\varphi\in\Psi}\varphi=\rho+
\langle\Psi\rangle$ that occurs in Kostant's theorem and such that
$\ab_\Psi$ is an abelian ideal in $\borel$ can be written as
$\rho+\langle\Psi\rangle=\widehat w\rho$ for some element $\widehat w$ in the
affine Weyl group $\widehat W$. Here, the affine Weyl group is the group
of affine isometries of $\cartan_\RR^*$ generated by the finite Weyl
group $W$---which is itself
generated by the simple reflections $s_1,\dots,s_l$ along the
simple roots $\alpha_1,\dots,\alpha_l$, that is,
$$s_i:\lambda\longmapsto\lambda-
\frac{2\,(\lambda|\alpha_i)}{(\alpha_i|\alpha_i)}\,\alpha_i=\lambda-\langle
\lambda,\alpha_i^\vee\rangle\,\alpha_i$$
---and, in addition, the affine reflection
$$s_0:\lambda\mapsto\lambda-\Bigl(\frac{2\,(\lambda|\theta)}{(\theta|\theta)}
-g\Bigr)\theta
=\lambda-\bigl(\langle\lambda,\theta^\vee\rangle-g\bigr)\theta
=s_\theta\lambda+g\theta.$$
Here, $\langle\phantom{\alpha},\phantom{\alpha}\rangle:\cartan_\RR^*\times
\cartan_\RR\to\RR$ is the natural pairing, and $\alpha_1^\vee,\dots,\alpha_l^\vee,
\theta^\vee$ are the coroots corresponding to $\alpha_1,\dots,\alpha_l,\theta$.
More generally, for any root $\varphi\in\Phi$ the corresponding coroot
$\varphi^\vee\in\cartan_\RR$ is defined by
$$\langle\lambda,\varphi^\vee\rangle
=\frac{2\,(\lambda|\varphi)}{(\varphi|\varphi)}\quad\forall\,\lambda\in\cartan_\RR^*.$$

The affine Weyl group $\widehat W$ is a Coxeter group with Coxeter generators
$s_0,\dots,s_l$.
Let $\ell:\widehat W\to\ZZ_{\geqslant0}$ be the usual length function,
$\ell(\widehat w)=r$ if $\widehat w=s_{i_1}\dots s_{i_r}$ with
${i_1,\dots,i_r}\in\{{0,\dots,l}\}$ and $r$ minimal. Similarly, denoting
again by $\ell:W\to\ZZ_{\geqslant0}$ the length function of the parabolic
subgroup $W\subseteq\widehat W$, one knows that it coincides with the restriction
of the length function of $\widehat W$.

\begin{myremark}
The definition of the affine Weyl group is not exactly the standard but
a scaled one, to the effect that $s_0\rho=\rho+\theta$.
One has the well-known decomposition $\widehat W\cong gM\rtimes W$ of
$\widehat W$ as a semidirect product of $W$ acting on the normal subgroup
$gM$, the lattice spanned by the long roots and dilated by the factor $g$,
in the obvious way. Each element $\mu\in gM$ acts as the translation
$\lambda\mapsto\lambda+\mu$.

There is of course also the linear version of the affine Weyl group acting
on $\cartan_\RR^*\oplus\RR\delta\oplus\RR\Lambda_0$ as in the book~\cite{Ka}.
One extends the inner product in $\cartan_\RR^*$ to a nondegenerate symmetric
bilinear form, again denoted $(\phantom{\alpha}|\phantom{\alpha})$, by
declaring that $\delta$ and $\Lambda_0$ are isotropic vectors perpendicular
to $\cartan_\RR^*$ and such that $(\delta|\Lambda_0)=1$.
We denote by abuse of notation the reflections $s_0,\dots,s_l\in
\operatorname{O}\bigl(\cartan_\RR^*\oplus\RR\delta\oplus\RR\Lambda_0,
(\phantom{\alpha}|\phantom{\alpha})\bigr)$ given by the formula
$$s_i:\lambda\longmapsto\lambda-\frac{2\,(\lambda|\alpha_i)}{(\alpha_i|\alpha_i)}
\,\alpha_i,$$
where $\alpha_1,\dots,\alpha_l$ are the simple roots as usual but now
$\alpha_0=\delta-\theta$. The group generated by $s_0,\dots,s_l$ will again
by abuse of notation be denoted by $\widehat W$. Each affine hyperplane
$\cartan_\RR^*\oplus\RR\delta+c\Lambda_0$ is mapped into itself by the
reflections $s_0,\dots,s_l$. The action of $\widehat W$ on $\cartan_\RR^*$ defined
previously comes from the action of $\widehat W$ on the subquotient
$\cartan_\RR^*\oplus\RR\delta+\frac12\Lambda_0\pmod{\RR\delta}$
if one identifies this subquotient with $\cartan_\RR^*$ in the evident way.
\end{myremark}

The fundamental weights $\varpi_1,\dots,\varpi_l\in\cartan_\RR^*$ are the basis
dual to the basis $\alpha_1^\vee,\dots,\alpha_l^\vee$ of $\cartan_\RR$.
Recall also that $\rho=\sum\limits_{i=1}^l\varpi_i$.
The next definition is slightly
non-standard: define $\stackrel{\vee}{\varpi}_1,\dots,\stackrel{\vee}{\varpi}_l
\in\cartan_\RR^*$ by
$\bigl(\stackrel{\vee}{\varpi}_i\bigm|\!\!\alpha_j\bigr)=\delta_{ij}\,\frac12$, that is,
$\stackrel{\vee}{\varpi}_i=\frac{1}{\Vert\alpha_i\Vert^2}\,\varpi_i$.

The map $\widehat W\to\cartan_\RR^*$, $\widehat w\mapsto\widehat w\rho$ is
injective. Its image will be termed the set of
\textbf{\boldmath{$\rho$}-points}. Let
$$A=\bigl\{\lambda\in\cartan_\RR^*\bigm|(\lambda|\alpha)\geqslant0\mbox{ for
all $\alpha\in\Pi$ and }\langle\lambda,\theta^\vee\rangle\leqslant g\bigr\}$$
be the (closed) fundamental alcove, which is a fundamental domain for
$\widehat W$ acting on $\cartan_\RR^*$. The fundamental alcove $A$ is the
simplex whose vertices are
$0,\frac{\stackrel{\vee}{\varpi}_1}{n_1},\dots,\frac{\stackrel{\vee}{\varpi}_l}{n_l}$
where $n_1,\dots,n_l$ are the marks, i.\,e., the (positive integer) coefficients in
$\theta=\sum\limits_{i=1}^ln_i\,\alpha_i$.

The cone with apex~$0$ spanned by $A$ is the dominant chamber
$$C=\bigl\{\lambda\in\cartan_\RR^*\bigm|(\lambda|\alpha)\geqslant0\mbox{ for
all $\alpha\in\Pi$}\bigr\}.$$
It is a fundamental domain for the finite Weyl group~$W$.

The $\widehat W$-translates of the fundamental alcove are called alcoves.
For each $i={0,\dots,l}$ one has the (affine if $i=0$) hyperplane
$$H_i:=\bigl\{\lambda\in\cartan_\RR^*\bigm|s_i\lambda=\lambda\bigr\},$$
and their $\widehat W$-translates are termed walls.

The $\rho$-points are precisely the integral weights in the interior of an
alcove. So there are the natural bijections
\begin{eqnarray*}
\widehat W\stackrel{\cong}{\longleftrightarrow}&
\{\mbox{alcoves}\}&\stackrel{\cong}{\longleftrightarrow}
\{\mbox{$\rho$-points}\}\\
\widehat w\,\longleftrightarrow&\widehat wA&\longleftrightarrow\,\widehat w\rho.
\end{eqnarray*}

We have already mentioned above that $\rho+\theta$ is the $\rho$-point of the
alcove $s_0A$. The $\rho$-points of the other neighbours $s_1A,\dots,s_lA$
of the fundamental alcove are $\rho-\alpha_1,\dots,\rho-\alpha_l$.
This follows because $s_i$ (for $i={1,\dots,l}$) permutes all positive
roots other than $\alpha_i$ and $s_i(\alpha_i)=-\alpha_i$.

The following picture shows part of the tessellation of the plane
by alcoves for type $\sfG_2$. The shaded region marks the fundamental alcove $A$.
The boundaries of the four alcoves in $2A=\{2\lambda\,|\,
\lambda\in A\}$ are drawn in solid lines.

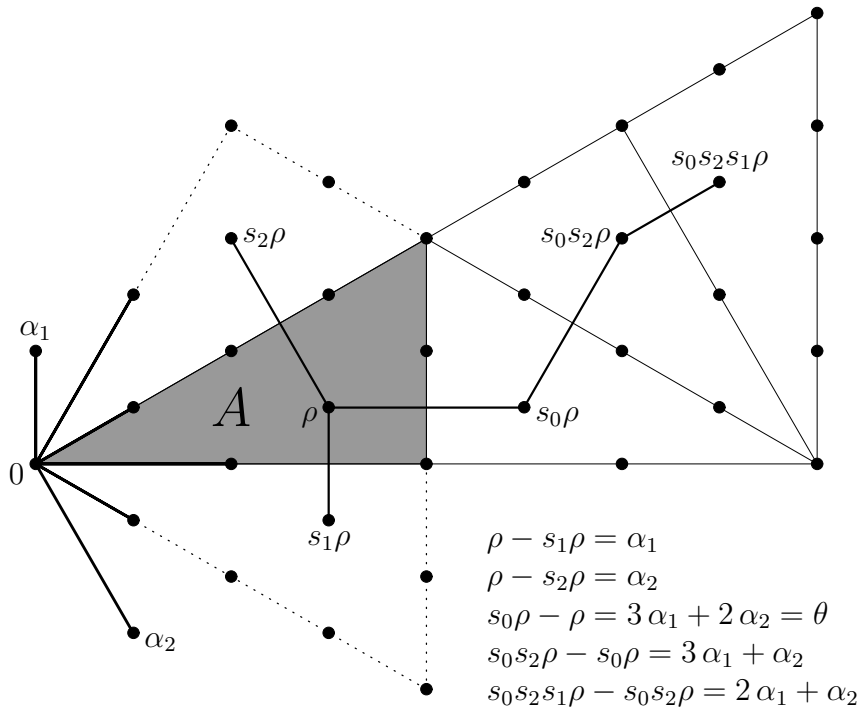
\setlength{\unitlength}{0.015mm}
\begin{picture}(7000,6700)(0,-2500)
\texture{55888888 88555555 5522a222 a2555555 55888888 88555555 552a2a2a 2a555555 
        55888888 88555555 55a222a2 22555555 55888888 88555555 552a2a2a 2a555555 
        55888888 88555555 5522a222 a2555555 55888888 88555555 552a2a2a 2a555555 
        55888888 88555555 55a222a2 22555555 55888888 88555555 552a2a2a 2a555555 }
\shade\path(0,0)(3464,0)(3464,2000)(0,0)
\Thicklines
\path(0,0)(0,1000)(0,0)(866,1500)(0,0)(866,500)(0,0)(1732,0)(0,0)(866,-500)
 (0,0)(866,-1500)
\thinlines
\path(0,0)(6928,0)(6928,4000)(0,0)
\path(5196,3000)(6928,0)(3464,2000)(3464,0)
\dottedline{80}(0,0)(1732,3000)(3464,2000)
\dottedline{80}(3464,0)(3464,-2000)(0,0)
\multiput(0,0)(0,1000){2}{\circle*{100}}
\multiput(866,-1500)(0,1000){4}{\circle*{100}}
\multiput(1732,-1000)(0,1000){5}{\circle*{100}}
\multiput(2598,-1500)(0,1000){5}{\circle*{100}}
\multiput(3464,-2000)(0,1000){5}{\circle*{100}}
\multiput(4330,500)(0,1000){3}{\circle*{100}}
\multiput(5196,0)(0,1000){4}{\circle*{100}}
\multiput(6062,500)(0,1000){4}{\circle*{100}}
\multiput(6928,0)(0,1000){5}{\circle*{100}}
\thicklines
\path(2598,500)(4330,500)(5196,2000)(6062,2500)
\path(1732,2000)(2598,500)(2598,-500)
\put(2498,500){\makebox(0,0)[tr]{$\rho$}}
\put(4430,500){\makebox(0,0)[tl]{$s_0\rho$}}
\put(5096,2000){\makebox(0,0)[r]{$s_0s_2\rho$}}
\put(6062,2600){\makebox(0,0)[b]{$s_0s_2s_1\rho$}}
\put(1832,2000){\makebox(0,0)[l]{$s_2\rho$}}
\put(2598,-600){\makebox(0,0)[t]{$s_1\rho$}}
\put(966,-1500){\makebox(0,0)[tl]{$\alpha_2$}}
\put(0,1100){\makebox(0,0)[b]{$\alpha_1$}}
\put(1732,500){\makebox(0,0){\LARGE$A$}}
\put(-100,0){\makebox(0,0)[tr]{$0$}}
\put(4000,-500){\makebox(0,0)[tl]{$\begin{array}{@{}l}
\rho-s_1\rho=\alpha_1\\
\rho-s_2\rho=\alpha_2\\
s_0\rho-\rho=3\,\alpha_1+2\,\alpha_2=\theta\\
s_0s_2\rho-s_0\rho=3\,\alpha_1+\alpha_2\\
s_0s_2s_1\rho-s_0s_2\rho=2\,\alpha_1+\alpha_2
\end{array}$}}
\label{pictureG2}
\end{picture}
\setlength{\unitlength}{3.5mm}

Passing from the sum $\rho+\langle\Psi\rangle=\widehat w\rho$ to the
corresponding alcove $\widehat wA$, one can rephrase Kostant's theorem
by  saying that the alcoves belonging to an
abelian ideal are exactly those lying in $2A$ (so there are
$2^l$ of them).

Let us also recall that there is a close connexion between reduced expressions
for elements  $\widehat w\in\widehat W$ and minimal galleries going from $A$ to
$\widehat wA$. In fact, in general, one has the bijection
\begin{align*}
\{\mbox{words in $s_0,\dots,s_l$}\}&\stackrel{\cong}{\longleftrightarrow}
\{\mbox{non-stuttering galleries beginning at $A$}\}\\
s_{i_1}s_{i_2}\dots s_{i_r}&\longleftrightarrow
A,\,s_{i_1}A,\,s_{i_1}s_{i_2}A,\,\dots,\,s_{i_1}s_{i_2}\dots s_{i_r}A,
\end{align*}
and reduced words correspond to minimal galleries. The length $\ell(\widehat w)$
is the number of walls for which $A$ and $\widehat wA$ lie on opposite
sides of the wall.

\section{Explicit description of the abelian ideals}
Our approach to describing the abelian ideals
hinges on the observation that for each abelian ideal
$\ab\unlhd\borel$ its subspace $\ab^\anc$ spanned by the root subspaces
for the roots
that are not perpendicular to the highest root $\theta$ is again an abelian
ideal in $\borel$. This is the content of the next proposition.

\begin{myproposition}
Let $\ab\unlhd\borel$ be an abelian ideal with $\Phi(\ab)\subseteq
\Phi_+$ the corresponding set of positive roots. Then
$$\Phi^\anc(\ab):=\bigl\{\varphi\in\Phi(\ab)\bigm|(\varphi|\theta)>0\bigr\}$$
is also the set of positive roots for an abelian ideal $\ab^\anc
\unlhd\borel$\,:
$$\Phi^\anc(\ab)=\Phi(\ab^\anc)\ \bigl({}=\Phi^\anc(\ab^\anc)
\bigr).$$
\end{myproposition}
\begin{proof}
The abelianess is clear: $\Phi^\anc(\ab)\stackrel{.}{+}\Phi^\anc(\ab)
\subseteq\Phi(\ab)\stackrel{.}{+}\Phi(\ab)=\varnothing$.\\
That the ideal property holds is also easy to show. For that we must see that
$\Phi^\anc(\ab)\stackrel{.}{+}\Phi_+\subseteq\Phi^\anc(\ab)$.
Let $\varphi\in\Phi^\anc(\ab)$ and $\varphi'\in\Phi_+$. Of the four a priori
possibilities
(1)~$\varphi+\varphi'\notin\Phi_+$,
(2)~$\varphi+\varphi'\in\Phi_+-\Phi(\ab)$,
(3)~$\varphi+\varphi'\in\Phi(\ab)-\Phi^\anc(\ab)$,
(4)~$\varphi+\varphi'\in\Phi^\anc(\ab)$,
we must exclude the cases (2) and (3).
That (2) is impossible follows from $\Phi^\anc(\ab)\subseteq\Phi(\ab)$
and the fact that $\ab$ is an ideal. Case~(3) cannot occur because
$(\varphi'|\theta)\geqslant0$ and hence $({\varphi+\varphi'}|\theta)>0$
by the definition of $\Phi^\anc(\ab)$.
\end{proof}

Now the problem of describing the abelian ideals in $\borel$ decomposes into
two problems according to the disjoint union decomposition
$$\bigl\{\mbox{abelian ideals $\ab\unlhd\borel$}\bigr\}
=\coprod_{\ab'}\bigl\{\mbox{abelian ideals $\ab\unlhd\borel$ with
$\ab^\anc=\ab'$}\bigr\}.$$
The two tasks are
\begin{itemize}
\item[](1) \ describe the index set $\bigl\{\ab'\bigm|\mbox{$\ab'\unlhd\borel$
abelian ideal with $\ab'=\ab'^\anc$}\bigr\}$;
\item[](2) \ for each $\ab'=\ab'^\anc$ describe the set of abelian ideals $\ab$
with $\ab'=\ab^\anc$.
\end{itemize}

We will first deal with task~(1) and show that there is a canonical one-to-one
correspondence
$$\bigl\{\ab^\anc\bigm|0\neq\ab\unlhd\borel
\mbox{ abelian ideal}\bigr\}\stackrel{\cong}{\longleftrightarrow}
\Phi_+^{\textup{(long)}}$$
(see Theorem~\ref{root} below).
This will then extend and give an a priori explanation for the observation that the
maximal abelian ideals are in canonical one-to-one correspondence with the
simple long roots, as it was recorded in~\cite{PR}.
We need some preparation.

Recall that the height $\operatorname{ht}(\varphi)$ of a root $\varphi=
\sum\limits_{i=1}^lc_i\,\alpha_i$ is defined to be $\operatorname{ht}(\varphi)=
\sum\limits_{i=1}^lc_i$, its coefficient sum with respect to the basis of simple
roots. The simple roots are those having height $1$, and the highest root
$\theta$ is the root whose height is maximal, namely,
$\operatorname{ht}(\theta)=h-1$, which is $1$ less than the Coxeter number.

For our purpose a modification of the height will be important. We define the
affine functional $L:\cartan_\RR^*\to\RR$ by
\begin{equation}\label{defL}
L(\varphi):=\frac{2\,(\theta-\varphi|\rho)}{(\theta|\theta)}\,.
\end{equation}
Whereas $\operatorname{ht}(\varphi)>0$ for $\varphi\in\Phi_+$ and
$\operatorname{ht}(\varphi)<0$ for $\varphi\in\Phi_-$, we have $L(\varphi)
\geqslant0$ for all $\varphi\in\Phi$; more precisely, $L(\theta)=0$ and
$L(\varphi)>0$ for all $\varphi\in\Phi-\{\theta\}$.
A second modification concerns the root lengths. Let us write again $\varphi=
\sum\limits_{i=1}^lc_i\,\alpha_i$. Then
\begin{align*}
L(\varphi)&=\frac{2\,(\theta-\varphi|\rho)}{(\theta|\theta)}=g-1-
\frac{2}{(\theta|\theta)}\textstyle\Bigl(\sum\limits_{i=1}^lc_i\,\alpha_i\Bigm|
\sum\limits_{j=1}^l\varpi_j\Bigr)\\
&=g-1-\sum_{i=1}^lc_i\frac{(\alpha_i|\alpha_i)}{(\theta|\theta)}
\end{align*}
because $(\alpha_i|\varpi_j)=\frac12(\alpha_i|\alpha_i)\,\delta_{ij}$.
Note also that $L(\varphi)=g-1-\langle\rho,\varphi^\vee\rangle\in
\ZZ_{\geqslant0}$ if $\varphi$ is a long root.
In particular, for $\varphi\in\Phi_+^{\textup{(long)}}$ we have
$L(\varphi)=g-2$ if and only if $\varphi$ is a long simple root.
The affine functional $L$ shows its importance in the following proposition.

\begin{myproposition}\label{L}
For each positive long root $\varphi\in\Phi_+^{\textup{(long)}}$ there is
a unique Weyl group element $w\in W$ of length $\ell(w)=L(\varphi)$
such that $w\varphi=\theta$ is the highest root. Moreover, $w'\varphi\neq\theta$
for all $w'\in W$ with $\ell(w')<L(\varphi)$.
\end{myproposition}
\begin{proof}
We first show the minimality that is expressed in the second sentence.
Let $s_i$ be a simple reflection with $s_i\varphi$ positive, too, hence
$\varphi\neq\pm\alpha_i$.
We compute
\begin{align}
L(\varphi)-L(s_i\varphi)&=\frac{2\,(\theta-\varphi|\rho)}{(\theta|\theta)}
-\frac{2\,(\theta-s_i\varphi|\rho)}{(\theta|\theta)}
\tag{by definition (\ref{defL})}\\
&=-\frac{2\,(\varphi|\rho)}{(\theta|\theta)}+
\frac{2\,(\varphi|s_i\rho)}{(\theta|\theta)}
\tag{by orthogonality}\\
&=-\frac{2\,(\varphi|\rho)}{(\theta|\theta)}+
\frac{2\,(\varphi|\rho-\alpha_i)}{(\theta|\theta)}
\tag{$s_i$ is a simple reflection}\\
&=-\frac{2\,(\varphi|\alpha_i)}{(\theta|\theta)}
=-\frac{2\,(\varphi|\alpha_i)}{(\varphi|\varphi)}
\tag{$\varphi$ is a long root}\\
\tag{$\varphi\neq\pm\alpha_i$ is a long root}
&=-\langle\alpha_i,\varphi^\vee\rangle\in\{0,\pm1\}.
\end{align}
It follows that $L(\varphi)=L(\varphi)-L(\theta)\leqslant\ell(w)$ if $w\varphi
=\theta$.

The same calculation shows that given $\varphi\in\Phi_+^{\textup{(long)}}
-\{\theta\}$,
there exists a simple reflection $s_i$ such that $L(\varphi)
-L(s_i\varphi)=1$. Otherwise, by the previous computation, we would
get $\langle\alpha_i,\varphi^\vee\rangle\geqslant0$ for all $i={1,\dots,l}$,
so that $\varphi$ would lie in the dominant chamber. This is absurd
because $\theta$ is the only long dominant root and $\varphi\neq\theta$
by assumption.
Hence there is a sequence $s_{i_1},\dots,s_{i_{L(\varphi)}}$
of simple reflections such that $L(s_{i_p}s_{i_{p-1}}\dots s_{i_1}\varphi)
=L(\varphi)-p$ (for $p={0,\dots,L(\varphi)}$). In particular,
$s_{i_{L(\varphi)}}s_{i_{L(\varphi)-1}}\dots s_{i_1}\varphi=\theta$.

Uniqueness follows from the uniqueness of coset representatives of
minimal length for standard parabolic subgroups. In fact, let $W_{\perp\theta}$
be the standard parabolic
subgroup generated by the simple reflections that fix the highest root~$\theta$.
Its Coxeter-Dynkin graph is the subgraph of the Coxeter-Dynkin graph
of $W$ induced by those nodes that are not adjacent to the affine node.
The quotient in question is then the set of right cosets $W_{\perp\theta}
\backslash W$.
\end{proof}

The following table compiles for each long simple root $\alpha_i$
the Weyl group element $w$ with $\ell(w)=g-2$ and such that
$w\alpha_i=\theta$. The labeling coincides with the labeling
in the table following page~\pageref{longtable}.

{\small
\renewcommand{\arraystretch}{1.2}%
\begin{longtable}{|@{\,}l@{\,}|l|l|}\hline
$\sfX$&$i$&$w$ such that $w\alpha_i=\theta$\\\hline\hline
$\sfA_l$
&$i$&$s_1\stackrel{\nearrow}{\dots}s_{i-1}\,s_l\stackrel{\searrow}{\dots}s_{i+1}$\\\hline
$\sfC_l$
&$l$&$s_1\stackrel{\nearrow}{\dots}s_{l-1}$\\\hline
$\sfB_l$
&$i$&$s_2\stackrel{\nearrow}{\dots}s_l\,
s_1\stackrel{\nearrow}{\dots}s_{i-1}\,
s_{l-1}\stackrel{\searrow}{\dots}s_{i+1}\quad(i={1,\dots,l-1})$\\\hline
$\sfD_l$
&$i$&$s_2\stackrel{\nearrow}{\dots}s_{l-2}\,
s_1\stackrel{\nearrow}{\dots}s_{i-1}\,
s_l\stackrel{\searrow}{\dots}s_{i+1}\quad(i={1,\dots,l-2})$\\*
&$i$&$s_2\stackrel{\nearrow}{\dots}s_{l-2}\,
s_1\stackrel{\nearrow}{\dots}s_{l-3}\,
s_{2l-i-1}\,s_{l-2}\quad(i={l-1,l})$\\\hline
$\sfE_6$
&$1$&$s_1\,s_2\,s_3\,s_4\,s_2\,s_5\,s_3\,s_6\,s_4\,s_2$\\*
&$2$&$s_1\,s_2\,s_3\,s_4\,s_2\,s_1\,s_5\,s_3\,s_6\,s_4$\\*
&$3$&$s_1\,s_2\,s_3\,s_4\,s_2\,s_1\,s_5\,s_6\,s_4\,s_2$\\*
&$4$&$s_1\,s_2\,s_3\,s_4\,s_2\,s_1\,s_5\,s_3\,s_2\,s_6$\\*
&$5$&$s_1\,s_2\,s_3\,s_4\,s_2\,s_1\,s_6\,s_4\,s_2\,s_3$\\*
&$6$&$s_1\,s_2\,s_3\,s_4\,s_2\,s_1\,s_5\,s_3\,s_2\,s_4$\\\hline
$\sfE_7$
&$1$&$s_1\,s_2\,s_3\,s_4\,s_5\,s_3\,s_2\,s_6\,s_4\,s_3\,s_5\,s_7\,s_6\,s_4\,
s_3\,s_2$\\*
&$2$&$s_1\,s_2\,s_3\,s_4\,s_5\,s_3\,s_2\,s_1\,s_6\,s_4\,s_3\,s_5\,s_7\,s_6\,
s_4\,s_3$\\*
&$3$&$s_1\,s_2\,s_3\,s_4\,s_5\,s_3\,s_2\,s_1\,s_6\,s_4\,s_3\,s_2\,s_5\,s_7\,
s_6\,s_4$\\*
&$4$&$s_1\,s_2\,s_3\,s_4\,s_5\,s_3\,s_2\,s_1\,s_6\,s_4\,s_3\,s_2\,s_5\,s_3\,
s_7\,s_6$\\*
&$5$&$s_1\,s_2\,s_3\,s_4\,s_5\,s_3\,s_2\,s_1\,s_6\,s_4\,s_3\,s_2\,s_7\,s_6\,
s_4\,s_3$\\*
&$6$&$s_1\,s_2\,s_3\,s_4\,s_5\,s_3\,s_2\,s_1\,s_6\,s_4\,s_3\,s_2\,s_5\,s_3\,
s_4\,s_7$\\*
&$7$&$s_1\,s_2\,s_3\,s_4\,s_5\,s_3\,s_2\,s_1\,s_6\,s_4\,s_3\,s_2\,s_5\,s_3\,
s_4\,s_6$\\\hline
$\sfE_8$&$1$&$
s_1\,s_2\,s_3\,s_4\,s_5\,s_6\,s_7\,s_5\,s_4\,s_3\,s_2\,s_8\,s_6\,s_5\,
s_4\,s_3\,s_7\,s_5\,s_4\,s_6\,s_5\,s_7\,s_8\,s_6\,s_5\,s_4\,s_3\,s_2$\\*
&$2$&$
s_1\,s_2\,s_3\,s_4\,s_5\,s_6\,s_7\,s_5\,s_4\,s_3\,s_2\,s_1\,s_8\,s_6\,
s_5\,s_4\,s_3\,s_7\,s_5\,s_4\,s_6\,s_5\,s_7\,s_8\,s_6\,s_5\,s_4\,s_3$\\*
&$3$&$
s_1\,s_2\,s_3\,s_4\,s_5\,s_6\,s_7\,s_5\,s_4\,s_3\,s_2\,s_1\,s_8\,s_6\,
s_5\,s_4\,s_3\,s_2\,s_7\,s_5\,s_4\,s_6\,s_5\,s_7\,s_8\,s_6\,s_5\,s_4$\\*
&$4$&$
s_1\,s_2\,s_3\,s_4\,s_5\,s_6\,s_7\,s_5\,s_4\,s_3\,s_2\,s_1\,s_8\,s_6\,
s_5\,s_4\,s_3\,s_2\,s_7\,s_5\,s_4\,s_3\,s_6\,s_5\,s_7\,s_8\,s_6\,s_5$\\*
&$5$&$
s_1\,s_2\,s_3\,s_4\,s_5\,s_6\,s_7\,s_5\,s_4\,s_3\,s_2\,s_1\,s_8\,s_6\,
s_5\,s_4\,s_3\,s_2\,s_7\,s_5\,s_4\,s_3\,s_6\,s_5\,s_4\,s_7\,s_8\,s_6$\\*
&$6$&$
s_1\,s_2\,s_3\,s_4\,s_5\,s_6\,s_7\,s_5\,s_4\,s_3\,s_2\,s_1\,s_8\,s_6\,
s_5\,s_4\,s_3\,s_2\,s_7\,s_5\,s_4\,s_3\,s_6\,s_5\,s_4\,s_7\,s_5\,s_8$\\*
&$7$&$
s_1\,s_2\,s_3\,s_4\,s_5\,s_6\,s_7\,s_5\,s_4\,s_3\,s_2\,s_1\,s_8\,s_6\,
s_5\,s_4\,s_3\,s_2\,s_7\,s_5\,s_4\,s_3\,s_6\,s_5\,s_4\,s_8\,s_6\,s_5$\\*
&$8$&$
s_1\,s_2\,s_3\,s_4\,s_5\,s_6\,s_7\,s_5\,s_4\,s_3\,s_2\,s_1\,s_8\,s_6\,
s_5\,s_4\,s_3\,s_2\,s_7\,s_5\,s_4\,s_3\,s_6\,s_5\,s_4\,s_7\,s_5\,s_6$\\\hline
$\sfF_4$
&$1$&$s_1\,s_2\,s_3\,s_2\,s_4\,s_3\,s_2$\\
&$2$&$s_1\,s_2\,s_3\,s_2\,s_1\,s_4\,s_3$\\\hline
$\sfG_2$
&$2$&$s_2\,s_1$\\\hline
\multicolumn{3}{p{123mm}}{
{\footnotesize
The notations $\stackrel{\nearrow}{\dots}$ and $\stackrel{\searrow}{\dots}$
mean that one has to interpolate by increasing and decreasing indices,
respectively, with the obvious conventions understood, e.\,g.,
$s_1\stackrel{\nearrow}{\dots}s_4=s_1\,s_2\,s_3\,s_4$ and also
$s_1\stackrel{\nearrow}{\dots}s_1=s_1$ and
$s_1\stackrel{\nearrow}{\dots}s_0=1$ (empty index set).
}}
\end{longtable}
\renewcommand{\arraystretch}{1}%
}

\begin{myremark}
Note that if $s_{j_1}\dots s_{j_{g-2}}\alpha_i=\theta$, then for each $r=1,\dots,g-2$,
$$\theta-\sum_{k=1}^r\frac{\Vert\theta\Vert^2}{\Vert\alpha_{j_k}\Vert^2}\,
\alpha_{j_k}$$
is a positive root (and for $r=g-2$ equals $\alpha_i$).

Surely, one can count the number of reduced decompositions, and one gets generalizations
of binomial coefficients. We do not elaborate one this point here which seems to
be known anyway.
\end{myremark}

\begin{myremark}
The lengths of the Weyl group elements that occurred in Proposition~\ref{L}
have the following description. Let
$$q(t):=\sum_{\substack{\varphi\in\Phi^{\textup{(long)}}\\
\makebox[0pt]{\footnotesize$\left\{\begin{array}{@{}l@{}}
w\varphi=\theta\textrm{ with}\\
\ell(w)\textrm{ minimal}
\end{array}\right.$}}}t^{\ell(w)}
=\sum_{\varphi\in\Phi_+^{\textup{(long)}}}t^{L(\varphi)}.$$
Since $s_i(\alpha_i)=-\alpha_i$, we have the corresponding sum
$$\sum_{\substack{\varphi\in\Phi^{\textup{(long)}}\\
\makebox[0pt]{\footnotesize$\left\{\begin{array}{@{}l@{}}
w\varphi=\theta\textrm{ with}\\
\ell(w)\textrm{ minimal}
\end{array}\right.$}}}t^{\ell(w)}=
q(t)+t^{2g-3}q(t^{-1})$$
for all the (positive and negative) long roots, which is the Poincar\'e polynomial
for the set of minimal coset representatives for $W_{\perp\theta}\backslash W$.
The usual Poincar\'e polynomial $W(t)$ of $W$ is
$W(t)=\sum\limits_{w\in W}t^{\ell(w)}$. It can be expressed by a product formula.
Namely, if $m_1,\dots,m_l$ are the exponents\footnote{A word about the
labeling: the numbers $m_1,\dots,m_l$ are not naturally associated
to the nodes of the Coxeter-Dynkin graph.}
of $W$ (or of its type $\sfX$), then using the abbreviation
$[n]:=(1-t^n)/(1-t)$ one can write $W(t)=\prod\limits_{i=1}^l[m_i+1]$.

The Poincar\'e polynomial of $W_{\perp\theta}\backslash W$ is the quotient
$W(t)/W_{\perp\theta}(t)$
of the corresponding Poincar\'e polynomials.
The rightmost column in the next table
contains the numbers $\nu(\sfX)={}$ the number of positive
long roots in a root system of type~$\sfX$. The Poincar\'e polynomial
evaluated at $t=1$ equals $2\,\nu(\sfX)$.
\label{Poincare}
\renewcommand{\arraystretch}{2.4}%
$$\begin{array}{|c|c|c|c|c|}\hline
\sfX&\sfX_{\perp\theta}&\mbox{exponents of $\sfX$}&W(t)/W_{\perp\theta}(t)&\nu(\sfX)\\\hline
\sfA_l&\sfA_{l-2}&1,2,\dots,l&[l][l+1]&\dfrac{l(l+1)}{2}\\
\sfC_l&\sfC_{l-1}&1,3,\dots,2l-1&[2l]&l\\
\sfB_l&\sfB_{l-2}+\sfA_1&1,3,\dots,2l-1&\dfrac{[2l-2][2l]}{[2]}&l(l-1)\\
\sfD_l&\sfD_{l-2}+\sfA_1&1,3,\dots,2l-3,l-1&\dfrac{[l][2l-4][2l-2]}{[2][l-2]}&l(l-1)\\
\sfE_6&\sfA_5&1,4,5,7,8,11&\dfrac{[8][9][12]}{[3][4]}&36\\
\sfE_7&\sfD_6&1,5,7,9,11,13,17&\dfrac{[12][14][18]}{[4][6]}&63\\
\sfE_8&\sfE_7&1,7,11,13,17,19,23,29&\dfrac{[20][24][30]}{[6][10]}&120\\
\sfF_4&\sfB_3&1,5,7,11&\dfrac{[8][12]}{[4]}&12\\
\sfG_2&\sfA_1&1,5&[6]&3\\\hline
\multicolumn{5}{p{115mm}}{
{\footnotesize
The usual conventions are employed for the entries in the column marked $\sfX_{\perp\theta}$,
namely, $\sfA_{-1}=\sfA_0=\sfB_0=\varnothing$,
$\sfC_1=\sfB_1=\sfA_1$, $\sfC_2=\sfB_2$, $\sfD_2=\sfA_1+\sfA_1$.
}}
\end{array}$$
\renewcommand{\arraystretch}{1}%
\end{myremark}

\begin{myremark}
Putting $[n]_k:=(1-t^{kn})/(1-t^k)$ one can write, e.\,g.,
$\dfrac{[8][9][12]}{[3][4]}=[2]_4[3]_3[12]_1=[2]_4[9]_1[4]_3=[8]_1[3]_3[3]_4$.\par
In all cases, $W(t)/W_{\perp\theta}(t)$ is of the form
$\dfrac{[a][c][e]}{[b][d]}$ (since \mbox{$1=[1]$}). Of course, one can always take $e=h$.
\par
Here is a little numerological table for the types $\sfD_4$, $\sfE_6$, $\sfE_7$,
and $\sfE_8$.
\renewcommand{\arraystretch}{2.4}
$$\begin{array}{|c|c|c|c|}\hline
\sfX&r:=\dfrac h6=\dfrac{l+2}{10-l}&\dfrac{[a][c][e]}{[b][d]}
=\dfrac{[4r][5r-1][6r]}{[r+1][2r]}
&|\mbox{group}|=(r+1)2r\\\hline
\sfD_4&1&\dfrac{[4][4][6]}{[2][2]}&|\Dih_2|=2\cdot2\\
\sfE_6&2&\dfrac{[8][9][12]}{[3][4]}&|\operatorname{Alt}_4|=3\cdot4\\
\sfE_7&3&\dfrac{[12][14][18]}{[4][6]}&|\operatorname{Sym}_4|=4\cdot6\\
\sfE_8&5&\dfrac{[20][24][30]}{[6][10]}&|\operatorname{Alt}_5|=6\cdot10\\\hline
\end{array}$$
\renewcommand{\arraystretch}{1}%
One could extend the table above to the types $\sfA_1$ and $\sfA_2$
but without the entries for the last column.
The six types $\sfA_1$, $\sfA_2$, $\sfD_4$, $\sfE_6$, $\sfE_7$, and
$\sfE_8$ are precisely the simply laced ones in Deligne's family
(see~\cite{De} and follow-up papers by various authors%
).
Further numerology pertaining to the types $\sfE_6$, $\sfE_7$, $\sfE_8$
can be found in the paper of Arnold \cite{Ar} about trinities.
\end{myremark}

For $w\in W$ one defines $\Phi_w:=\Phi_+\cap w\Phi_-$, the set of positive
roots which are of the form $w\varphi$ for a negative root $\varphi$.
The following fundamental lemma is well-known.

\begin{mylemma}\label{Phi}
For a Weyl group element $w\in W$ with reduced decomposition $w=s_{i_1}\dots s_{i_k}$
($i_j\in\{1,\dots,l\}$) the set $\Phi_w$ consists of the $k$ distinct positive roots
$$\alpha_{i_1},\,s_{i_1}(\alpha_{i_2}),\,s_{i_1}s_{i_2}(\alpha_{i_3}),\,
\dots,\,s_{i_1}\dots s_{i_{k-1}}(\alpha_{i_k}).$$
\end{mylemma}
\begin{proof}
Clearly $\Phi_1=\varnothing$. Now one uses induction to show that
$\Phi_{s_iw}=s_i\Phi_w\cup\{\alpha_i\}$ if $\ell(s_iw)=\ell(w)+1$
($\alpha_i\notin\Phi_w$) using the fact that $s_i$ (for $i={1,\dots,l}$)
permutes all positive roots other than $\alpha_i$ and $s_i(-\alpha_i)=\alpha_i$.
\end{proof}

One can define $s(w):=\sum\limits_{\varphi\in\Phi_w}\varphi$.
The function $s:W\to\cartan_\RR^*$ satisfies the $1$-cocycle condition
$s(ww')=ws(w')+s(w)$, and in fact, $s(w)=\rho-w\rho$.

\begin{mylemma}\label{length}
$\begin{array}[t]{@{}l}
\ell(s_iw)=\ell(w)\pm1\,\Longleftrightarrow\,w^{-1}\alpha_i\in\Phi_\pm,\\
\ell(ws_i)=\ell(w)\pm1\,\Longleftrightarrow\,w\alpha_i\in\Phi_\pm.
\end{array}$
\end{mylemma}

Now we work in the affine Kac-Moody algebra context. For that we need the
affine root system and its partition into the sets of positive and negative
roots. The so-called imaginary roots are fixed by $\widehat W$ and play no
role here. So we disregard them. Let
$\widehat\Phi_+:=\Phi_+\cup\bigl\{\varphi+n\delta\bigm|\varphi\in\Phi,\
n\in\ZZ_{>0}\bigr\}$ and $\widehat\Phi_-:=-\widehat\Phi_+$. We further define
for $\widehat w\in\widehat W$ the set $\widehat\Phi_{\widehat w}:=
\widehat\Phi_+\cap\widehat w\widehat\Phi_-$
of cardinality $\ell(\widehat w)$. The sum of the elements of
$\widehat\Phi_{\widehat w}$ is $\widehat\rho-\widehat w\widehat\rho$ where
$\widehat\rho=\rho+\frac12\Lambda_0$. Similarly, one can extend Lemma~\ref{Phi}
and Lemma~\ref{length} to the affine Kac-Moody algebra context.

One knows from Peterson's work that $\Psi=\Phi(\ab)$ for an abelian ideal
$\ab\unlhd\borel$ if and only if the set $\{\delta-\varphi|\varphi\in\Psi\}$ is
of the form $\widehat\Phi_{\widehat w}$.
The affine point of view explains Kostant's theorem (Theorem~\ref{criterion}).
For instance if $\ab$ is a $k$-dimensional abelian ideal in $\borel$
with $\Psi:=\Phi(\ab)=\{\varphi_1,\dots,\varphi_k\}$, then we have
$\widehat\Phi_{\widehat w}=\{\delta-\varphi_1,\dots,\delta-\varphi_k\}$ for
some $\widehat w\in\widehat W$. From
\begin{align*}
k\delta-\langle\Psi\rangle&=(\delta-\varphi_1)+\dots+(\delta-\varphi_k)
=\widehat\rho-\widehat w\widehat\rho
\intertext{and since $\widehat w$ is orthogonal, we get}
0&=\Vert\widehat w\widehat\rho\Vert^2-\Vert\widehat\rho\Vert^2=
\Vert\rho+\langle\Psi\rangle-k\delta+\tfrac12\Lambda_0\Vert^2-
\Vert\rho+\tfrac12\Lambda_0\Vert^2\\
&=\Vert\rho+\langle\Psi\rangle\Vert^2-k-\Vert\rho\Vert^2,
\end{align*}
so that in fact
$$\textstyle\Bigl\Vert\rho+\sum\limits_{\varphi\in\Psi}\varphi\Bigr\Vert^2
-\Vert\rho\Vert^2=|\Psi|.$$

\begin{mytheorem}\label{root}
Let $\varphi\in\Phi_+^{\textup{(long)}}$ be a positive long root. Let
$w\in W$ be the Weyl group element such that $w\varphi=\theta$ and with
$\ell(w)=L(\varphi)$ as in Proposition~\ref{L}. Then for all
$\psi\in\Phi_w$, $\theta-\psi$ is a positive root and
$$\ab^{\varphi,\min}:=\lie_\theta\oplus\bigoplus_{\psi\in\Phi_w}
\lie_{\theta-\psi}$$
is an abelian ideal of $\borel$. The $\rho$-point of the alcove
corresponding to $\ab^{\varphi,\min}$ is $s_0w\rho$.
\end{mytheorem}

\begin{myremark}
We shall see later that each nonzero abelian ideal $\ab\unlhd\borel$ that
satisfies $\ab=\ab^\anc$ is of the form $\ab=\ab^{\varphi,\min}$ for a unique
positive long root~$\varphi$.
\end{myremark}

\begin{proof}
For the proof we use Kostant's criterion.
First let us write $w$ as a reduced decomposition $w=s_{i_{L(\varphi)}}\dots
s_{i_1}$ as in the proof of Proposition~\ref{L}. Each root $\psi\in\Phi_w$ is of
the form $\psi=s_{i_{L(\varphi)}}\dots s_{i_{q+1}}(\alpha_{i_q})$ and we compute
\begin{align*}
\frac{2\,(\psi|\theta)}{(\theta|\theta)}&=\frac{2}{(\theta|\theta)}
\bigl(s_{i_{L(\varphi)}}\dots s_{i_{q+1}}(\alpha_{i_q})\bigm|
s_{i_{L(\varphi)}}\dots s_{i_1}\varphi\bigr)\\
&=\frac{2}{(\theta|\theta)}\bigl(\alpha_{i_q}\bigm|s_{i_q}\dots s_{i_1}\varphi
\bigr)
=\frac{2}{(\theta|\theta)}\bigl(-\alpha_{i_q}\bigm|s_{i_{q-1}}\dots s_{i_1}\varphi
\bigr)\\
&=-\bigl\langle\alpha_{i_q},(s_{i_{q-1}}\dots s_{i_1}\varphi)^\vee\bigr\rangle\\
&=L(s_{i_{q-1}}\dots s_{i_1}\varphi)-L(s_{i_q}\dots s_{i_1}\varphi)=1.
\end{align*}
Hence we have $s_\theta(\psi)=\psi-\theta$ and
$\theta-\psi$ is a positive root. Now we put $\ell:=\ell(w)=L(\varphi)$
for abbreviation, so that $|\Phi_w|=\ell$.
We check that the $\ell+1$ element set $\Psi:=\{\theta\}
\cup\{\theta-\psi|\psi\in\Phi_w\}\subseteq\Phi_+$ satisfies Kostant's
criterion (Theorem~\ref{criterion}) for an abelian ideal. Using
$\sum\limits_{\psi\in\Phi_w}\psi=\rho-w\rho$ we have
\begin{align*}
\rho+\theta+\sum_{\psi\in\Phi_w}(\theta-\psi)
&=w\rho+(\ell+1)\theta
=w\rho+\bigl(L(\varphi)+1\bigr)\theta\\
&=w\rho+\bigl(g-\langle\rho,\varphi^\vee\rangle\bigr)\theta\\
&=w\rho+\bigl(g-\langle w\rho,(w\varphi)^\vee\rangle\bigr)\theta
\intertext{and because $w\varphi=\theta$ we get}
&=w\rho-\langle w\rho,\theta^\vee\rangle\theta+g\theta=s_0w\rho.
\end{align*}
This proves the assertion about the $\rho$-point.
Now we compute
\begin{align*}
\bigl\Vert s_0w\rho\bigr\Vert^2-\bigl\Vert\rho\bigr\Vert^2
&=\bigl\Vert(\ell+1)\theta+w\rho\bigr\Vert^2-\bigl\Vert\rho\bigr\Vert^2\\
&=(\ell+1)^2\Vert\theta\Vert^2+(\ell+1)\,2\,(\theta|w\rho)
\intertext{and with $\ell\,\Vert\theta\Vert^2=2\,(\theta-\varphi|\rho)$
and $w^{-1}\theta=\varphi$ the calculation continues}
&=(\ell+1)\bigl(2\,(\theta-\varphi|\rho)+\Vert\theta\Vert^2+2\,(\varphi|\rho)\bigr)\\
&=(\ell+1)\bigl(\Vert\rho+\theta\Vert^2-\Vert\rho\Vert^2\bigr)=\ell+1.
\end{align*}
This completes the proof of the theorem.
\end{proof}

\begin{myremark}
As the notation $\ab^{\varphi,\min}$ suggests there will also be abelian ideals $\ab^{\varphi,\max}$.
In fact, each nonzero abelian ideal $\ab$ satisfies
$\ab^{\varphi,\min}\subseteq\ab\subseteq\ab^{\varphi,\max}$ for some positive
long root $\varphi$ which is characterized by
$\ab^\anc=\ab^{\varphi,\min}$.
If $\varphi$ is not perpendicular to the highest root $\theta$,
then $\ab^{\varphi,\max}=\ab^{\varphi,\min}$.
\end{myremark}

Let us now look closer at the case where $\varphi\perp\theta$. Before
giving the general picture, we state a preliminary result.

\begin{myproposition}\label{min+}
Let $\varphi\in\Phi_+^{\textup{(long)}}$ and $w\in W$ be as in
Theorem~\ref{root} and suppose in addition that
$\varphi$ is perpendicular to the highest root $\theta$. Then
$$\ab^{\varphi,\min^+}:=\ab^{\varphi,\min}\oplus\lie_{w\theta}$$
is an abelian ideal of $\borel$.
\end{myproposition}
\begin{proof}
We first show that $w\theta$ is a positive root perpendicular to $\theta$.
In fact, $(w\theta|\theta)=(w\theta|w\varphi)=(\theta|\varphi)=0$. Hence
$w\theta$ is a long root spanned by the simple roots $\alpha_i$ that are
perpendicular to $\theta$. (To make this assertion clear, let us write
$w\theta=\sum\limits_{i=1}^la_i\,\alpha_i$. Here the coefficients $a_i$
are either all nonnegative or all nonpositive. Now we take the inner product
with $\theta$ and use $(\alpha_i|\theta)\geqslant0$ because $\theta$
lies in the dominant chamber and $\alpha_i$ is a positive root.)
For each such root $\alpha_i\perp\theta$ we have
$s_iw\varphi=s_i\theta=\theta$. Hence $\ell(s_iw)\geqslant\ell(w)$ by the
minimality of $\ell(w)$. Lemma~\ref{length} shows that $w^{-1}\alpha_i\in
\Phi_+$. Since $\theta$ lies in the dominant chamber, we get
$0\leqslant(w^{-1}\alpha_i|\theta)=(\alpha_i|w\theta)$ for all simple
roots $\alpha_i\perp\theta$.
This means that $w\theta$ lies in the dominant chamber for the root subsystem
$\Phi_{\perp\theta}$ (spanned by the simple roots $\alpha_i\perp\theta$). One
can then see that $w\theta$ is in fact the highest root of the $\varphi$-component
of $\Phi_{\perp\theta}$.

Putting again $\ell:=\ell(w)$ we define the set $\Psi$ of cardinality $\ell+2$ as
$\Psi:=\{\theta\}
\cup\{\theta-\psi|\psi\in\Phi_w\}\cup\{w\theta\}\subseteq\Phi_+$. ($w\theta$
is of course different from the elements $\theta-\psi$ because only $w\theta$
is perpendicular to $\theta$.) Now we employ Kostant's criterion as in
the proof of Theorem~\ref{root}. Looking back at that proof we see that
we must show that
$$\Vert(\ell+1)\theta+w\rho+w\theta\Vert^2-\Vert(\ell+1)\theta+w\rho\Vert^2=1.$$
This follows from $\theta\perp w\theta$ and the $W$-invariance of the inner
product together with the identity
$\Vert\rho+\theta\Vert^2-\Vert\rho\Vert^2=1$.
\end{proof}

\subsection*{Minimal coset representatives and Poincar\'e polynomials}
In this section we define for each positive root $\varphi\in\Phi_+$
the polynomial $P_\varphi(t)\in\ZZ_{\geqslant0}[t]$ by setting
$$P_\varphi(t):=\frac{\widehat W_{\perp\varphi}(t)}{W_{\perp\varphi}(t)}\,.$$
Here $\widehat W_{\perp\varphi}$ is the standard parabolic subgroup of
the affine Weyl group $\widehat W$ generated by those reflections
$s_i$ ($i=0,\dots,l$) for which $\alpha_i\perp\varphi$ (here
$\alpha_0=-\theta$). Note that $\widehat W_{\perp\varphi}$ is a finite Coxeter group.
Similarly, $W_{\perp\varphi}$ is the standard parabolic subgroup of the
finite Weyl group $W$ generated by those simple reflections $s_i$
($i=1,\dots,l$) for which $\alpha_i\perp\varphi$. 
In particular, $\widehat W_{\perp\varphi}=W_{\perp\varphi}$ if
$\varphi\not\perp\theta$. The expressions
$\widehat W_{\perp\varphi}(t)$
and $W_{\perp\varphi}(t)$ stand for the Poincar\'e polynomials of the Coxeter groups
in question, and the quotient
$\frac{\widehat W_{\perp\varphi}(t)}{W_{\perp\varphi}(t)}$ is the Poincar\'e
polynomial for the set of minimal coset representatives in
$W_{\perp\varphi}\backslash\widehat W_{\perp\varphi}$.

Let $\widehat w_\circ^\varphi$ be the longest element of $\widehat W_{\perp\varphi}$
and $w_\circ^\varphi$ the longest element of $W_{\perp\varphi}$. The set of minimal
coset representatives in $W_{\perp\varphi}\backslash\widehat W_{\perp\varphi}$ is
the interval $\bigl[1,w_\circ^\varphi\widehat w_\circ^\varphi\bigr]$ in the
Bruhat-Chevalley order (note that $(w_\circ^\varphi)^2=1$
and also $(\widehat w_\circ^\varphi)^2=1$).
In particular, the longest element in
$\bigl[1,w_\circ^\varphi\widehat w_\circ^\varphi\bigr]$ has length
$\ell(w_\circ^\varphi\widehat w_\circ^\varphi)=\ell(\widehat w_\circ^\varphi)
-\ell(w_\circ^\varphi)$.

In the following long table we show the polynomials $P_\alpha(t)$ for all simple
roots $\alpha\in\Pi$. The polynomials $P_\varphi(t)$ can be extracted from
this piece of information. This is clear for simple types different from
$\sfA_l$ because then the affine vertex of the Coxeter-Dynkin graph is
a leaf in a tree and hence $P_\varphi(t)=P_\alpha(t)$ for an appropriate
simple root
$\alpha\in\Pi$. For type $\sfA_l$ we can reduce to the case of a simple root by
looking at $\sfA_k$ for appropriate $k$, namely,
$P^{\sfA_l}_{\alpha_i+\dots+\alpha_{i+j}}(t)=P^{\sfA_{l-j}}_{\alpha_i}(t)$.

\begin{mydefinition}
For a nonnegative integer $n$ let us recall the definition of the polynomial
$$[n]:=\frac{1-t^n}{1-t}\in\ZZ_{\geqslant0}[t].$$
Moreover, we define the factorials
$$[n]!:=\prod_{i=1}^n\,[i]$$
and their relatives
$$[2n]!!:=\prod_{i=1}^n\,[2i]\quad\mbox{and}\quad
[2n+1]!!:=\prod_{i=0}^n\,[2i+1].$$
Of course, $[0]!=[0]!!=1$.
\end{mydefinition}

The following table shows the polynomials $P_{\alpha_i}(t)
=\widehat W_{\perp\alpha_i}(t)/W_{\perp\alpha_i}(t)$
and the minimal coset representatives for $W_{\perp\alpha_i}\backslash
\widehat W_{\perp\alpha_i}$ (the latter for the classical series in the rank~$5$
case).
Above or beneath each node marked by the simple reflection $s_i$ we have depicted
along with the polynomial $P_{\alpha_i}(t)$ the Hasse graph of the Bruhat-Chevalley poset
of $W_{\perp\alpha_i}\backslash\widehat W_{\perp\alpha_i}$. To read a
minimal coset representative we have to start at the lower node and read upwards along
the edges. E.\,g., the minimal coset representatives for
$W_{\perp\alpha_3}\backslash\widehat W_{\perp\alpha_3}$ for type $\sfA_5$
are $1$, $s_0$, $s_0s_1$, $s_0s_5$, $s_0s_1s_5=s_0s_5s_1$,
$s_0s_1s_5s_0=s_0s_5s_1s_0$.

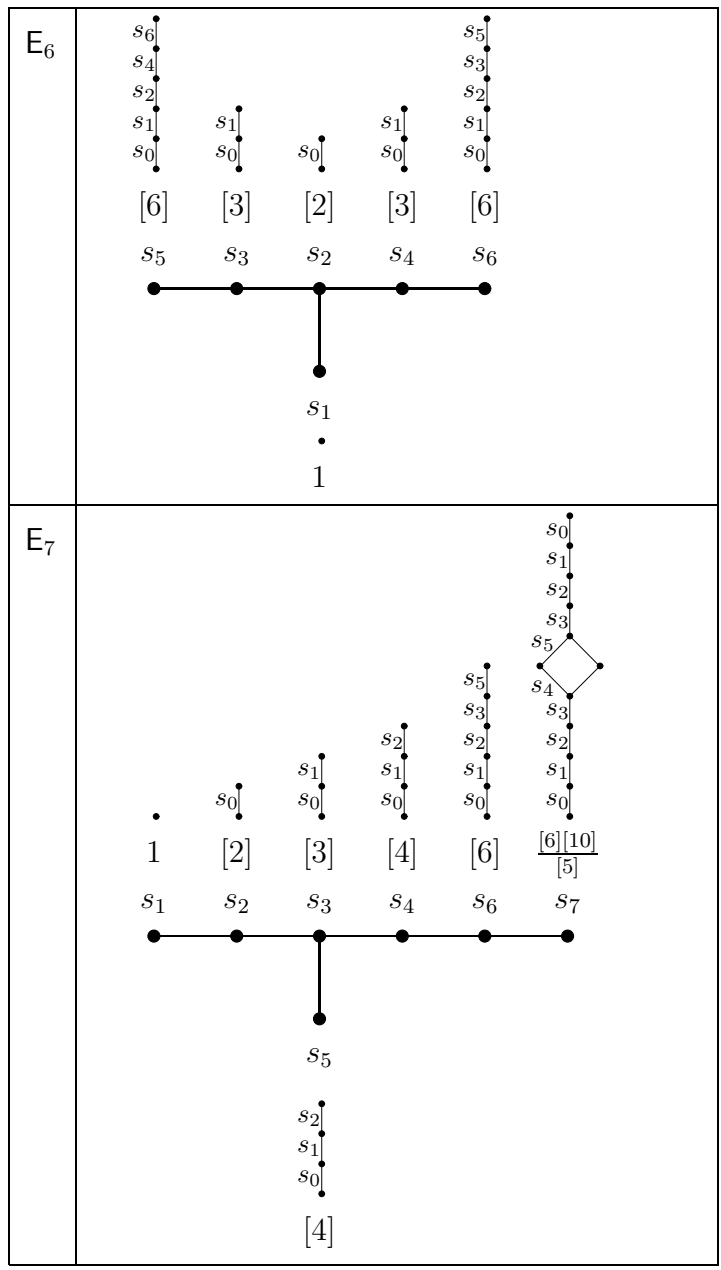
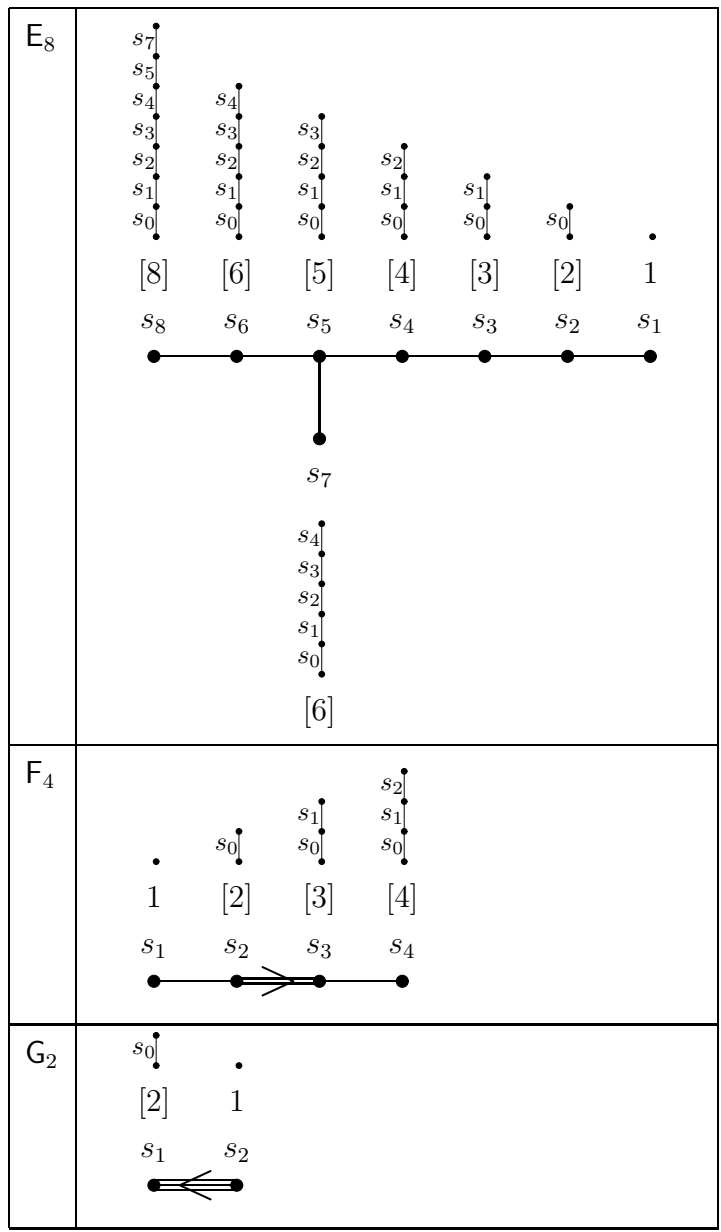
\label{longtable}
\begin{longtable}{|l|l|}\hline
$\sfA_l$&\setlength{\unitlength}{11mm}\Agenl{$\alpha_1$}{$\alpha_2$}{$\alpha_{l-1}$}{$\alpha_l$}\\*
&$P_{\alpha_i}(t)=\dfrac{[l-1]!}{[i-1]!\,[l-i]!}\quad(i=1,\dots,l)$\\*
\raisebox{30mm}{$\sfA_5$}&\raisebox{0pt}[39mm][2mm]{\setlength{\unitlength}{11mm}
\Afuenf{\setcounter{linearcounter}{0}\begin{tabular}{c}\linear{}{}{}{}{}{}{}\\[1mm]\phantom{[}1\phantom{]}\\[1mm]$s_1$\end{tabular}}%
{\setcounter{linearcounter}{3}\begin{tabular}{c}\linear{$s_0$}{$s_5$}{$s_4$}{}{}{}{}\\[1mm]{[4]}\\[1mm]$s_2$\end{tabular}}%
{\begin{tabular}{c}\Aloz\\[1mm]\phantom{[}\raisebox{0pt}[0pt][0pt]{$\frac{[3][4]}{[2]}$}\phantom{]}\\[1mm]$s_3$\end{tabular}}%
{\setcounter{linearcounter}{3}\begin{tabular}{c}\linear{$s_0$}{$s_1$}{$s_2$}{}{}{}{}\\[1mm]{[4]}\\[1mm]$s_4$\end{tabular}}%
{\setcounter{linearcounter}{0}\begin{tabular}{c}\linear{}{}{}{}{}{}{}\\[1mm]\phantom{[}1\phantom{]}\\[1mm]$s_5$\end{tabular}}
}
\\\hline
$\sfC_l$&\setlength{\unitlength}{11mm}\Cgenl{$\alpha_1$}{$\alpha_2$}{$\alpha_{l-2}$}{$\alpha_{l-1}$}{$\alpha_l$}\\*
&$P_{\alpha_i}(t)=\dfrac{[2i-2]!!}{[i-1]!}\quad(i=1,\dots,l)$\\*
\raisebox{55mm}{$\sfC_5$}&\raisebox{0pt}[60mm][2mm]{\setlength{\unitlength}{11mm}
\Cfuenf{\setcounter{linearcounter}{0}\begin{tabular}{c}\linear{}{}{}{}{}{}{}\\[1mm]\phantom{[}1\phantom{]}\\[1mm]$s_1$\end{tabular}}%
{\setcounter{linearcounter}{1}\begin{tabular}{c}\linear{$s_0$}{}{}{}{}{}{}\\[1mm]{[2]}\\[1mm]$s_2$\end{tabular}}%
{\setcounter{linearcounter}{3}\begin{tabular}{c}\linear{$s_0$}{$s_1$}{$s_0$}{}{}{}{}\\[1mm]{[4]}\\[1mm]$s_3$\end{tabular}}%
{\begin{tabular}{c}\Cloz\\[1mm]\phantom{[}\raisebox{0pt}[0pt][0pt]{$\frac{[4][6]}{[3]}$}\phantom{]}\\[1mm]$s_4$\end{tabular}}
{\begin{tabular}{c}\Clong\\[1mm]\phantom{[}\raisebox{0pt}[0pt][0pt]{$\frac{[6][8]}{[3]}$}\phantom{]}\\[1mm]$s_5$\end{tabular}}
}
\\\hline
$\sfB_l$&\setlength{\unitlength}{11mm}\Bgenl{$\alpha_1$}{$\alpha_2$}{$\alpha_{l-2}$}{$\alpha_{l-1}$}{$\alpha_l$}\\*
&$P_{\alpha_1}(t)=[2]$\\*
&$P_{\alpha_i}(t)=\dfrac{[2i-4]!!}{[i-2]!}\quad(i=2,\dots,l)$\\*
\raisebox{40mm}{$\sfB_5$}&\raisebox{0pt}[45mm][2mm]{\setlength{\unitlength}{11mm}
\Bfuenf{\setcounter{linearcounter}{1}\begin{tabular}{c}\linear{$s_0$}{}{}{}{}{}{}\\[1mm]{[2]}\\[1mm]$s_1$\end{tabular}}%
{\setcounter{linearcounter}{0}\begin{tabular}{c}\linear{}{}{}{}{}{}{}\\[1mm]\phantom{[}1\phantom{]}\\[1mm]$s_2$\end{tabular}}%
{\setcounter{linearcounter}{1}\begin{tabular}{c}\linear{$s_0$}{}{}{}{}{}{}\\[1mm]{[2]}\\[1mm]$s_3$\end{tabular}}%
{\setcounter{linearcounter}{3}\begin{tabular}{c}\linear{$s_0$}{$s_2$}{$s_1$}{}{}{}{}\\[1mm]{[4]}\\[1mm]$s_4$\end{tabular}}%
{\begin{tabular}{c}\Bloz\\[1mm]\phantom{[}\raisebox{0pt}[0pt][0pt]{$\frac{[4][6]}{[3]}$}\phantom{]}\\[1mm]$s_5$\end{tabular}}
}
\\\hline
$\sfD_l$&\setlength{\unitlength}{11mm}\Dgenl{$\alpha_1$}{$\alpha_2$}{$\alpha_{l-3}$}{$\alpha_{l-2}$}{$\alpha_{l-1}$}{$\alpha_l$}\\*
&$P_{\alpha_1}(t)=[2]$\\*
&$P_{\alpha_i}(t)=\dfrac{[2i-4]!!}{[i-2]!}\quad(i=2,\dots,l-1)$\\*
&$P_{\alpha_l}(t)=P_{\alpha_{l-1}}(t)$\\*
\raisebox{40mm}{$\sfD_5$}&\raisebox{0pt}[45mm][25mm]{\setlength{\unitlength}{11mm}
\Dfuenf{\setcounter{linearcounter}{1}\begin{tabular}{c}\linear{$s_0$}{}{}{}{}{}{}\\[1mm]{[2]}\\[1mm]$s_1$\end{tabular}}%
{\setcounter{linearcounter}{0}\begin{tabular}{c}\linear{}{}{}{}{}{}{}\\[1mm]\phantom{[}1\phantom{]}\\[1mm]$s_2$\end{tabular}}%
{\setcounter{linearcounter}{1}\begin{tabular}{c}\linear{$s_0$}{}{}{}{}{}{}\\[1mm]{[2]}\\[1mm]$s_3$\end{tabular}}%
{\setcounter{linearcounter}{3}\begin{tabular}{c}\linear{$s_0$}{$s_2$}{$s_1$}{}{}{}{}\\[1mm]{[4]}\\[1mm]$s_4$\end{tabular}}%
{\setcounter{linearcounter}{3}\begin{tabular}{c}$s_5$\\[-14mm]\linear{$s_0$}{$s_2$}{$s_1$}{}{}{}{}\\[1mm]{[4]}\end{tabular}}
}
\\\hline
\raisebox{28mm}{$\sfE_6$}&\raisebox{0pt}[34mm][32mm]{\setlength{\unitlength}{11mm}
\Esechs{\setcounter{linearcounter}{5}\begin{tabular}{c}\linear{$s_0$}{$s_1$}{$s_2$}{$s_4$}{$s_6$}{}{}\\[1mm]{[6]}\\[1mm]$s_5$\end{tabular}}%
{\setcounter{linearcounter}{0}\begin{tabular}{c}$s_1$\\[-26mm]\linear{}{}{}{}{}{}{}\\[1mm]\phantom{[}1\phantom{]}\end{tabular}}%
{\setcounter{linearcounter}{2}\begin{tabular}{c}\linear{$s_0$}{$s_1$}{}{}{}{}{}\\[1mm]{[3]}\\[1mm]$s_3$\end{tabular}}%
{\setcounter{linearcounter}{1}\begin{tabular}{c}\linear{$s_0$}{}{}{}{}{}{}\\[1mm]{[2]}\\[1mm]$s_2$\end{tabular}}%
{\setcounter{linearcounter}{2}\begin{tabular}{c}\linear{$s_0$}{$s_1$}{}{}{}{}{}\\[1mm]{[3]}\\[1mm]$s_4$\end{tabular}}%
{\setcounter{linearcounter}{5}\begin{tabular}{c}\linear{$s_0$}{$s_1$}{$s_2$}{$s_3$}{$s_5$}{}{}\\[1mm]{[6]}\\[1mm]$s_6$\end{tabular}}
}
\\\hline
\raisebox{48mm}{$\sfE_7$}&\raisebox{0pt}[54mm][47mm]{\setlength{\unitlength}{11mm}
\Esieben{\setcounter{linearcounter}{0}\begin{tabular}{c}\linear{}{}{}{}{}{}{}\\[1mm]\phantom{[}1\phantom{]}\\[1mm]$s_1$\end{tabular}}%
{\setcounter{linearcounter}{3}\begin{tabular}{c}$s_5$\\[-12mm]\linear{$s_0$}{$s_1$}{$s_2$}{}{}{}{}\\[1mm]{[4]}\end{tabular}}%
{\setcounter{linearcounter}{1}\begin{tabular}{c}\linear{$s_0$}{}{}{}{}{}{}\\[1mm]{[2]}\\[1mm]$s_2$\end{tabular}}%
{\setcounter{linearcounter}{2}\begin{tabular}{c}\linear{$s_0$}{$s_1$}{}{}{}{}{}\\[1mm]{[3]}\\[1mm]$s_3$\end{tabular}}%
{\setcounter{linearcounter}{3}\begin{tabular}{c}\linear{$s_0$}{$s_1$}{$s_2$}{}{}{}{}\\[1mm]{[4]}\\[1mm]$s_4$\end{tabular}}%
{\setcounter{linearcounter}{5}\begin{tabular}{c}\linear{$s_0$}{$s_1$}{$s_2$}{$s_3$}{$s_5$}{}{}\\[1mm]{[6]}\\[1mm]$s_6$\end{tabular}}%
{\begin{tabular}{c}\Eloz\\[1mm]\phantom{[}\raisebox{0pt}[0pt][0pt]{$\frac{[6][10]}{[5]}$}\phantom{]}\\[1mm]$s_7$\end{tabular}}
}
\\\hline
\raisebox{38mm}{$\sfE_8$}&\raisebox{0pt}[43mm][55mm]{\setlength{\unitlength}{11mm}
\Eacht{\setcounter{linearcounter}{7}\begin{tabular}{c}\linear{$s_0$}{$s_1$}{$s_2$}{$s_3$}{$s_4$}{$s_5$}{$s_7$}\\[1mm]{[8]}\\[1mm]$s_8$\end{tabular}}%
{\setcounter{linearcounter}{5}\begin{tabular}{c}$s_7$\\[-4mm]\linear{$s_0$}{$s_1$}{$s_2$}{$s_3$}{$s_4$}{}{}\\[1mm]{[6]}\end{tabular}}%
{\setcounter{linearcounter}{5}\begin{tabular}{c}\linear{$s_0$}{$s_1$}{$s_2$}{$s_3$}{$s_4$}{}{}\\[1mm]{[6]}\\[1mm]$s_6$\end{tabular}}%
{\setcounter{linearcounter}{4}\begin{tabular}{c}\linear{$s_0$}{$s_1$}{$s_2$}{$s_3$}{}{}{}\\[1mm]{[5]}\\[1mm]$s_5$\end{tabular}}%
{\setcounter{linearcounter}{3}\begin{tabular}{c}\linear{$s_0$}{$s_1$}{$s_2$}{}{}{}{}\\[1mm]{[4]}\\[1mm]$s_4$\end{tabular}}%
{\setcounter{linearcounter}{2}\begin{tabular}{c}\linear{$s_0$}{$s_1$}{}{}{}{}{}\\[1mm]{[3]}\\[1mm]$s_3$\end{tabular}}%
{\setcounter{linearcounter}{1}\begin{tabular}{c}\linear{$s_0$}{}{}{}{}{}{}\\[1mm]{[2]}\\[1mm]$s_2$\end{tabular}}%
{\setcounter{linearcounter}{0}\begin{tabular}{c}\linear{}{}{}{}{}{}{}\\[1mm]\phantom{[}1\phantom{]}\\[1mm]$s_1$\end{tabular}}
}
\\\hline
\raisebox{23mm}{$\sfF_4$}&\raisebox{0pt}[28mm][9mm]{\setlength{\unitlength}{11mm}
\Fvier{\setcounter{linearcounter}{0}\begin{tabular}{c}\linear{}{}{}{}{}{}{}\\[1mm]\phantom{[}1\phantom{]}\\[1mm]$s_1$\end{tabular}}%
{\setcounter{linearcounter}{1}\begin{tabular}{c}\linear{$s_0$}{}{}{}{}{}{}\\[1mm]{[2]}\\[1mm]$s_2$\end{tabular}}%
{\setcounter{linearcounter}{2}\begin{tabular}{c}\linear{$s_0$}{$s_1$}{}{}{}{}{}\\[1mm]{[3]}\\[1mm]$s_3$\end{tabular}}%
{\setcounter{linearcounter}{3}\begin{tabular}{c}\linear{$s_0$}{$s_1$}{$s_2$}{}{}{}{}\\[1mm]{[4]}\\[1mm]$s_4$\end{tabular}}
}
\\\hline
\raisebox{13mm}{$\sfG_2$}&\raisebox{0pt}[18mm][9mm]{\setlength{\unitlength}{11mm}
\Gzwei{\setcounter{linearcounter}{1}\begin{tabular}{c}\linear{$s_0$}{}{}{}{}{}{}\\[1mm]{[2]}\\[1mm]$s_1$\end{tabular}}%
{\setcounter{linearcounter}{0}\begin{tabular}{c}\linear{}{}{}{}{}{}{}\\[1mm]\phantom{[}1\phantom{]}\\[1mm]$s_2$\end{tabular}}
}
\\\hline
\end{longtable}

\begin{mylemma}\label{dominantchamber}
Let $\widehat w\in\widehat W_{\perp\varphi}$ be a minimal coset representative
for a coset in $W_{\perp\varphi}\backslash\widehat W_{\perp\varphi}$. Then
$\widehat w\rho$ lies in the dominant chamber.
\end{mylemma}
\begin{proof}
We have $\ell(s_i\widehat w)>\ell(\widehat w)$ for all $i={1,\dots,l}$, namely,
for those $i$ for which $\alpha_i\perp\varphi$ by the minimality of
$\widehat w$, and for the remaining $i$ because
$s_i\notin\widehat W_{\perp\varphi}$. Hence $\widehat w$ is a minimal coset
representative for a coset in $W\backslash\widehat W$. The assertion
$\widehat w\rho\in C$ is now clear. (Since $C$ is a fundamental domain for $W$,
there is a unique $w\in W$ and a (minimal) gallery from the fundamental alcove
$A$ to $w\widehat wA$ which stays inside the dominant chamber $C$ and so
does not cross any of the hyperplanes $H_1,\dots,H_l$. By the minimality of
$\widehat w$ we get $w=1$, i.\,e., $\widehat wA\subseteq C$, or equivalently, $\widehat w\rho
\in C$.)
\end{proof}

The next lemma generalizes the orthogonality statement $w\theta\perp\theta$
in the proof of Proposition~\ref{min+}, which corresponds to $\widehat w=s_0$
and requires $\varphi\perp\theta$ in Lemma~\ref{orthogonality}.

\begin{mylemma}\label{orthogonality}
Let $\varphi\in\Phi_+^{\textup{(long)}}$. Let $w\in W$ be such that
$w\varphi=\theta$ and let $\widehat w\in\widehat W_{\perp\varphi}$. Then
$s_0w\widehat w\rho-s_0w\rho$ is perpendicular to the highest root $\theta$.
\end{mylemma}
\begin{proof}
We compute (recall that $s_\theta$ is the linear part of $s_0$)
\begin{align*}
\bigl(s_0w\widehat w\rho-s_0w\rho\bigm|\theta\bigr)&=
\bigl(s_\theta(w\widehat w\rho-w\rho)\bigm|\theta\bigr)\\
&=\bigl(w(\widehat w\rho-\rho)\bigm|{-\theta}\bigr)
=\bigl(\widehat w\rho-\rho\bigm|{-\varphi}\bigr).
\end{align*}
It remains to show that $\bigl(\widehat w\rho-\rho\bigm|{-\varphi}\bigr)=0$.
Now $\widehat W_{\perp\varphi}$ is generated by those reflections $s_i$ for
which $\alpha_i\perp\varphi$ (here $\alpha_0=-\theta$), and hence we have
$s_i(\lambda+\varphi^\perp)\subseteq\lambda+\varphi^\perp$ for all
$\lambda\in\cartan_\RR^*$. In particular, it follows that
$\widehat w\rho\in\rho+\varphi^\perp$ so that
$\bigl(\widehat w\rho-\rho\bigm|{-\varphi}\bigr)=0$.
\end{proof}

\begin{myremark}
The following computation suggests the existence of some ``shell-like structure''.
In the classical $\sfA_l$ case it means that every Young diagram decomposes into
a set of hooks that are stacked together.
\begin{align*}
\makebox[2cm][l]{$\bigl\Vert s_0w\widehat w\rho\bigr\Vert^2-\bigl\Vert s_0w\rho\bigr\Vert^2$}\\
&=\bigl\Vert s_0w\widehat w\rho-s_0w\rho\bigr\Vert^2
+2\,\bigl(\underbrace{s_0w\widehat w\rho
-s_0w\rho}_{\textstyle{}\perp\theta\makebox[0pt][l]{ (by Lemma~\ref{orthogonality})}}\bigm|s_0w\rho\bigr)\\
&=\bigl\Vert s_\theta w(\widehat w\rho-\rho)\bigr\Vert^2+2\,\big(\underbrace{
s_\theta w(\widehat w\rho-\rho)}_{\textstyle{}\perp\theta}\bigm|s_\theta w\rho
+g\theta\bigr)\\
&=\bigl\Vert\widehat w\rho-\rho\bigr\Vert^2+2\,\bigl(\widehat w\rho-\rho\bigm|
\rho\bigr)=\bigl\Vert\widehat w\rho\bigr\Vert^2-\bigl\Vert\rho\bigr\Vert^2.
\end{align*}
\end{myremark}

Here is the main theorem of this paper.

\begin{mytheorem}\label{MainTheorem}
Let $\varphi\in\Phi_+^{\textup{(long)}}$ be a positive long root and
$\widehat w\in\widehat W_{\perp\varphi}$ be the minimal coset representative
for a coset in $W_{\perp\varphi}\backslash\widehat W_{\perp\varphi}$.
To the pair $(\varphi,\widehat w)$ we associate the $\rho$-point $s_0w\widehat w\rho$ where
$w\in W$ is the Weyl group element such that $w\varphi=\theta$ and with
$\ell(w)=L(\varphi)$ as in Proposition~\ref{L}.
Then $s_0w\widehat w\rho\in2A-A$ and hence $s_0w\widehat w\rho$
is the $\rho$-point of an alcove that corresponds to a
nonzero abelian ideal $\ab^{\varphi,\widehat w}\unlhd\borel$.
Moreover, each nonzero abelian ideal occurs in this way.

Hence we have the parametrization we looked for
\begin{align*}
\bigl\{\ab\bigm|0\neq\ab\unlhd\borel,\mbox{ $\ab$ abelian}\bigr\}
&\stackrel{1:1}{\longleftrightarrow}
\coprod_{\varphi\in\Phi_+^{\textup{(long)}}}
W_{\perp\varphi}\backslash\widehat W_{\perp\varphi}
\end{align*}
where the right coset $W_{\perp\varphi}\widehat w$ in the component for
$\varphi$ on the right hand side with
$\widehat w\in\widehat W_{\perp\varphi}$ its minimal coset representative
corresponds to the abelian ideal $\ab^{\varphi,\widehat w}$.
\end{mytheorem}

\begin{myremark}
The abelian ideals $\ab^{\varphi,\min}$, $\ab^{\varphi,\min^+}$ if $\varphi\perp\theta$, and
$\ab^{\varphi,\max}$ that were mentioned earlier have the following descriptions.
\begin{align*}
\ab^{\varphi,\min}&=\ab^{\varphi,1}\\
\ab^{\varphi,\min^+}&=\ab^{\varphi,s_0}\mbox{ if $\varphi\perp\theta$}\\
\ab^{\varphi,\max}&=\ab^{\varphi,w_\circ^\varphi\widehat w_\circ^\varphi}
\end{align*}
\end{myremark}

\noindent
\begin{proof}
Since $s_0w\widehat w\neq1$---in fact, $\ell(s_0w\widehat w)=1+\ell(w)
+\ell(\widehat w)$---we must show that $s_0w\widehat w\rho\in2A$. Recall that
$$2A=\bigl\{\lambda\in\cartan_\RR^*\bigm|(\lambda|\alpha_i)\geqslant0
\mbox{ ($i={1,\dots,l}$) and }(\lambda|\theta)\leqslant1\bigr\}.$$

Let us first show that $(s_0w\widehat w\rho|\theta)\leqslant1$. For that we compute
\begin{align*}
\bigl(s_0w\widehat w\rho\bigm|\theta\bigr)&=\underbrace{\bigl(s_0w\widehat w\rho
-s_0w\rho\bigm|\theta\bigr)}_{\textstyle{}=0
\mbox{ (by Lemma~\ref{orthogonality})}}
+\bigl(s_0w\rho\bigm|\theta\bigr)\\
&=\bigl(s_\theta w\rho+g\theta\bigm|\theta\bigr)
=1-\bigl(w\rho\bigm|\theta\bigr)
=1-(\rho|\varphi)\leqslant1
\end{align*}
because $\rho$ lies in the dominant chamber and $\varphi$ is a positive root.

Next, for a simple root $\alpha_i$ with $\alpha_i\perp\theta$ we have
$\ell(s_iw)=\ell(w)+1$ by the minimality of $\ell(w)$ because
$s_iw\varphi=s_i\theta=\theta$. Hence $w^{-1}\alpha_i\in\Phi_+$ by
Lemma~\ref{length}. Now we compute
\begin{align*}
\bigl(s_0w\widehat w\rho\bigm|\alpha_i\bigr)&=\bigl(s_\theta
w\widehat w\rho+g\theta\bigm|\alpha_i\bigr)
=\bigl(w\widehat w\rho\bigm|s_\theta\alpha_i\bigr)+g\bigl(\theta\bigm|\alpha_i
\bigr)\\
&=\bigl(w\widehat w\rho\bigm|\alpha_i\bigr)
=\bigl(\widehat w\rho\bigm|w^{-1}\alpha_i\bigr)\geqslant0
\end{align*}
because $\widehat w\rho$ lies in the dominant chamber (by
Lemma~\ref{dominantchamber}) and $w^{-1}\alpha_i$
is a positive root.

Lastly, we have to deal with the simple roots that are not perpendicular to the
highest root. In most cases there is only one such root. The details will
appear in the final version of this paper.

Finally, we must show that our construction is exhaustive. It would be
nice to have a geometric argument for that. This would then prove
the First Sum Formula (Theorem~\ref{FirstSumFormula}) stated below. Momentarily
the situation is different. We prove the First Sum Formula directly and
thereby see that our construction yields all the $2^l-1$ nonzero abelian ideals.
\end{proof}

\begin{mytheorem}[First Sum Formula]\label{FirstSumFormula}
$$\sum_{\varphi\in\Phi_+^{\textup{(long)}}}P_\varphi(1)=2^l-1.$$
\end{mytheorem}
\begin{proof}
As mentioned above, the proof of the First Sum Formula is the last step
in proving Theorem~\ref{MainTheorem}.
We provide a case by case proof.

Let $S:=\sum\limits_{\varphi\in\Phi_+^{\textup{(long)}}}
P_\varphi(1)$.

\noindent
\fbox{$\sfA_l$}
Using $P_{\alpha_i+\dots+\alpha_{i+j}}(1)=\frac{(l-j-1)!}{(i-1)!\,(l-j-i)!}$
($i={1,\dots,l}$, $j={0,\dots,l-i}$) we get\\
\phantom{\fbox{$\sfA_l$}}
$S=\sum\limits_{i=1}^l\sum\limits_{j=0}^{l-i}P_{\alpha_i+\dots+\alpha_{i+j}}(1)
=\sum\limits_{i=1}^l\sum\limits_{j=0}^{l-i}\frac{(l-j-1)!}{(i-1)!\,(l-j-i)!}=2^l-1.$

For the other simple types the affine Coxeter-Dynkin graph has a tree as its underlying simple
graph. Each root $\varphi$ is of the form $\varphi=\sum\limits_{i=1}^la_i\,
\alpha_i$ and has support $\operatorname{supp}\varphi:=\bigl\{\alpha_i\bigm|a_i\neq0\bigr\}$.
Let $\operatorname{pr}(\varphi)\in\operatorname{supp}\varphi$ be the simple root which
is nearest to $\alpha_0$ when considered as nodes in the affine Coxeter-Dynkin tree.
It is clear that $P_\varphi(t)=P_{\operatorname{pr}(\varphi)}(t)$. For $i={1,\dots,l}$ let
$r_i$ be the number of positive long roots $\varphi$ for which $\operatorname{pr}(\varphi)
=\alpha_i$.
The sum $S$ can now be rewritten as $S=\sum\limits_{i=1}^lr_i\,P_{\alpha_i}(1)$.
The numbers $r_i$ can be expressed via the numbers $\nu(\sfX)={}$ the number of positive
long roots of a root system of type~$\sfX$ as in the table on page~\pageref{Poincare}. 

\renewcommand{\arraystretch}{1.2}%
\noindent
\fbox{$\sfC_l$} $\begin{array}[t]{@{}l@{}}
r_i=1\mbox{ ($i={1,\dots,l}$)}\\[1ex]
S=\sum\limits_{i=1}^l2^{i-1}=2^l-1
\end{array}$

\noindent
\fbox{$\sfB_l$} $\begin{array}[t]{@{}l@{}}
r_1=\nu(\sfA_1)=1\\
r_2=\nu(\sfB_l)-\nu(\sfB_{l-2})-\nu(\sfA_1)=4l-7\\
r_i=\nu(\sfB_{l-i+1})-\nu(\sfB_{l-i})=2l-2i\mbox{ ($i={3,\dots,l-2}$)}\\
r_{l-1}=\nu(\sfB_2)=2\mbox{ (for $l\geqslant4$)}\\
r_l=0\\[1ex]
S=1\cdot2+(4l-7)\cdot1+\sum\limits_{i=3}^{l-1}(2l-2i)\cdot2^{i-2}+0\cdot2^{l-2}=2^l-1
\end{array}$

\noindent
\fbox{$\sfD_l$} $\begin{array}[t]{@{}l@{}}
r_1=\nu(\sfA_1)=1\\
r_2=\nu(\sfD_l)-\nu(\sfD_{l-2})-\nu(\sfA_1)=4l-7\\
r_i=\nu(\sfD_{l-i+1})-\nu(\sfD_{l-i})=2l-2i\mbox{ ($i={3,\dots,l-3}$)}\\
r_{l-2}=\nu(\sfA_3)-2\,\nu(\sfA_1)=4\mbox{ (for $l\geqslant5$)}\\
r_{l-1}=r_l=\nu(\sfA_1)=1\\[1ex]
S=1\cdot2+(4l-7)\cdot1+\sum\limits_{i=3}^{l-2}(2l-2i)\cdot2^{i-2}+1\cdot2^{l-3}+1\cdot2^{l-3}\\
\phantom{S}=2^l-1
\end{array}$

For the exceptional types we write the numbers $r_i$ directly near the
corresponding node
in the Coxeter-Dynkin graph. It is clear how to compute them, e.\,g., for
$\sfE_6$, $r_1=\nu(\sfE_6)-\nu(\sfA_5)$, $r_2=\nu(\sfA_5)-2\,\nu(\sfA_2)$,
and so on.

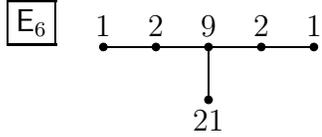
\setlength{\unitlength}{7mm}
\noindent
\fbox{$\sfE_6$} $\begin{array}[t]{@{}l@{}}
\Esechs{1}{21}{2}{9}{2}{1}\\[1ex]
S=21\cdot1+9\cdot2+2\cdot3+2\cdot3+1\cdot6+1\cdot6=2^6-1
\end{array}$

\noindent
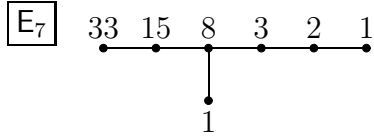
\fbox{$\sfE_7$} $\begin{array}[t]{@{}l@{}}
\Esieben{33}{1}{15}{8}{3}{2}{1}\\[1ex]
S=33\cdot1+15\cdot2+8\cdot3+3\cdot4+1\cdot4+2\cdot6+1\cdot12=2^7-1
\end{array}$

\noindent
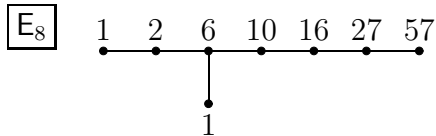
\fbox{$\sfE_8$} $\begin{array}[t]{@{}l@{}}
\Eacht{1}{1}{2}{6}{10}{16}{27}{57}\\[1ex]
S=57\cdot1+27\cdot2+16\cdot3+10\cdot4+6\cdot5+2\cdot6+1\cdot6+1\cdot8\\
\phantom{S}=2^8-1
\end{array}$

\noindent
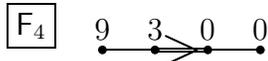
\fbox{$\sfF_4$} $\begin{array}[t]{@{}l@{}}
\Fvier{9}{3}{0}{0}\\[1ex]
S=9\cdot1+3\cdot2+0\cdot3+0\cdot4=2^4-1
\end{array}$

\noindent
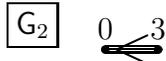
\fbox{$\sfG_2$} $\begin{array}[t]{@{}l@{}}
\Gzwei{0}{3}\\[1ex]
S=3\cdot1+0\cdot2=2^2-1
\end{array}$
\renewcommand{\arraystretch}{1}%

\end{proof}

Among the maximal abelian ideals---we shall look at them more closely in
the next section---are those whose dimension is maximal.
A.~Malcev \cite{Mal} calculated the dimension for each type. Our approach allows
to express these dimensions in a uniform way as
$g-1+N'-N''$ where $g$ is the dual Coxeter number of $\lie$ and $N',N''$
are the numbers of positive roots of certain root subsystems. 

\begin{mycorollary}
The maximal dimension of an abelian ideal in $\borel$ can be expressed as
$g-1+\widehat N_{\perp\alpha}-N_{\perp\alpha}$ where $\alpha
\in\Pi^{\textup{(long)}}$ is a long simple root such that
the difference $\widehat N_{\perp\alpha}-N_{\perp\alpha}$
is maximal. Here $\widehat N_{\perp\alpha}=\ell(\widehat w_\circ^\alpha)$
is the length of the longest element $\widehat w_\circ^\alpha
\in\widehat W_{\perp\alpha}$, or, equivalently, the number of
positive roots of a root system of type $\widehat\sfX_{\perp\alpha}$.
Analogously, we have $N_{\perp\alpha}=\ell(w_\circ^\alpha)$ for
$w_\circ^\alpha\in W_{\perp\alpha}$ its longest element.
\end{mycorollary}

\subsection*{Abelian ideals of maximal dimension}\label{maxdim}
The fourth column in the table shows the Coxeter-Dynkin graph
$\textrm{CD}_\sfX$ of type $\sfX$ and the affine Coxeter-Dynkin
graph $\textrm{CD}_{\widehat\sfX}$ with the node $\alpha_0$ encircled.
{\footnotesize
\renewcommand{\arraystretch}{2.4}%
$$\begin{array}{|l|c|c|c|c|c|c|}\hline
\multicolumn{1}{|c|}{\sfX}&g_\sfX-1&N_\sfX&
\textrm{CD}_\sfX\quad\textrm{CD}_{\widehat\sfX}\quad\alpha
&\widehat\sfX_{\perp\alpha}&\sfX_{\perp\alpha}
&\max\dim\\\hline
\sfA_1&1&1&\AoneAff{$\alpha$}{}&\varnothing&\varnothing&1\\
\sfA_2&2&3&\AtwoAff{$\alpha$}{}{}&\varnothing&\varnothing&2\\
\begin{array}[t]{@{}l@{}}\sfA_l\\[-3ex]
\substack{n\geqslant3\hfill\\l\textrm{ odd}}\end{array}
&l&\dfrac{l(l+1)}{2}&\AlAffod{}{$\alpha$}{}{}&
\sfA_{l-2}&\sfA_{\frac{l-3}{2}}+\sfA_{\frac{l-3}{2}}
&\dfrac{(l+1)^2}{4}\\
\begin{array}[t]{@{}l@{}}\sfA_l\\[-3ex]
\substack{l\geqslant4\hfill\\l\textrm{ even}}\end{array}
&l&\dfrac{l(l+1)}{2}&\AlAffev{}{$\alpha$}{}{}{}&
\sfA_{l-2}&\sfA_{\frac{l-2}{2}}+\sfA_{\frac{l-4}{2}}
&\dfrac{l^2+2l}{4}\\
\begin{array}[t]{@{}l@{}}\sfC_l\\[-4ex]
\scriptstyle l\geqslant2\end{array}&l&l^2&
\ClAff{}{}{}{}{$\alpha$}{}
&\sfC_{l-1}&\sfA_{l-2}
&\dfrac{l^2+l}{2}\\
\sfB_3&4&9&\BthreeAff{}{}{}{$\alpha$}&\sfA_1+\sfA_1&\sfA_1&5\\
\begin{array}[t]{@{}l@{}}\sfB_l\\[-4ex]
\scriptstyle l\geqslant4\end{array}&2l-2&l^2&
\BlAff{}{}{}{}{$\alpha$}{}{}
&\sfD_{l-2}&
\sfA_{l-3}&\dfrac{l^2-l+2}{2}\\
\begin{array}[t]{@{}l@{}}\sfD_l\\[-4ex]
\scriptstyle l\geqslant4\end{array}&2l-3&l(l-1)&
\DlAff{}{}{}{}{}{$\alpha$}{}{}
&\sfD_{l-2}+\sfA_1
&\sfA_{l-3}+\sfA_1&\dfrac{l^2-l}{2}\\
\sfE_6&11&36&\EaAff{$\alpha$}{}{}{}{}{}{}&\sfA_5&\sfA_4&16\\
\sfE_7&17&63&\EbAff{}{}{}{}{}{}{$\alpha$}{}&\sfD_6&\sfD_5&27\\
\sfE_8&29&120&\EcAff{$\alpha$}{}{}{}{}{}{}{}{}&\sfA_7&\sfA_6&36\\
\sfF_4&8&24&\FAff{}{$\alpha$}{}{}{}&\sfA_1+\sfA_1&\sfA_1&9\\
\sfG_2&3&6&\GAff{}{$\alpha$}{}&\varnothing&\varnothing&3\\\hline
\multicolumn{7}{p{122mm}}
{\footnotesize
The usual conventions apply, namely, $\sfA_0=\varnothing$,
$\sfC_1=\sfA_1$, $\sfB_2=\sfC_2$, $\sfD_2=\sfA_1+\sfA_1$.}
\end{array}$$
\renewcommand{\arraystretch}{1}%
}

\begin{myremark}
In the previous table, in some cases there are several possibilities
for the simple long root $\alpha$ that yields an abelian ideal of
maximal dimension. By inspection we see that the number of abelian ideals
in $\borel$ of maximal dimension is
\begin{enumerate}
\item[] $3$ for type $\sfD_4$,
\item[] $2$ for types $\sfA_l$ ($l$ even), $\sfD_l$ ($l>4$), and $\sfE_6$, and
\item[] $1$ for the other types.
\end{enumerate}
\end{myremark}

\begin{myremark}
Instead of taking $\widehat\sfX_{\perp\alpha}$ and $\sfX_{\perp\alpha}$ one
could already delete the common components (nonvoid for the
types $\sfB_3$, $\sfD_l$, and $\sfF_4$).
\end{myremark}

\begin{myremark}
As already mentioned the numbers in the rightmost column of the table above
were first computed case by case by A.~Malcev~\cite{Mal}. In the recent paper
\cite{Bo} B.~Boe computed, again case by case,
the maximal length $\ell(\widehat w)$ of an affine Weyl group element $\widehat w$
such that $\widehat wA\subseteq(k+1)A$;
see \cite[Table~1]{Bo} but with the types $\sfC_l$ and $\sfB_l$ interchanged
because there the highest \emph{short} root is used to define the tessellation by alcoves.
Neither Boe's paper nor its review paper \cite{Sr} mentions
the connexion with Malcev's result about abelian ideals.
\end{myremark}

\section{Maximal abelian ideals}\label{maximal}
We are now interested in the maximal abelian ideals in
$\borel$. It has been observed in \cite{PR} that the number of maximal abelian
ideals in a fixed Borel subalgebra of $\lie$ equals the number of
long simple roots. A canonical one-to-one correspondence was exhibited between
the two sets. However, the proof was based on a case by case consideration
and was therefore rather unsatisfactory.
Here we will give a geometric approach which makes the whole picture very
transparent.

We know that each abelian ideal $\ab\unlhd\borel$ corresponds to an alcove
$\widehat wA\subseteq2A$. If $\widehat wA$ has no facet\footnote{We follow
the traditional terminology which speaks of ``facets'' for ``faces of
codimension one''. This disagrees with the English translation of the
Bourbaki volume~\cite{B}. As far as I understand and remembering A.~Borel's
course (taught in German) about Tits buildings many years ago, the terminology
in French is ``face'' for ``facette de codimension une''.}
lying in the wall $2H_0$, then $\ab$ cannot be a maximal ideal.
Hence each maximal abelian ideal has an alcove with one facet lying in the
wall $2H_0$. It is convenient to have some terminology which describes this
geometric situation.

\begin{mydefinition}
An \textbf{upper alcove} $\widehat wA$ is an alcove in $2A$ such that
one facet of $\widehat wA$ lies in the wall $2H_0$. For an upper alcove
$\widehat wA$ the \textbf{lower vertex} is the vertex that sticks out,
i.\,e., does not lie in the wall $2H_0$.
\end{mydefinition}

Let us look at some examples. For type $\sfA_2$ there are two upper alcoves,
namely, those with $\rho$-points $s_0s_2\rho$ and $s_0s_1\rho$. Both belong
to maximal abelian ideals, namely, $\ab^{\alpha_1,1}$ and $\ab^{\alpha_2,1}$.
For the former alcove, the lower vertex has type~$1$ and for the latter type~$2$.
For type $\sfC_2$ there are again two upper alcoves, with $\rho$-points
$s_0s_1\rho$ (ideal $\alpha^{\alpha_2,1}$) and $s_0s_1s_0\rho$
(ideal $\ab^{\alpha_2,s_0}$), both with the same lower vertex of type~$2$.
Only the latter belongs to the maximal abelian ideal.
For type $\sfG_2$ (see the picture on page~\pageref{pictureG2}) there is only
one upper alcove, with $\rho$-point $s_0s_2s_1\rho$ and lower vertex of
type~$2$.
\\[2mm]

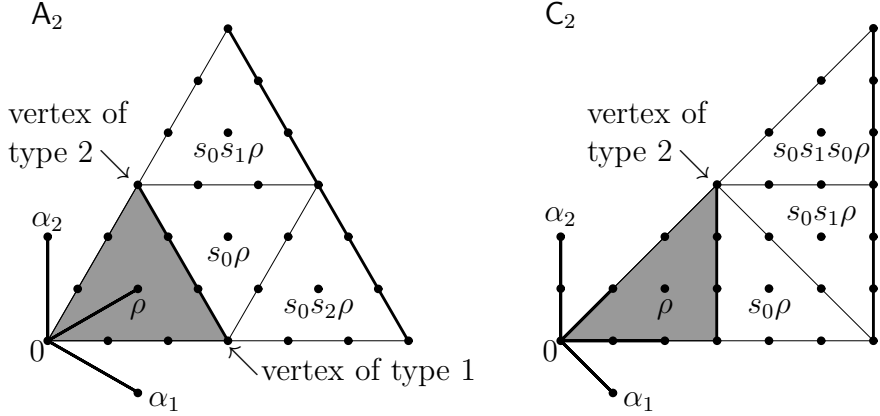
\setlength{\unitlength}{.0000666667\textwidth}
\setlength{\unitlength}{0.008mm}
\begin{picture}(7500,6062)(0,-866)
\texture{55888888 88555555 5522a222 a2555555 55888888 88555555 552a2a2a 2a555555 
        55888888 88555555 55a222a2 22555555 55888888 88555555 552a2a2a 2a555555 
        55888888 88555555 5522a222 a2555555 55888888 88555555 552a2a2a 2a555555 
        55888888 88555555 55a222a2 22555555 55888888 88555555 552a2a2a 2a555555 }
\shade\path(0,0)(3000,0)(1500,2598)(0,0)
\Thicklines
\path(0,0)(0,1732)(0,0)(1500,866)(0,0)(1500,-866)(0,0)
\thinlines
\path(0,0)(6000,0)(3000,5196)(0,0)
\path(3000,0)(4500,2598)(1500,2598)(3000,0)
\Thicklines
\path(3000,0)(1500,2598)
\path(6000,0)(3000,5196)
\put(0,0){\circle*{100}}
\multiput(0,1732)(500,-866){4}{\circle*{100}}
\multiput(1000,1732)(500,-866){3}{\circle*{100}}
\multiput(1500,2598)(500,-866){4}{\circle*{100}}
\multiput(2000,3464)(500,-866){5}{\circle*{100}}
\multiput(2500,4330)(500,-866){6}{\circle*{100}}
\multiput(3000,5196)(500,-866){7}{\circle*{100}}
\put(0,1832){\makebox(0,0)[b]{$\alpha_2$}}
\put(1673,-866){\makebox(0,0)[lt]{$\alpha_1$}}
\put(3000,-100){\makebox(0,0)[tl]{${}^{\textstyle\nwarrow}$vertex of type $1$}}
\put(1400,2698){\makebox(0,0)[br]{
\begin{tabular}{l@{}}vertex of\\type $2\,_{\textstyle\searrow}$\end{tabular}}}
\put(1500,693){\makebox(0,0)[t]{$\rho$}}
\put(3000,1559){\makebox(0,0)[t]{$s_0\rho$}}
\put(3000,3291){\makebox(0,0)[t]{$s_0s_1\rho$}}
\put(4500,693){\makebox(0,0)[t]{$s_0s_2\rho$}}
\put(-173,0){\makebox(0,0)[t]{$0$}}
\put(0,5196){\makebox(0,0)[b]{$\sfA_2$}}
\end{picture}
\quad
\setlength{\unitlength}{.0577334\textwidth}
\setlength{\unitlength}{6.92801mm}
\begin{picture}(6,7)(0,-1)
\texture{55888888 88555555 5522a222 a2555555 55888888 88555555 552a2a2a 2a555555 
        55888888 88555555 55a222a2 22555555 55888888 88555555 552a2a2a 2a555555 
        55888888 88555555 5522a222 a2555555 55888888 88555555 552a2a2a 2a555555 
        55888888 88555555 55a222a2 22555555 55888888 88555555 552a2a2a 2a555555 }
\shade\path(0,0)(3,0)(3,3)(0,0)
\Thicklines
\path(0,0)(0,2)(0,0)(1,1)(0,0)(2,0)(0,0)(1,-1)(0,0)
\thinlines
\path(0,0)(6,0)(6,6)(0,0)
\path(3,0)(3,3)(6,3)
\path(3,3)(6,0)
\Thicklines
\path(3,0)(3,3)
\path(6,0)(6,6)
\multiput(0,0)(0,1){3}{\circle*{.115473}}
\multiput(1,-1)(0,1){3}{\circle*{.115473}}
\multiput(2,0)(0,1){3}{\circle*{.115473}}
\multiput(3,0)(0,1){4}{\circle*{.115473}}
\multiput(4,0)(0,1){5}{\circle*{.115473}}
\multiput(5,0)(0,1){6}{\circle*{.115473}}
\multiput(6,0)(0,1){7}{\circle*{.115473}}
\put(0,2.2){\makebox(0,0)[b]{$\alpha_2$}}
\put(1.2,-1){\makebox(0,0)[lt]{$\alpha_1$}}
\put(2,0.8){\makebox(0,0)[t]{$\rho$}}
\put(4,0.8){\makebox(0,0)[t]{$s_0\rho$}}
\put(5,2.2){\makebox(0,0)[b]{$s_0s_1\rho$}}
\put(5,3.8){\makebox(0,0)[t]{$s_0s_1s_0\rho$}}
\put(2.88453,3.11547){\makebox(0,0)[br]{
\begin{tabular}{l@{}}vertex of\\type $2\,_{\textstyle\searrow}$\end{tabular}}}
\put(-.2,0){\makebox(0,0)[t]{$0$}}
\put(0,6){\makebox(0,0)[b]{$\sfC_2$}}
\end{picture}
\setlength{\unitlength}{3.5mm}
\\

From previous results we already know that the lower vertices are in
one-to-one correspondence with the long simple roots.

Another way for proving that each lower vertex has the type of a long simple
root can be deduced from the following proposition which we also use for
our Second Sum Formula (Theorem~\ref{SecondSumFormula}).

\begin{myproposition}
$\vol_{l-1}(F_0):\dots:\vol_{l-1}(F_l)=
\Vert\alpha_0\Vert\,n_0:\dots:\Vert\alpha_l\Vert\,n_l$.
\end{myproposition}
\begin{proof}
Recall that the vertices of the fundamental alcove $A$ with facets
$F_0,\dots,F_l$ are $0,\frac{\stackrel{\vee}{\varpi}_1}{n_1},\dots,
\frac{\stackrel{\vee}{\varpi}_l}{n_l}$. We compute the volume of an
alcove in two different ways.\par
The volume of the pyramid $A$ over $F_0$ with apex $0$ is $\frac1l$ times
$\vol_{l-1}(F_0)$ times the distance of the apex $0$ from the hyperplane
$H_0$ supporting the face $F_0$. This distance is $\frac{1}{2\Vert\theta\Vert}$
because $\frac{\theta}{2\Vert\theta\Vert^2}=\frac12g\theta\in H_0$ is the orthogonal projection
of the apex $0$ to $H_0$. On the other hand, the volume of $A$ is $\frac{1}{l!}$
times the volume $D=\Bigl|\frac{\stackrel{\vee}{\varpi}_1}{n_1}\wedge\dots\wedge
\frac{\stackrel{\vee}{\varpi}_l}{n_l}\Bigr|$ of the parallelepiped spanned by the
vectors $\frac{\stackrel{\vee}{\varpi}_i}{n_i}$ ($i={1,\dots,l}$). Hence
$\vol_{l-1}(F_0)=2\,\Vert\theta\Vert\,\frac{D}{(l-1)!}
=2\,\Vert\alpha_0\Vert\,n_0\,\frac{D}{(l-1)!}$\,.\par
Now we compute the $(l-1)$-dimensional volume of an $(l-1)$-simplex $F_i$
($i={1,\dots,l}$) as the $l$-dimensional volume of the prism $F_i\times I$
where $I$ is a unit interval perpendicular to $F_i$. Hence
\begin{align*}
\textstyle\vol_{l-1}(F_i)&=\textstyle\frac{1}{(l-1)!}\Bigl|
\frac{\stackrel{\vee}{\varpi}_1}{n_1}\wedge\dots\wedge
\frac{\stackrel{\vee}{\varpi}_{i-1}}{n_{i-1}}\wedge
\frac{\alpha_i}{\Vert\alpha_i\Vert}\wedge
\frac{\stackrel{\vee}{\varpi}_{i+1}}{n_{i+1}}\wedge\dots\wedge
\frac{\stackrel{\vee}{\varpi}_l}{n_l}\Bigr|\\
&=\textstyle2\,\Vert\alpha_i\Vert\,n_i\,\frac{D}{(l-1)!}
\end{align*}
because
$\alpha_i=2\sum\limits_{k=1}^l\bigl(\alpha_i\bigm|\stackrel{\vee}{\varpi}_k
\bigr)\alpha_k=2\sum\limits_{k=1}^l\bigl(\alpha_i\bigm|\alpha_k\bigr)
\!\stackrel{\vee}{\varpi}_k$.\par
This proves the proposition.
\end{proof}

\label{connexiongh}
\begin{myremark}
The two formulae
\begin{align*}
\vol_l(A)&=\frac1l\cdot\vol_{l-1}(F_0)\cdot\frac{1}{2\,\Vert\theta\Vert}
\intertext{and}
\vol_l(A)&=\sum_{i=0}^l\vol_l(\mbox{pyramid with base $F_i$ and apex $\rho$})\\
&=\frac1l\cdot\sum_{i=0}^l\operatorname{dist}(\rho,F_i)\cdot
\vol_{l-1}(F_i)\\
&=\frac1l\cdot\sum_{i=0}^l\frac12\,\Vert\alpha_i\Vert\cdot
\frac{n_i\,\Vert\alpha_i\Vert}{n_0\,\Vert\alpha_0\Vert}\cdot\vol_{l-1}(F_0)\\
&=\frac1l\cdot\frac{1}{2\,\Vert\theta\Vert}\cdot\vol_{l-1}(F_0)\cdot
\sum_{i=0}^ln_i\,\Vert\alpha_i\Vert^2
\end{align*}
show that $\sum\limits_{i=0}^ln_i\,\Vert\alpha_i\Vert^2=1$.
\end{myremark}

The previous proposition makes clear that the lower vertex of an upper alcove
cannot have the type of a short simple root for commensurability reasons.
(Here the convention is that a root is long and not short if the root system
is simply laced.)
We next observe that no lower vertex can have type $0$. For volume
reasons such a vertex would have to lie in $F_0$ which is absurd.

\begin{mytheorem}[Second Sum Formula]\label{SecondSumFormula}
The following sum formula holds.
$$\sum_{i=1}^ln_i\,P_{\alpha_i}(1)=2^{l-1}$$
\end{mytheorem}
\begin{proof}
We look at $\vol_{l-1}(2F_0)$ and compute the volume in two ways. First,
of course, $\vol_{l-1}(2F_0)=2^{l-1}\vol_{l-1}(F_0)$. Second, consider the
tessellation of $2F_0$ induced by the tessellation of $\cartan_\RR^*$ by
the alcoves. Namely, $\vol_{l-1}(F_i)=n_i\cdot\vol_{l-1}(F_0)$ and for
each $\alpha_i\in\Pi^{\textup{(long)}}$ there are $P_{\alpha_i}(1)$
simplices of type $i$ in the tessellation of $2F_0$.
\end{proof}

\section{An example}
In this section we try to exemplify the remark about the ``shell-like'' structure
before Theorem~\ref{MainTheorem}. Take $\sfA_{11}$. Look at the abelian ideal
corresponding to the following Young diagram.

$$\Yboxdim5mm
\young(\refla\reflb\reflc\refld\refle,\refll\refla\reflb\reflc,%
\reflk\refll\refla\reflb,\reflj\reflk\refll\refla,%
\refli\reflj\reflk\refll,\reflh\refli\reflj,\reflg\reflh)$$

$$w=s_0\underbrace{s_1\,s_2\,s_3\,s_4\,s_{11}\,s_{10}\,s_9\,s_8\,s_7\,s_6}_{\textstyle{}=:w_1}
s_0\underbrace{s_1\,s_2\,s_{11}\,s_{10}\,s_9\,s_8\,s_7}_{\textstyle{}=:w_2}
s_0\underbrace{s_2\,s_{11}\,s_{10}\,s_9}_{\textstyle{}=:w_3}
s_0\underbrace{s_{11}}_{\textstyle\makebox[0pt]{${}=:w_4$}}$$

\renewcommand{\arraystretch}{1.2}%
$$\begin{array}{|c|c|c|c|c|}\hline
r&s_{i_r}&s_{i_1}\dots s_{i_r}\rho-s_{i_1}\dots s_{i_{r-1}}\rho&
s_{j_r}&s_{j_1}\dots s_{j_r}\rho-s_{j_1}\dots s_{j_{r-1}}\rho\\\hline
\phantom{0}1&s_0&1\,1\,1\,1\,1\,1\,1\,1\,1\,1\,1&s_0&1\,1\,1\,1\,1\,1\,1\,1\,1\,1\,1\\
\phantom{0}2&s_1&0\,1\,1\,1\,1\,1\,1\,1\,1\,1\,1&s_1&0\,1\,1\,1\,1\,1\,1\,1\,1\,1\,1\\
\phantom{0}3&s_2&0\,0\,1\,1\,1\,1\,1\,1\,1\,1\,1&s_{11}&1\,1\,1\,1\,1\,1\,1\,1\,1\,1\,0\\\cline{4-5}
\phantom{0}4&s_3&0\,0\,0\,1\,1\,1\,1\,1\,1\,1\,1&s_0&0\,1\,1\,1\,1\,1\,1\,1\,1\,1\,0\\
\phantom{0}5&s_4&0\,0\,0\,0\,1\,1\,1\,1\,1\,1\,1&s_2&0\,0\,1\,1\,1\,1\,1\,1\,1\,1\,1\\
\phantom{0}6&s_{11}&1\,1\,1\,1\,1\,1\,1\,1\,1\,1\,0&s_1&0\,0\,1\,1\,1\,1\,1\,1\,1\,1\,0\\
\phantom{0}7&s_{10}&1\,1\,1\,1\,1\,1\,1\,1\,1\,0\,0&s_{10}&1\,1\,1\,1\,1\,1\,1\,1\,1\,0\,0\\
\phantom{0}8&s_9&1\,1\,1\,1\,1\,1\,1\,1\,0\,0\,0&s_{11}&0\,1\,1\,1\,1\,1\,1\,1\,1\,0\,0\\\cline{4-5}
\phantom{0}9&s_8&1\,1\,1\,1\,1\,1\,1\,0\,0\,0\,0&s_0&0\,0\,1\,1\,1\,1\,1\,1\,1\,0\,0\\
10&s_7&1\,1\,1\,1\,1\,1\,0\,0\,0\,0\,0&s_3&0\,0\,0\,1\,1\,1\,1\,1\,1\,1\,1\\
11&s_6&1\,1\,1\,1\,1\,0\,0\,0\,0\,0\,0&s_2&0\,0\,0\,1\,1\,1\,1\,1\,1\,1\,0\\\cline{2-3}
12&s_0&0\,1\,1\,1\,1\,1\,1\,1\,1\,1\,0&s_1&0\,0\,0\,1\,1\,1\,1\,1\,1\,0\,0\\
13&s_1&0\,0\,1\,1\,1\,1\,1\,1\,1\,1\,0&s_9&1\,1\,1\,1\,1\,1\,1\,1\,0\,0\,0\\
14&s_2&0\,0\,0\,1\,1\,1\,1\,1\,1\,1\,0&s_{10}&0\,1\,1\,1\,1\,1\,1\,1\,0\,0\,0\\
15&s_{11}&0\,1\,1\,1\,1\,1\,1\,1\,1\,0\,0&s_{11}&0\,0\,1\,1\,1\,1\,1\,1\,0\,0\,0\\\cline{4-5}
16&s_{10}&0\,1\,1\,1\,1\,1\,1\,1\,0\,0\,0&s_0&0\,0\,0\,1\,1\,1\,1\,1\,0\,0\,0\\
17&s_9&0\,1\,1\,1\,1\,1\,1\,0\,0\,0\,0&s_4&0\,0\,0\,0\,1\,1\,1\,1\,1\,1\,1\\
18&s_8&0\,1\,1\,1\,1\,1\,0\,0\,0\,0\,0&s_8&1\,1\,1\,1\,1\,1\,1\,0\,0\,0\,0\\
19&s_7&0\,1\,1\,1\,1\,0\,0\,0\,0\,0\,0&s_9&0\,1\,1\,1\,1\,1\,1\,0\,0\,0\,0\\\cline{2-3}
20&s_0&0\,0\,1\,1\,1\,1\,1\,1\,1\,0\,0&s_{10}&0\,0\,1\,1\,1\,1\,1\,0\,0\,0\,0\\
21&s_1&0\,0\,0\,1\,1\,1\,1\,1\,1\,0\,0&s_{11}&0\,0\,0\,1\,1\,1\,1\,0\,0\,0\,0\\
22&s_{11}&0\,0\,1\,1\,1\,1\,1\,1\,0\,0\,0&s_7&1\,1\,1\,1\,1\,1\,0\,0\,0\,0\,0\\
23&s_{10}&0\,0\,1\,1\,1\,1\,1\,0\,0\,0\,0&s_8&0\,1\,1\,1\,1\,1\,0\,0\,0\,0\,0\\
24&s_9&0\,0\,1\,1\,1\,1\,0\,0\,0\,0\,0&s_9&0\,0\,1\,1\,1\,1\,0\,0\,0\,0\,0\\\cline{2-3}
25&s_0&0\,0\,0\,1\,1\,1\,1\,1\,0\,0\,0&s_6&1\,1\,1\,1\,1\,0\,0\,0\,0\,0\,0\\
26&s_{11}&0\,0\,0\,1\,1\,1\,1\,0\,0\,0\,0&s_7&0\,1\,1\,1\,1\,0\,0\,0\,0\,0\,0\\\hline
\end{array}$$
\renewcommand{\arraystretch}{1}%

\section{Symmetries of the Hasse graphs}
In this section we look at the Hasse graph of the poset of abelian ideals
in $\borel$ and determine its group of symmetries.
A natural geometric realization of the Hasse graph of abelian ideals in
$\borel$ lives in $\cartan_\RR^*$. The nodes are the $\rho$-points of the
alcoves contained in $2A$. Two $\rho$-points are connected if and only if
their alcoves are adjacent. Surely, the geometric symmetry group of this $1$-dimensional
complex is a subgroup of the abstract symmetry group of the Hasse graph.
In fact, it turns out that the two symmetry groups coincide unless $\lie$
has type $\sfC_3$ or $\sfG_2$. In the former case the abstract Hasse graph
has the following shape with symmetry group $\ZZ/2\ZZ\times\ZZ/2\ZZ$.
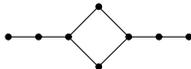
$${\setlength{\unitlength}{4mm}
\begin{picture}(6,2)
\path(0,1)(2,1)(3,2)(4,1)\path(2,1)(3,0)(4,1)(6,1)
\multiput(0,1)(1,0){3}{\circle*{.2}}\multiput(4,1)(1,0){3}{\circle*{.2}}
\multiput(3,0)(0,2){2}{\circle*{.2}}
\end{picture}}$$
In the natural geometric realization the cycle of length four is actually not a square
but a rectangle with side ratio $\sqrt2:1$. Thus the geometric symmetry group
collapses to $\ZZ/2\ZZ$. In the case of $\sfG_2$ the two groups are
$1$ and $\ZZ/2\ZZ$ (see page~\pageref{pictureG2}).

Loosely speaking, the geometric symmetry group is the symmetry group of $2A$,
hence isomorphic to the symmetry group of the affine Coxeter-Dynkin graph.
Going through the classification one sees that the abstract symmetry
group is the same as the geometric one, with the two exceptions mentioned above.

\section{Examples: rank $4$}
In the next few pages we show the Hasse graphs of the posets of abelian ideals
in $\borel$ for the five simple types of rank~$4$. Each node of the Hasse graph
consists of a diagram of a shape of which an enlarged version is drawn before
the Hasse graph. The boxes of the enlarged version are filled with the
nonforbidden\footnote{A forbidden positive root $\varphi$ is such that
$\theta-2\varphi$ is a sum of positive roots. Then the root space $\lie_\varphi$
cannot belong to an abelian ideal of $\borel$.}
positive roots. Each node in the Hasse graph corresponds to the abelian ideal
$\bigoplus\limits_\varphi\lie_\varphi$ where $\varphi$ runs over the positive roots marked
by a dot.

The arrows in the Hasse graphs have the following meaning. Each node which is not
the source of an arrow corresponds to an ideal of the form $\ab^{\varphi,\min}$
for some $\varphi\in\Phi_+^{\textup{(long)}}$. For $\varphi$ a long simple root,
we have labeled the node belonging to $\ab^{\varphi,\min}$. The passage from
$0\neq\ab$ to $\ab^\anc$ corresponds to following the arrows till one arrives at a sink.
Finally, an arrow points from the empty diagram ($\ab=0$) to the diagram filled
with one dot ($\ab=\lie_\theta$).
Disregard the arrows for the automorphism groups.

\newpage
\setlength{\unitlength}{11mm}
\subsection*{$\mathsf A_l$\quad
$\Agenl{$\alpha_1$}{$\alpha_2$}{$\alpha_{l-1}$}{$\alpha_l$}$}

\setlength{\unitlength}{3.5mm}
$$\Yboxdim14mm
\young(\Axxxx\Axxxo\Axxoo\Axooo,\Aoxxx\Aoxxo\Aoxoo,\Aooxx\Aooxo,\Aooox)$$

\begin{center}
\mbox{}{\tiny\Yboxdim1.5mm
\xymatrix@R=4mm@C=5mm{&&&{\young(\bullet\bullet\bullet\ ,\bullet\bullet\bullet,\ \ ,\ )}&&
{\young(\bullet\bullet\ \ ,\bullet\bullet\ ,\bullet\bullet,\ )}\\
&&&{\young(\bullet\bullet\bullet\ ,\bullet\bullet\ ,\ \ ,\ )}\ar@{<-}[u]^4&&
{\young(\bullet\bullet\ \ ,\bullet\bullet\ ,\bullet\ ,\ )}\ar@{<-}[u]^1\\
{\makebox[0pt][r]{\raisebox{3ex}{\normalsize$\alpha_1$\ }}\young(\bullet\bullet\bullet\bullet,\ \ \ ,\ \ ,\ )}&&
{\makebox[0pt][r]{\raisebox{3ex}{\normalsize$\alpha_2$\ }}\young(\bullet\bullet\bullet\ ,\bullet\ \ ,\ \ ,\ )}\ar@{<-}[ur]^0&&
{\young(\bullet\bullet\ \ ,\bullet\bullet\ ,\ \ ,\ )}\ar@{-}[ul]^3\ar@{-}[ur]^2&&
{\makebox[0pt][r]{\raisebox{3ex}{\normalsize$\alpha_3$\ }}\young(\bullet\bullet\ \ ,\bullet\ \ ,\bullet\ ,\ )}\ar@{<-}[ul]^0&&
{\makebox[0pt][r]{\raisebox{3ex}{\normalsize$\alpha_4$\ }}\young(\bullet\ \ \ ,\bullet\ \ ,\bullet\ ,\bullet)}\\
&&{\young(\bullet\bullet\bullet\ ,\ \ \ ,\ \ ,\ )}\ar@{-}[ull]^2\ar@{-}[u]^1&&
{\young(\bullet\bullet\ \ ,\bullet\ \ ,\ \ ,\ )}\ar@{-}[ull]^3
\ar@{<-}[u]^0\ar@{-}[urr]^2&&
{\young(\bullet\ \ \ ,\bullet\ \ ,\bullet\ ,\ )}\ar@{-}[u]^4\ar@{-}[urr]^3\\
&&&{\young(\bullet\bullet\ \ ,\ \ \ ,\ \ ,\ )}\ar@{-}[ul]^3\ar@{-}[ur]^1
&&{\young(\bullet\ \ \ ,\bullet\ \ ,\ \ ,\ )}\ar@{-}[ul]^4\ar@{-}[ur]^2\\
&&&&{\young(\bullet\ \ \ ,\ \ \ ,\ \ ,\ )}\ar@{-}[ul]^4\ar@{-}[ur]^1\\
&&&&{\young(\ \ \ \ ,\ \ \ ,\ \ ,\ )}\ar@{->}[u]^0}}
\end{center}
\mbox{}\\
$\Aut\bigl(\Hasse(\sfA_1)\bigr)\cong\ZZ/2\ZZ$\\
$\Aut\bigl(\Hasse(\sfA_l)\bigr)\cong\Dih_{l+1}\quad(l\geqslant2)$

\newpage
\setlength{\unitlength}{11mm}
\subsection*{$\mathsf C_l$\quad
$\Cgenl{$\alpha_1$}{$\alpha_2$}{$\alpha_{l-2}$}{$\alpha_{l-1}$}{$\alpha_l$}$}

\mbox{}
\setlength{\unitlength}{3.5mm}
$$\Yboxdim14mm
\young(\Cwwwx\Cxwwx\Cxxwx\Cxxxx,:\Cowwx\Coxwx\Coxxx,::\Coowx\Cooxx,:::\Cooox)$$

\begin{center}
\mbox{}{\tiny\Yboxdim1.5mm
\xymatrix@R=4mm{&{\young(\bullet\bullet\bullet\bullet,:\bullet\bullet\bullet,%
::\bullet\bullet,:::\bullet)}\\
&{\young(\bullet\bullet\bullet\bullet,:\bullet\bullet\bullet,%
::\bullet\bullet,:::\ )}\ar@{<-}[u]^0\\
&{\young(\bullet\bullet\bullet\bullet,:\bullet\bullet\bullet,%
::\bullet\ ,:::\ )}\ar@{<-}[u]^1\\
{\young(\bullet\bullet\bullet\bullet,:\bullet\bullet\bullet,%
::\ \ ,:::\ )}\ar@{<-}[ur]^0&&
{\young(\bullet\bullet\bullet\bullet,:\bullet\bullet\ ,%
::\bullet\ ,:::\ )}\ar@{<-}[ul]^2\\
{\young(\bullet\bullet\bullet\bullet,:\bullet\bullet\ ,::\ \ ,:::\ )}\ar@{<-}[u]^2\ar@{<-}[urr]^0&&
{\young(\bullet\bullet\bullet\ ,:\bullet\bullet\ ,::\bullet\ ,:::\ )}\ar@{-}[u]^3\\
{\young(\bullet\bullet\bullet\bullet,:\bullet\ \ ,::\ \ ,:::\ )}\ar@{<-}[u]^1&&
{\young(\bullet\bullet\bullet\ ,:\bullet\bullet\ ,::\ \ ,:::\ )}\ar@{-}[ull]^3\ar@{<-}[u]^0\\
{\makebox[0pt][r]{\raisebox{3ex}{\normalsize$\alpha_4$\ }}\young(\bullet\bullet\bullet\bullet,:\ \ \ ,::\ \ ,:::\ )}\ar@{<-}[u]^0&&
{\young(\bullet\bullet\bullet\ ,:\bullet\ \ ,::\ \ ,:::\ )}\ar@{-}[ull]^3\ar@{<-}[u]^1\\
{\young(\bullet\bullet\bullet\ ,:\ \ \ ,::\ \ ,:::\ )}\ar@{-}[u]^3\ar@{<-}[urr]^0&&
{\young(\bullet\bullet\ \ ,:\bullet\ \ ,::\ \ ,:::\ )}\ar@{-}[u]^2\\
&{\young(\bullet\bullet\ \ ,:\ \ \ ,::\ \ ,:::\ )}\ar@{-}[ul]^2\ar@{<-}[ur]^0\\
&{\young(\bullet\ \ \ ,:\ \ \ ,::\ \ ,:::\ )}\ar@{-}[u]^1\\
&{\young(\ \ \ \ ,:\ \ \ ,::\ \ ,:::\ )}\ar@{->}[u]^0\\
}}
\end{center}
\mbox{}\\[-1cm]
$\Aut\bigl(\Hasse(\sfC_2)\bigr)\cong\ZZ/2\ZZ$\\
$\Aut\bigl(\Hasse(\sfC_3)\bigr)\cong\ZZ/2\ZZ\times\ZZ/2\ZZ$\\
$\Aut\bigl(\Hasse(\sfC_l)\bigr)\cong\ZZ/2\ZZ\quad(l\geqslant4)$

\newpage
\setlength{\unitlength}{11mm}
\subsection*{$\mathsf B_l$\quad
$\Bgenl{$\alpha_1$}{$\alpha_2$}{$\alpha_{l-2}$}{$\alpha_{l-1}$}{$\alpha_l$}$}

\mbox{}
\setlength{\unitlength}{3.5mm}
$$\Yboxdim14mm
\young(\Bxwww\Bxxww\Bxxxw\Bxxxx\Bxxxo\Bxxoo\Bxooo,:\Boxww\Boxxw,::\Booxw)$$

\begin{center}
\mbox{}{\tiny\Yboxdim1.5mm
\xymatrix@R=7mm@C=3mm{
{\young(\bullet\bullet\bullet\bullet\bullet\bullet\bullet,:\ \ ,::\ )}&&&&&
{\young(\bullet\bullet\bullet\bullet\ \ \ ,:\bullet\bullet,::\bullet)}\\
{\makebox[0pt][r]{\raisebox{3ex}{\normalsize$\alpha_1$\ }}\young(\bullet\bullet\bullet\bullet\bullet\bullet\ ,:\ \ ,::\ )}\ar@{<-}[u]^0&&
{\makebox[0pt][r]{\raisebox{3ex}{\normalsize$\alpha_2$\ }}\young(\bullet\bullet\bullet\bullet\bullet\ \ ,:\bullet\ ,::\ )}&&
{\makebox[0pt][r]{\raisebox{3ex}{\normalsize$\alpha_3$\ }}\young(\bullet\bullet\bullet\bullet\ \ \ ,:\bullet\bullet,::\ )}\ar@{<-}[ur]^0&&
{\young(\bullet\bullet\bullet\ \ \ \ ,:\bullet\bullet,::\bullet)}\ar@{-}[ul]^4\\
&{\young(\bullet\bullet\bullet\bullet\bullet\ \ ,:\ \ ,::\ )}\ar@{-}[ul]^2\ar@{-}[ur]^1&&
{\young(\bullet\bullet\bullet\bullet\ \ \ ,:\bullet\ ,::\ )}\ar@{-}[ul]^3\ar@{-}[ur]^2&&
{\young(\bullet\bullet\bullet\ \ \ \ ,:\bullet\bullet,::\ )}\ar@{-}[ul]^4\ar@{<-}[ur]^0\\
&&{\young(\bullet\bullet\bullet\bullet\ \ \ ,:\ \ ,::\ )}\ar@{-}[ul]^3\ar@{-}[ur]^1&&
{\young(\bullet\bullet\bullet\ \ \ \ ,:\bullet\ ,::\ )}\ar@{-}[ul]^4\ar@{-}[ur]^2\\
&&{\young(\bullet\bullet\bullet\ \ \ \ ,:\ \ ,::\ )}\ar@{-}[u]^4\ar@{-}[urr]^1&&
{\young(\bullet\bullet\ \ \ \ \ ,:\bullet\ ,::\ )}\ar@{-}[u]^3\\
&&&{\young(\bullet\bullet\ \ \ \ \ ,:\ \ ,::\ )}\ar@{-}[ul]^3\ar@{-}[ur]^1\\
&&&{\young(\bullet\ \ \ \ \ \ ,:\ \ ,::\ )}\ar@{-}[u]^2\\
&&&{\young(\ \ \ \ \ \ \ ,:\ \ ,::\ )}\ar@{->}[u]^0\\
}}
\end{center}
\mbox{}\\
$\Aut\bigl(\Hasse(\sfB_l)\bigr)\cong\ZZ/2\ZZ\quad(l\geqslant2)$\\

\newpage
\setlength{\unitlength}{11mm}
\subsection*{$\mathsf D_l$\quad
$\Dgenl{$\alpha_1$}{$\alpha_2$}{$\alpha_{l-3}$}{$\alpha_{l-2}$}{$\alpha_{l-1}$}{$\alpha_l$}$}
\mbox{}\\[2mm]

\mbox{}
\setlength{\unitlength}{3.5mm}
$$\Yboxdim14mm
\young(\Dxwxx\Dxxxx\Dxxox\Dxxxo\Dxxoo\Dxooo,:\Doxxx\Doxox\Doxxo,::\Dooox\Dooxo)$$

\begin{center}
\mbox{}{\tiny\Yboxdim1.5mm
\xymatrix@R=7mm@C=3mm{&
{\young(\bullet\bullet\bullet\bullet\bullet\bullet,:\ \ \ ,::\ \ )}&&
{\young(\bullet\bullet\bullet\ \ \ ,:\bullet\bullet\ ,::\bullet\ )}&&
{\young(\bullet\bullet\ \bullet\ \ ,:\bullet\ \bullet,::\ \bullet)}\\
{\makebox[0pt][r]{\raisebox{3ex}{\normalsize$\alpha_1$\ }}\young(\bullet\bullet\bullet\bullet\bullet\ ,:\ \ \ ,::\ \ )}\ar@{<-}[ur]^0&&
{\makebox[0pt][r]{\raisebox{3ex}{\normalsize$\alpha_2$\ }}\young(\bullet\bullet\bullet\bullet\ \ ,:\bullet\ \ ,::\ \ )}&&
{\makebox[0pt][r]{\raisebox{3ex}{\normalsize$\alpha_3$\ }}\young(\bullet\bullet\bullet\ \ \ ,:\bullet\bullet\ ,::\ \ )}\ar@{<-}[ul]^0&&
{\makebox[0pt][r]{\raisebox{3ex}{\normalsize$\alpha_4$\ }}\young(\bullet\bullet\ \bullet\ \ ,:\bullet\ \bullet,::\ \ )}\ar@{<-}[ul]^0\\
&{\young(\bullet\bullet\bullet\bullet\ \ ,:\ \ \ ,::\ \ )}\ar@{-}[ul]^2\ar@{-}[ur]^1&&
{\young(\bullet\bullet\bullet\ \ \ ,:\bullet\ \ ,::\ \ )}\ar@{-}[ul]^4\ar@{-}[ur]^(0.6)2&&
{\young(\bullet\bullet\ \bullet\ \ ,:\bullet\ \ ,::\ \ )}
\ar@{-}[ulll]|{\phantom{f}}^(0.4)3\ar@{-}[ur]^2\\
&{\young(\bullet\bullet\bullet\ \ \ ,:\ \ \ ,::\ \ )}\ar@{-}[u]^4\ar@{-}[urr]^(0.6)1&&
{\young(\bullet\bullet\ \bullet\ \ ,:\ \ \ ,::\ \ )}
\ar@{-}[ull]|{\phantom{f}}^(0.4)3\ar@{-}[urr]|{\phantom{f}}^(0.6)1&&
{\young(\bullet\bullet\ \ \ \ ,:\bullet\ \ ,::\ \ )}\ar@{-}[ull]^(0.4)3\ar@{-}[u]^4\\
&&&{\young(\bullet\bullet\ \ \ \ ,:\ \ \ ,::\ \ )}\ar@{-}[ull]^3\ar@{-}[u]^4\ar@{-}[urr]^1\\
&&&{\young(\bullet\ \ \ \ \ ,:\ \ \ ,::\ \ )}\ar@{-}[u]^2\\
&&&{\young(\ \ \ \ \ \ ,:\ \ \ ,::\ \ )}\ar@{->}[u]^0
}}
\end{center}
\mbox{}\\
$\Aut\bigl(\Hasse(\sfD_4)\bigr)\cong\operatorname{Sym}_4$\\
$\Aut\bigl(\Hasse(\sfD_l)\bigr)\cong\Dih_4\quad(l\geqslant5)$

\newpage
\setlength{\unitlength}{11mm}
\subsection*{$\mathsf F_4$\quad
$\Fvier{$\alpha_1$}{$\alpha_2$}{$\alpha_3$}{$\alpha_4$}$}

\setlength{\unitlength}{3.5mm}
$$\Yboxdim14mm
\young(\Fwhfw\Fxhfw\Fxwfw\Fxwww\Fxxww\Foxww,:\Fxwhw\Fxwhx\Fxwwx,:\Fxwwo)$$

\begin{center}
\mbox{}{\tiny\Yboxdim1.5mm
\xymatrix@R=7mm@C=3mm{&
{\young(\bullet\bullet\bullet\bullet\bullet\bullet,:\bullet\bullet\bullet,:\ )}\\
{\young(\bullet\bullet\bullet\bullet\bullet\bullet,:\bullet\bullet\ ,:\ )}\ar@{-}[ur]^3&&
{\makebox[0pt][r]{\raisebox{3ex}{\normalsize$\alpha_2$\ }}\young(\bullet\bullet\bullet\bullet\bullet\ ,:\bullet\bullet\bullet,:\ )}\ar@{<-}[ul]^0&&
{\makebox[0pt][r]{\raisebox{3ex}{\normalsize$\alpha_1$\ }}\young(\bullet\bullet\bullet\bullet\ \ ,:\bullet\bullet\bullet,:\bullet)}\\
{\young(\bullet\bullet\bullet\bullet\bullet\bullet,:\bullet\ \ ,:\ )}\ar@{-}[u]^4&&
{\young(\bullet\bullet\bullet\bullet\bullet\ ,:\bullet\bullet\ ,:\ )}\ar@{<-}[ull]^0\ar@{-}[u]^3&&
{\young(\bullet\bullet\bullet\bullet\ \ ,:\bullet\bullet\bullet,:\ )}\ar@{-}[ull]^1\ar@{-}[u]^2\\
&{\young(\bullet\bullet\bullet\bullet\bullet\ ,:\bullet\ \ ,:\ )}\ar@{<-}[ul]^0\ar@{-}[ur]^4&&
{\young(\bullet\bullet\bullet\bullet\ \ ,:\bullet\bullet\ ,:\ )}\ar@{-}[ul]^1\ar@{-}[ur]^3\\
&{\young(\bullet\bullet\bullet\bullet\ \ ,:\bullet\ \ ,:\ )}\ar@{-}[u]^1\ar@{-}[urr]^4&&
{\young(\bullet\bullet\bullet\ \ \ ,:\bullet\bullet\ ,:\ )}\ar@{-}[u]^2\\
&&{\young(\bullet\bullet\bullet\ \ \ ,:\bullet\ \ ,:\ )}\ar@{-}[ul]^2\ar@{-}[ur]^4\\
&&{\young(\bullet\bullet\bullet\ \ \ ,:\ \ \ ,:\ )}\ar@{-}[u]^3\\
&&{\young(\bullet\bullet\ \ \ \ ,:\ \ \ ,:\ )}\ar@{-}[u]^2\\
&&{\young(\bullet\ \ \ \ \ ,:\ \ \ ,:\ )}\ar@{-}[u]^1\\
&&{\young(\ \ \ \ \ \ ,:\ \ \ ,:\ )}\ar@{->}[u]^0
}}
\end{center}
\mbox{}\\
$\Aut\bigl(\Hasse(\sfF_4)\bigr)=1$\\

\newpage
\section{$\sfE_6$}
The picture shows the Hasse graph of the poset of abelian ideals
for type~$\sfE_6$. I chose to draw it in a way in which the
symmetry becomes manifest.
The nodes marked by $\varphi=\theta,\alpha_1,\dots,\alpha_6$ carry
the abelian ideals $\ab^{\varphi,\min}$. The encircled nodes mark
the maximal abelian ideals $\ab^{\alpha_i,\max}$ ($i={1,\dots,6}$).

\begin{center}
\setlength{\unitlength}{0.00063333in}
\begin{picture}(6816,7031)(0,-10)
\path(3933,3033)(2433,3333)(2433,5133)
	(3933,4833)(5133,5733)(3633,6033)
	(3633,4233)(5133,3933)(3933,3033)
\path(4533,3483)(3033,3783)(3033,5583)
	(4533,5283)(4533,3483)
\path(3933,3933)(2433,4233)(3633,5133)
	(5133,4833)(3933,3933)
\path(3183,3183)(3183,4983)(4383,5883)
	(4383,4083)(3183,3183)
\path(3933,4833)(3933,3033)(3783,2733)
	(3033,2883)(3183,3183)
\path(4533,3483)(4383,3183)(3633,3333)(3783,3633)
\path(3633,3333)(2433,2433)(3183,2283)
	(4383,3183)(5133,3033)(3933,2133)(3183,2283)
\path(4533,2583)(3783,2733)
\path(3933,2133)(3783,1833)(3783,33)
\path(3483,3933)(3633,4233)(2433,3333)
\path(5133,3933)(5133,5733)(4983,5433)
	(4233,5583)(4383,5883)
\path(5133,4833)(4983,4533)(4233,4683)(4383,4983)
\path(4233,4683)(4233,6483)(5733,6183)
	(5733,4383)(4983,4533)
\path(4983,4533)(4983,6333)
\path(4983,5433)(5733,5283)
\path(4533,4383)(3033,4683)(2883,4383)
	(2883,5283)(3033,5583)
\path(2433,4233)(2283,3933)(2283,4833)
	(2433,5133)(3633,6033)
\path(2883,4383)(1683,3483)(1683,5283)
	(2883,6183)(2883,5283)
\path(2883,5283)(1683,4383)
\path(2283,4833)(2283,5733)
\path(1683,5283)(1533,4983)(33,5283)
\path(4383,4983)(3183,4083)
\path(3783,3633)(3783,5433)
\path(5733,6183)(5583,5883)(6783,6783)
\put(3683,483){\makebox(0,0)[lb]{\tiny$0$}}
\put(3683,1383){\makebox(0,0)[lb]{\tiny$1$}}
\put(3783,1983){\makebox(0,0)[lb]{\tiny$2$}}
\put(4308,2283){\makebox(0,0)[lb]{\tiny$3$}}
\put(4908,2733){\makebox(0,0)[lb]{\tiny$5$}}
\put(4908,3633){\makebox(0,0)[lb]{\tiny$3$}}
\put(5208,4308){\makebox(0,0)[lb]{\tiny$1$}}
\put(5808,4758){\makebox(0,0)[lb]{\tiny$0$}}
\put(5808,5658){\makebox(0,0)[lb]{\tiny$1$}}
\put(5583,6033){\makebox(0,0)[lb]{\tiny$2$}}
\put(5958,6033){\makebox(0,0)[lb]{\tiny$3$}}
\put(6558,6483){\makebox(0,0)[lb]{\tiny$5$}}
\put(5358,6333){\makebox(0,0)[lb]{\tiny$4$}}
\put(4683,6483){\makebox(0,0)[lb]{\tiny$6$}}
\put(4008,6033){\makebox(0,0)[lb]{\tiny$4$}}
\put(3258,5883){\makebox(0,0)[lb]{\tiny$3$}}
\put(2508,6033){\makebox(0,0)[lb]{\tiny$5$}}
\put(1908,5583){\makebox(0,0)[lb]{\tiny$3$}}
\put(1533,5133){\makebox(0,0)[lb]{\tiny$2$}}
\put(1083,4908){\makebox(0,0)[lb]{\tiny$4$}}
\put(408,5058){\makebox(0,0)[lb]{\tiny$6$}}
\put(1583,4758){\makebox(0,0)[lb]{\tiny$1$}}
\put(1583,3933){\makebox(0,0)[lb]{\tiny$0$}}
\put(2333,3708){\makebox(0,0)[lb]{\tiny$1$}}
\put(2808,3108){\makebox(0,0)[lb]{\tiny$4$}}
\put(2808,2208){\makebox(0,0)[lb]{\tiny$6$}}
\put(3408,2058){\makebox(0,0)[lb]{\tiny$4$}}
\put(3783,2883){\makebox(0,0)[lb]{\tiny$2$}}
\put(5058,5483){\makebox(0,0)[lb]{\tiny$2$}}
\put(2283,4983){\makebox(0,0)[lb]{\tiny$2$}}
\put(3583,4008){\makebox(0,0)[lb]{\tiny$2$}}
\put(3858,858){\makebox(0,0)[lb]{\small$\theta$}}
\put(3483,3733){\makebox(0,0)[b]{\small$\alpha_1$}}
\put(3658,4908){\makebox(0,0)[b]{\small$\alpha_2$}}
\put(2903,4183){\makebox(0,0)[b]{\small$\alpha_3$}}
\put(4233,4483){\makebox(0,0)[b]{\small$\alpha_4$}}
\put(1683,3283){\makebox(0,0)[b]{\small$\alpha_5$}}
\put(5733,4183){\makebox(0,0)[b]{\small$\alpha_6$}}
\put(6783,6783){\circle*{50}}\put(6783,6783){\circle{90}}
\put(6183,6333){\circle*{50}}
\put(5583,5883){\circle*{50}}
\put(5733,6183){\circle*{50}}
\put(4983,6333){\circle*{50}}
\put(4233,6483){\circle*{50}}\put(4233,6483){\circle{90}}
\put(5733,5283){\circle*{50}}
\put(5733,4383){\circle*{50}}
\put(4983,5433){\circle*{50}}
\put(4983,4533){\circle*{50}}
\put(4233,5583){\circle*{50}}
\put(4233,4683){\circle*{50}}
\put(5133,4833){\circle*{50}}
\put(5133,5733){\circle*{50}}
\put(4383,4983){\circle*{50}}
\put(4383,5883){\circle*{50}}
\put(5133,3933){\circle*{50}}
\put(4383,4083){\circle*{50}}
\put(3633,6033){\circle*{50}}\put(3633,6033){\circle{90}}
\put(3633,5133){\circle*{50}}
\put(3633,4233){\circle*{50}}
\put(4533,3483){\circle*{50}}
\put(4533,4383){\circle*{50}}
\put(4533,5283){\circle*{50}}
\put(3783,5433){\circle*{50}}
\put(3783,4533){\circle*{50}}
\put(3783,3633){\circle*{50}}
\put(3033,3783){\circle*{50}}
\put(3033,4683){\circle*{50}}
\put(3033,5583){\circle*{50}}
\put(3483,3933){\circle*{50}}\put(3483,3933){\circle{90}}
\put(3933,3033){\circle*{50}}
\put(3933,3933){\circle*{50}}
\put(3933,4833){\circle*{50}}
\put(3183,4983){\circle*{50}}
\put(3183,4083){\circle*{50}}
\put(3183,3183){\circle*{50}}
\put(2433,3333){\circle*{50}}
\put(2433,4233){\circle*{50}}
\put(2433,5133){\circle*{50}}
\put(4383,3183){\circle*{50}}
\put(3633,3333){\circle*{50}}
\put(5133,3033){\circle*{50}}
\put(4533,2583){\circle*{50}}
\put(3783,2733){\circle*{50}}
\put(3033,2883){\circle*{50}}
\put(2433,2433){\circle*{50}}
\put(3183,2283){\circle*{50}}
\put(3933,2133){\circle*{50}}
\put(3783,1833){\circle*{50}}
\put(3783,933){\circle*{50}}
\put(2883,4383){\circle*{50}}
\put(2883,5283){\circle*{50}}
\put(2883,6183){\circle*{50}}\put(2883,6183){\circle{90}}
\put(2283,5733){\circle*{50}}
\put(2283,4833){\circle*{50}}
\put(2283,3933){\circle*{50}}
\put(1683,3483){\circle*{50}}
\put(1683,4383){\circle*{50}}
\put(1683,5283){\circle*{50}}
\put(1533,4983){\circle*{50}}
\put(783,5133){\circle*{50}}
\put(33,5283){\circle*{50}}\put(33,5283){\circle{90}}
\put(3783,33){\circle*{50}}
\put(3408,-400){\makebox(0,0){
$\Aut\bigl(\Hasse(\sfE_6)\bigr)\cong\operatorname{Sym}_3$}}
\end{picture}
\setlength{\unitlength}{8mm}
\begin{picture}(4,4.1)(0,-.3)
\multiput(0,1)(1,0){5}{\circle*{0.2}}\multiput(2,0)(0,-1){2}{\circle*{0.2}}
\put(0,1){\line(1,0){4}}\put(2,-1){\line(0,1){2}}
\put(0,1.333){\makebox(0,0)[b]{\footnotesize$5$}}
\put(1.667,0){\makebox(0,0)[r]{\footnotesize$1$}}
\put(1,1.333){\makebox(0,0)[b]{\footnotesize$3$}}
\put(2,1.333){\makebox(0,0)[b]{\footnotesize$2$}}
\put(3,1.333){\makebox(0,0)[b]{\footnotesize$4$}}
\put(4,1.333){\makebox(0,0)[b]{\footnotesize$6$}}
\put(1.667,-1){\makebox(0,0)[r]{\footnotesize$0$}}
\end{picture}
\end{center}

\end{document}